\documentclass[a4paper,oneside]{article}

\title{Killing (super)algebras for generalised spin manifolds}
\author{Andrew D.K. Beckett}
\date{July 2026}

\usepackage[utf8]{inputenc} % allows use of characters like é,ü,@,£ etc. in input
\usepackage[hmargin=1in]{geometry}
\usepackage[british]{babel} % british english typographical rules

\usepackage[
	backend = biber, % default is biber, which is much slower - change to biber for final version
	style = numeric-comp, % for better numerical citations, replaces cite pkg
	giveninits = true,
	sorting = none, % sort in order of citation
	maxbibnames = 5, % prevent ``et al'' on most items
	url = false, % suppress URLs since most items have DOI or arXiv links; doesn't affect "online" type items
	backref = true % print pages where item is cited in bibliography
	]{biblatex}
\AtEveryBibitem{\ifentrytype{article}{\clearfield{issn}}{}}
\AtEveryBibitem{\ifentrytype{article}{\clearfield{editor}}{}}

\usepackage{csquotes}  % recommended with babel, stops biblatex throwing errors
\usepackage{mathtools} % includes amsmath plus extras, such as multlined
\usepackage{amsthm,amscd,amssymb} % extra AMS packages, may not be needed
\usepackage{stmaryrd} % improved math symbols and some new ones
\usepackage{physics} % physics symbols, includes amsmath
\usepackage{tensor} % index typesetting
\usepackage[makeroom]{cancel} % allows \cancel and related commands in math mode
\usepackage{faktor} % better display math quotients
\usepackage{mathrsfs} % mathscr font
\usepackage{bbold} % more symbols in mathbb
\usepackage{tikz} % for creating graphics, diagrams etc.
\usetikzlibrary{cd}%calc,arrows % commutative diagrams
\usepackage{enumitem} % layout options for lists
\setenumerate[1]{label={(\arabic*)}}
\usepackage[titletoc,title]{appendix}  % additional options for appendices
\usepackage[unicode]{hyperref}
\usepackage{xurl} % this and breaklinks below requried to prevent links over-running in arXiv version
\hypersetup{%
  pdftitle   = {ksa-gen-spin},
  pdfkeywords = {Lie superalgebra, Killing superalgebra, Killing spinor, connection, superconnection, isometry, Spencer cohomology, filtration, deformation, homogeneous, R-symmetry, gauged supergravity, spin-G, generalised spin},
  pdfauthor  = {Andrew DK Beckett}
  pdfcreator = {\LaTeX\ with package \flqq hyperref\frqq},
  colorlinks = true,
  linkcolor = blue,
  menucolor = blue,
  urlcolor = [rgb]{0.7,0,0},
  citecolor = [rgb]{0,0.7,0},
  breaklinks = true
}

\theoremstyle{plain}
\newtheorem{lemma}{Lemma}
\newtheorem{proposition}[lemma]{Proposition}
\newtheorem{theorem}[lemma]{Theorem}
\newtheorem{corollary}[lemma]{Corollary}
\newtheorem{definition}[lemma]{Definition}

\theoremstyle{definition}
\newtheorem{remark}[lemma]{Remark}
\newtheorem{example}[lemma]{Example}
\newcommand{\fa}{\mathfrak{a}}

\newcommand{\fg}{\mathfrak{g}}
\newcommand{\fh}{\mathfrak{h}}
\newcommand{\fk}{\mathfrak{k}}

\newcommand{\fs}{\mathfrak{s}}
\newcommand{\fgl}{\mathfrak{gl}}

\newcommand{\fso}{\mathfrak{so}}
\newcommand{\fiso}{\mathfrak{iso}}

\newcommand{\fX}{\mathfrak{X}}
\newcommand{\fK}{\mathfrak{K}}
\newcommand{\fr}{\mathfrak{r}}
\newcommand{\fR}{\mathfrak{R}}
\newcommand{\fS}{\mathfrak{S}}
\newcommand{\fV}{\mathfrak{V}}
\newcommand{\fD}{\mathfrak{D}}
\newcommand{\fann}{\mathfrak{ann}}
\newcommand{\fstab}{\mathfrak{stab}}
\newcommand{\fsymm}{\mathfrak{symm}}
\newcommand{\RR}{\mathbb{R}}
\newcommand{\CC}{\mathbb{C}}
\newcommand{\HH}{\mathbb{H}}
\newcommand{\ZZ}{\mathbb{Z}}

\newcommand{\eA}{\mathscr{A}}

\newcommand{\eC}{\mathscr{C}}
\newcommand{\eD}{\mathscr{D}}
\newcommand{\eE}{\mathscr{E}}

\newcommand{\eL}{\mathscr{L}}

\newcommand{\cH}{\mathcal{H}}
\newcommand{\cK}{\mathcal{K}}

\newcommand{\cHhat}{\widehat{\cH}}
\newcommand{\cKhat}{\widehat{\cK}}
\newcommand{\ssB}{\mathsf{B}}
\newcommand{\ssC}{\mathsf{C}}
\newcommand{\ssH}{\mathsf{H}}
\newcommand{\ssZ}{\mathsf{Z}}

\newcommand{\Sbundle}{\underline{S}}
\newcommand{\Shatbundle}{\underline{\widehat{S}}}

\newcommand{\Wbundle}{\underline{W}}

\newcommand{\Qbar}{\overline{Q}}

\newcommand{\Rbar}{\overline{R}}
\newcommand{\fshat}{\widehat{\fs}}

\newcommand{\nablahat}{\widehat{\nabla}}
\newcommand{\eLhat}{\widehat{\eL}}
\newcommand{\dhat}{\widehat{\dd}}
\newcommand{\fKhat}{\widehat{\fK}}
\newcommand{\fVhat}{\widehat{\fV}}
\newcommand{\fShat}{\widehat{\fS}}
\newcommand{\Rhat}{\widehat{R}}
\newcommand{\Phat}{\widehat{P}}

\newcommand{\gammahat}{\widehat{\gamma}}

\newcommand{\eDhat}{\widehat{\eD}}
\newcommand{\eEhat}{\widehat{\eE}}
\newcommand{\Phihat}{\widehat{\Phi}}
\newcommand{\Psihat}{\widehat{\Psi}}
\newcommand{\varphihat}{\widehat{\varphi}}
\newcommand{\sigmahat}{\widehat{\sigma}}
\newcommand{\fhhat}{\widehat{\fh}}
\newcommand{\Hhat}{\widehat{H}}
\DeclareFontFamily{U}{mathx}{}
\DeclareFontShape{U}{mathx}{m}{n}{<-> mathx10}{}
\DeclareSymbolFont{mathx}{U}{mathx}{m}{n}
\DeclareMathAccent{\widecheck}{0}{mathx}{"71}
\newcommand{\Sch}{\widecheck{S}}
\newcommand{\fsch}{\widecheck{\fs}}
\newcommand{\frch}{\widecheck{\fr}}
\newcommand{\Rch}{\widecheck{R}}
\newcommand{\fhch}{\widecheck{\fh}}
\newcommand{\Htilde}{\widetilde{H}}

\newcommand{\betatilde}{\widetilde{\beta}}
\newcommand{\gammatilde}{\widetilde{\gamma}}
\newcommand{\deltatilde}{\widetilde{\delta}}
\newcommand{\rhotilde}{\widetilde{\rho}}
\newcommand{\thetatilde}{\widetilde{\theta}}
\newcommand{\fatilde}{\widetilde{\fa}}

\newcommand{\varphitilde}{\widetilde{\varphi}}

\newcommand{\kappabold}{\boldsymbol{\kappa}}
\newcommand{\betabold}{\boldsymbol{\beta}}
\newcommand{\gammabold}{\boldsymbol{\gamma}}
\newcommand{\rhobold}{\boldsymbol{\rho}}
\newcommand{\pihat}{\widehat{\pi}}

\newcommand{\varpihat}{\widehat{\varpi}}

\newcommand{\Wedge}{\mathchoice{{\textstyle\bigwedge}}{\bigwedge}{\bigwedge}{\bigwedge}}	% display version is {\texstyle\bigwedge}; text, script and scriptscript versions are \bigwedge
\newcommand{\Odot}{\mathchoice{{\textstyle\bigodot}}{\bigodot}{\bigodot}{\bigodot}}
\newcommand{\Otimes}{\mathchoice{{\textstyle\bigotimes}}{\bigotimes}{\bigotimes}{\bigotimes}}

\newcommand{\1}{\mathbb{1}}
\newcommand{\ccomm}[2]{\qty[\comm{#1}{#2}]}
\newcommand{\pair}[2]{\left\langle #1,#2\right\rangle} %comma inner product
\DeclareMathOperator{\Id}{Id}

\DeclareMathOperator{\End}{End}
\DeclareMathOperator{\Hom}{Hom}
\DeclareMathOperator{\GL}{GL}
\DeclareMathOperator{\Spin}{Spin}

\DeclareMathOperator{\Orth}{O}
\DeclareMathOperator{\SO}{SO}

\DeclareMathOperator{\Sp}{Sp}
\DeclareMathOperator{\USp}{USp}
\DeclareMathOperator{\U}{U}
\DeclareMathOperator{\SU}{SU}
\DeclareMathOperator{\Lie}{Lie}

\DeclareMathOperator{\Stab}{Stab}

\DeclareMathOperator{\Ad}{Ad}
\DeclareMathOperator{\ad}{ad}

\DeclareMathOperator{\Gr}{Gr}

\DeclareMathOperator{\diag}{diag}
\DeclareMathOperator{\evaluate}{ev}
\DeclareMathOperator{\pr}{pr}

\begin{document}

\maketitle

\begin{center}
	\textit{In memory of Paul de Medeiros.}
\end{center}

\begin{abstract}
	\noindent
	We define the notion of a Killing (super)algebra for a connection on a spinor bundle associated to a generalised spin structure on a pseudo-Riemannian manifold of any signature. We are led naturally to include in the even subspace not only Killing vectors but also certain infinitesimal gauge transformations, and we show that the definition of the (super)algebra requires, in addition to the spinor connection and a Dirac current, a map to pair spinor fields into infinitesimal gauge transformations. We show that these (super)algebras are filtered subdeformations of (an analogue of) the Poincaré superalgebra extended by the \(R\)-symmetry algebra. By employing Spencer cohomology, we study such deformations from a purely algebraic point of view and, at least in the case of Lorentzian signature and high supersymmetry, identify the subclass of deformations to which the Killing superalgebras belong. Finally, we show that, with some caveats, one can reconstruct a supersymmetric background geometry from such a deformation as a homogeneous space on which the deformation is realised as a subalgebra of the Killing superalgebra.
\end{abstract}

\tableofcontents

\section{Introduction}
\label{sec:introduction}

A Killing superalgebra is a Lie superalgebra with an even subspace consisting of Killing vector fields on some pseudo-Riemannian manifold \((M,g)\) and odd subspace consisting of ``Killing spinors''; that is to say sections of a bundle of spinors on \((M,g)\) which are parallel with respect to a connection \(D\), sometimes with additional algebraic constraints. These are not, in general, the geometric Killing spinors encountered in the Riemannian geometry literature, but belong to a more general class termed \emph{differential spinors} \cite{Shahbazi2024_1,Shahbazi2024_2} also studied in general form in e.g. \cite{Lazaroiu2016,Cortes2021}. These notions of Killing spinors and superalgebras arise in supergravity theory, where they describe the supersymmetries of certain solutions, but in the study of the algebraic structure of the superalgebras via Spencer cohomology (which began with \cite{Figueroa-OFarrill2017,Figueroa-OFarrill2017_1}), it has been shown that one can write down connections on spinor bundles which are more general than those known to arise from supergravity but which nonetheless allow one to define supersymmetric backgrounds and their Killing superalgebras \cite{Beckett2021,deMedeiros2018}. Moreover, earlier work inspired by the Killing superalgebras of supergravity showed that the \emph{geometric} Killing spinors and Killing vectors of higher spheres can be arranged in an analogous way into Lie \emph{algebras}, providing a geometric realisation of some exceptional algebras as well as \(\fso(8)\)-triality \cite{Figueroa-OFarrill2008_except}.

In recent work \cite{Beckett2024_ksa}, we developed a framework for understanding Killing (super)algebras which unifies and systematises these examples, which had, until then, generally been studied on a rather ad-hoc basis, with the exception of work on the Homogeneity Theorem for highly supersymmetric solutions of supergravity by Hustler \cite{Hustler2016}. Remaining agnostic as to signature and other choices, we reported a set of necessary and sufficient conditions for the existence of a Killing (super)algebra associated to a connection on a spinor bundle, and, generalising results in the case of 11-dimensional supergravity \cite{Figueroa-OFarrill2017_1}, described the algebraic structure of those (super)algebras -- they are \emph{filtered subdeformations} of the Poincaré superalgebra (or an appropriate analogue thereof, depending upon signature). We analysed this structure in more detail through the lens of Spencer cohomology, finding that in Lorentzian signature, the Homogeneity Theorem provided significant simplifications in the highly supersymmetric case. We further exploited these simplifications to describe the reconstruction of highly supersymmetric backgrounds as homogeneous spaces from their Killing superalgebras.

In the present work, our goal is to generalise the results of \cite{Beckett2024_ksa} to the case of Killing superalgebras associated to manifolds equipped with a \emph{generalised spin structure} (also known as a \emph{spin-\(G\) structure}); that is, an analogue of a spin structure where the spin group is ``twisted'' with some ``internal'' symmetry group \(G\) with a central \(\ZZ_2\) subgroup, allowing topological obstructions to the existence of ordinary spin structures to be overcome. More precisely, instead of the structure group of the special orthonormal frame bundle of \((M,g)\) being lifted from \(\SO(s,t)\) to \(\Spin(s,t)\), it is lifted to \(\Spin^G(s,t):=\Spin(s,t)\times_{\ZZ_2} G\). The physical relevance of such structures was first mooted by Hawking and Pope \cite{Hawking1978}, then by Back, Freund and Forger \cite{Back1978}, who noted that they provide a way for fermions to be coupled to spacetime backgrounds which are not spin at the expense of those fermions carrying charge quantum numbers for the group \(G\) which are algebraically related to their spin quantum numbers. A more mathematical perspective was then provided by Forger, Hess, Avis and Isham \cite{Forger1979,Avis1980}. The most well-known example, with \(G=\U(1)\), is the spin-\(c\) structure introduced by Atiyah, Bott and Shapiro \cite{Atiyah1964} (and in Lorentzian signature by Whiston \cite{Whiston1975}; see also \cite{Atiyah1973,Lawson1989}), pre-dating Hawking--Pope. The next is spin-\(h\) (or spin-\(q\)), where \(G=\SU(2)\cong\Sp(1)\), which has gained the attention of geometers \cite{Nagase1995,Lawson2023,Albanese2021} although was already considered (not under this name) in \cite{Back1978,Avis1980}. A natural generalisation which has received significant recent interest from geometers is the case \(G=\Spin(k)\) for \(k\in\mathbb{N}\) (giving spin, spin-\(c\) and spin-\(h\) for \(k=1,2,3\) respectively), known as \emph{spinorially twisted spin} or \emph{(generalised) spin-\(k\)} structures \cite{Espinosa2016,Albanese2021,Artacho2023}. Coming from a different perspective, the closely-related \emph{real Lipschitz structures}, more precisely in the guise of \emph{canonical spin structures}, have been used to provide necessary and sufficient conditions for the existence of bundles of spinors in any signature \cite{Lazaroiu2019_1}; a summary with a focus on application to supergravity is given in \cite{Lazaroiu2019}. Other recent work on physics applications includes \cite{Janssens2018}. See \cite{Beckett2025_gen_spin} and the references therein for more on generalised spin structures.

It will be necessary for our purposes for \(G\) to be an \emph{\(R\)-symmetry group}; that is, the group of automorphisms of a spinor representation which leaves the Dirac current -- the map which allows spinors to be paired into vectors -- invariant.
The \(R\)-symmetry group also acts via automorphisms on the Poincaré superalgebra (and its generalisations to other signatures), with the action being effective on the odd part of the superalgebra -- it ``rotates the supercharges'' in the physics parlance -- and trivial on the even part; the superalgebra thus admits an extension which includes the \(R\)-symmetry Lie algebra. 
Such extended Poincaré superalgebras first appeared in the seminal work of Haag--Łopuszański--Sohnius \cite{Haag1975}, where it was shown that, under a standard set of assumptions, the supersymmetry algebra of a 4-dimensional Poincaré-invariant quantum field theory must be the direct sum of the Poincaré superalgebra and a Lie algebra of ``internal'' symmetries, but that the Poincaré component could be \(N\)-extended (i.e. it could contain \(N\) minimal supercharge multiplets) and extended by the \(R\)-symmetry algebra -- ``central charge'' terms in odd-odd brackets are also admitted but will not play a role for us. The notion was generalised to arbitrary dimension by Strathdee \cite{Strathdee1986}, and a brief summary appears in the review of Van Proeyen \cite{VanProeyen1999}.
A complete treatment in arbitrary signature and with various choices of Dirac current did not appear in the literature until fairly recently \cite{Gall2021}, but see also \cite{DAuria2001_2}.

Killing superalgebras ``twisted'' by \(R\)-symmetry were first introduced by de Medeiros, Figueroa-O'Farrill and Santi in the 6-dimensional Lorentzian case by considering filtered subdeformations of the \(R\)-symmetry-extended Poincaré superalgebra \cite{deMedeiros2018}.
While we take a similar tack in Section~\ref{sec:filtered-def-flat-model-rsymm}, we first provide a more geometrical motivation in Section~\ref{sec:killing-rsymm}, arguing that such twisted superalgebras are more natural than the ordinary type studied in \cite{Beckett2024_ksa} when one works with generalised spin manifolds.
In the work of de Medeiros \textit{et al}., two different treatments of the twisted Killing superalgebra were given. The first \cite[\S6.3]{deMedeiros2018} takes place in a more geometric formalism but somewhat obscures the gauge-theoretic nature of the construction by assuming a flat \(R\)-symmetry connection; the resulting definition of the Killing superalgebra is clear, but not clearly motivated. The other approach \cite[Sec.7]{deMedeiros2018} is an entirely local one which makes the gauge theory very explicit and builds a definition of the Killing superalgebra which is much more clearly motivated but essentially elides the question of global structure. In addition to placing these superalgebras in their proper geometric setting, we will synthesise the two approaches, taking the emphasis on geometric formality of the former and motivating the definition in a similar way to the latter.
As an added benefit, out approach admits \(R\)-symmetry connections which are not flat.

We note that the Killing superalgebras discussed here are expected to arise as the symmetry superalgebras of bosonic backgrounds of supergravity theories in which the \(R\)-symmetry algebra is gauged (such that the gravitini acquire colour charge) -- see for example \cite{Ceresole2000,Gunaydin2001} for the minimal and matter-coupled 5-dimensional theories and \cite{Bellorin2007,Bellorin2009} for a discussion of the Killing spinors in those theories. This point is argued in the author's PhD thesis \cite[Sec.7.2, Sec.8.2]{Beckett2024_thesis} and in related upcoming work.

As in \cite{Beckett2024_ksa}, we take a very general approach here, defining Killing (super)algebras on generalised spin manifolds independently of various choices such as dimension, signature, extension of the spinor module and the Dirac current, aiming to generalise previous work \cite{Figueroa-OFarrill2017,deMedeiros2018} as far as possible, but we will eventually specialise to Lorentzian signature and a restricted choice of Dirac currents, motivated by the Homogeneity Theorem.

This work is organised as follows. Section~\ref{sec:background} provides background on ``flat model'' (super)algebras -- generalisations of the Poincaré superalgebra to general signature and choice of Dirac current as discussed in \cite{Alekseevsky1997,Beckett2024_ksa} -- as well as their \(R\)-symmetry groups and extensions by the \(R\)-symmetry Lie algebra. Also discussed are spin-\(R\) structures -- generalised spin structures twisted by the \(R\)-symmetry group. In Section~\ref{sec:killing-rsymm} we construct a definition of the Killing (super)algebra for such a structure, finding necessary and sufficient conditions on the defining data for it to be a Lie (super)algebra (Theorem~\ref{thm:killing-algebra-exist-rsymm}) and showing that it is a filtered subdeformation of the flat model (Theorem~\ref{thm:killing-algebra-filtered-rsymm}). In Section~\ref{sec:filtered-def-flat-model-rsymm}, we consider the purely algebraic perspective, discussing the filtered subdeformations abstractly using Spencer cohomology and identifying, at least in the highly supersymmetric Lorentzian case, a subclass which may be realised as Killing superalgebras (Theorem~\ref{thm:admiss-cocycle-integr-rsymm}). In Section~\ref{sec:high-susy-spin-mfld-rsymm}, we then consider the problem of reconstructing a geometric background and its Killing superalgebra as a homogeneous space of such a deformation, culminating in Theorem~\ref{thm:reconstruction-homog-spin-mfld-rsymm}.

Much of the background material on spin-\(R\) structures contained herein (\S\ref{sec:spin-R}, \S\ref{sec:alg-structure-killing-rsymm}, \S\ref{sec:homogeneous-lorentz-spin-mfld-rsymm}) applies equally well to all generalised spin structures and is expanded upon in the more general context and in more detail in a separate set of notes \cite{Beckett2025_gen_spin}. To the author's knowledge, the covariant Cartan calculus and symmetry algebra of a generalised spin structure (in particular Proposition~\ref{prop:symm-alg-spin-r}) developed here and in \textit{loc. cit.} do not appear elsewhere in the literature, but the latter is at least related to the symmetry algebra of a regular Cartan geometry discussed by Čap and Neusser \cite{Cap2009}.
Our major results,  Theorems~\ref{thm:killing-algebra-exist-rsymm}, \ref{thm:killing-algebra-filtered-rsymm}, \ref{thm:admiss-cocycle-integr-rsymm} and \ref{thm:reconstruction-homog-spin-mfld-rsymm}, are all new, but generalise various results in \cite{Beckett2024_ksa,Figueroa-OFarrill2017_1,deMedeiros2018}, with new technicalities having to be overcome. To aid comparison, this paper is structured in a similar way to \cite{Beckett2024_ksa}, which it most directly generalises, and we use the same notation and conventions unless otherwise indicated.

\newpage

\section{Background}
\label{sec:background}

We will make use of the terminology and notation presented in Section~2 and Appendix~A of \cite{Beckett2024_ksa}.

\subsection{Flat model (super)algebras and \(R\)-symmetry}
\label{sec:flat-models-R-symm}

In this section, we recall the definitions of Dirac currents and (generalised) Poincaré superalgebras, define their \(R\)-symmetry groups and discuss ``extending'' them by their \(R\)-symmetry algebras.

\subsubsection{Ordinary flat model (super)algebras}
\label{sec:flat-model-algs}

In \cite[\S2.1]{Beckett2024_ksa}, we gave the following definition for a generalisation of the Poincaré superalgebra which, guided by links to Cartan geometry, we named \emph{flat model (super)algebras}, though other sources use various other names (including \emph{\(N\)-(super)-extended Poincaré algebras}, \emph{\(\pm N\)-extended Poincaré algebras} \cite{Alekseevsky1997} and \emph{Poincaré Lie superalgebras} \cite{Cortes2020,Gall2021}). We use conventions for spinor modules explained in \cite[\S{}A.1]{Beckett2024_ksa}.

\begin{definition}[Dirac current, flat model (super)algebra {\cite[Def.2.1]{Beckett2024_ksa}}]\label{def:flat-model-alg}
	Let \((V,\eta)\) be a (pseudo-)inner product space and let \(S\)  be a (possibly \(N\)-extended) real spinor module. 
	We define the \(\ZZ\)-graded vector space
	\begin{align}\label{eq:grading-fs}
		&\fs=\fs_{-2}\oplus\fs_{-1}\oplus\fs_{0}
		&&\text{where}
		&\fs_{-2} &= V,
		&\fs_{-1} &= S,
		&\fs_0 &= \fso(V),
		&\fs_i &= 0 \text{ otherwise}.
	\end{align}
	An \(\fso(V)\)-module map \(\kappa:\Odot^2 S\to V\) (resp. \(\kappa:\Wedge^2 S\to V\)) is called a \emph{symmetric} (resp. \emph{skew-symmetric}) \emph{Dirac current}. Given such a map, we define the \emph{flat model superalgebra} (resp. \emph{flat model algebra}) to be the \(\ZZ\)-graded Lie superalgebra (resp. Lie algebra) \((\fs,\comm{-}{-})\) where the bracket is defined as follows:
	\begin{equation}
	\begin{aligned}\label{eq:flat-model-bracket}
	    \comm{A}{B} &= AB - BA,
	    & \comm{A}{v} &= Av,
	    & \comm{v}{w} &= 0, \\
	    \comm{\epsilon}{\zeta} &= \kappa(\epsilon,\zeta),
	    &\comm{A}{\epsilon} &= A\cdot\epsilon,
	    &\comm{v}{\epsilon} &= 0,
	\end{aligned}
	\end{equation}
	where \(A,B\in\fso(V)\), \(v,w\in V\), \(\epsilon,\zeta\in S\).
\end{definition}

\subsubsection{Schur algebra and \(R\)-symmetry group}
\label{sec:R-symm}

Let \((V,\eta)\) be an inner product space, let \(S\) be a (possibly \(N\)-extended) spinor representation of \(\Spin(V)\) and define the \emph{Schur algebra} \(\eC(S)\) as the algebra of \(\fso(V)\)-equivariant endomorphisms of \(S\) and the \emph{Schur group} \(\eC^\times(S)\) as its group of units \cite{Alekseevsky1997,Gall2021}:
\begin{align}
	\eC(S) &= \End_{\fso(V)}(S) = \qty{ a\in\End(S) \mid \forall A\in\fso(V),\forall\epsilon\in S,\,a(A\cdot\epsilon) = A\cdot(a\epsilon) },	
	\label{eq:schur-alg}\\
	\eC^\times(S) &= \GL_{\fso(V)}(S) = \qty{ a\in\GL(S) \mid \forall A\in\fso(V),\forall\epsilon\in S,\,a(A\cdot\epsilon) = A\cdot(a\epsilon) }.
	\label{eq:schur-group}
\end{align}
Now suppose we have an \(\fso(V)\)-equivariant map \(\kappa:\Odot^2S\to V\) or \(\kappa:\Wedge^2S\to V\) (i.e. a Dirac current). We define the \emph{\(R\)-symmetry group of \((S,\kappa)\)} as the subgroup of the Schur group preserving \(\kappa\):
\begin{equation}\label{eq:r-symm-group}
	R_\kappa = \{ a\in\eC^\times(S) \mid \forall\epsilon,\epsilon'\in S,\,\kappa(a\epsilon,a\epsilon') = \kappa(\epsilon,\epsilon')\};
\end{equation}
its Lie algebra, the \emph{\(\fr\)-symmetry algebra of \((S,\kappa)\)}, is of course
\begin{equation}\label{eq:r-symm-alg}
	\fr_\kappa = \{ a\in\eC(S) \mid \forall\epsilon,\epsilon'\in S,\,\kappa(a\epsilon,\epsilon') + \kappa(\epsilon,a\epsilon') = 0\}.
\end{equation} 
We will omit the subscript \(\kappa\) when the choice of Dirac current is unambiguous. We note that some sources give slightly different definitions for the \(R\)-symmetry group; our definition is the same as in \cite{Gall2021}, allowing for the \(R\)-symmetry group to be disconnected.

Now let \(\fs\) be the flat model (super)algebra defined by \((V,\eta,S,\kappa)\). Since the action of \(R_\kappa\) preserves \(\kappa\), it acts by \(\ZZ\)-graded Lie (super)algebra automorphisms on \(\fs\) where the action on \(\fs_{\overline 0}=V\oplus\fso(V)\) is trivial. Likewise, \(\fr_\kappa\) acts by \(\ZZ\)-graded derivations. In fact, at least for \(n=\dim V>2\), \(\fr_\kappa\) acts by \emph{outer} derivations, meaning that if \(a\in\fr_\kappa\), there is no \(Y\in\fs\) such that \(a\cdot X=\ad_Y X\) for all \(X\in\fs\). 

The \(R\)-symmetry groups for symmetric Dirac currents (hence superalgebras) were determined in \cite{Gall2021}. In that work, the choice of Dirac current is factored into a choice of reality conditions on complex spinor modules carrying a complex Dirac current. Consulting \cite[Table 10]{Gall2021}, we see that \(R\)-symmetry groups are classical Lie groups:\footnote{see \cite[App.A2]{Morris2015-arxiv} for the definitions of these groups.} as in the standard physics treatment in Lorentzian signature \cite{Strathdee1986,VanProeyen1999}, the groups \(\Orth(N)\), \(\U(N)\), \(\Sp(N)=\USp(2N)\) appear, but with more unusual choices of reality conditions, \(\Orth(r,s)\), \(\U(r,s)\), \(\Sp(r,s)\) (and products) are possible; in other signatures one can also have \(\GL(N,\RR)\), \(\SO(N,\CC)\), \(\Sp(N,\RR)\), \(\Sp(N,\CC)\), \(\U^*(2N)=\GL(N,\HH)\), \(\SO^*(2N)=\SO(N,\HH)\), and certain products. To the author's knowledge, this analysis has not been done for the skew-symmetric Dirac current (Lie algebra) case, but one expects a similar list. We observe here that in Lorentzian signature, it is always possible to choose a symmetric Dirac current with \emph{compact} \(R\)-symmetry group, which we will make use of in later sections (see e.g. Lemma~\ref{lemma:hhat-faithful-wlog}).

\subsubsection{Extending the flat model}
\label{sssec:isometry-alg-r-alg}

As previously alluded to, the flat model (super)algebra \(\fs\) admits an extension by the \(R\)-symmetry algebra \(\fr\) (where we now suppress the \(\kappa\) subscript) which we define as follows.

\begin{definition}\label{def:flat-model-alg-rsymm}
	Let \((V,\eta)\) be a (pseudo-)inner product space, let \(S\)  be a (possibly \(N\)-extended) spinor module and let \(\kappa:\Odot^2 S\to V\) (resp. \(\kappa:\Wedge^2 S\to V\)) be a Dirac current as in Definition~\ref{def:flat-model-alg}. Then the \emph{\(\fr\)-extended flat model superalgebra} (resp. \emph{algebra}) is the \(\ZZ\)-graded Lie superalgebra (resp. Lie algebra) with underlying \(\ZZ\)-graded vector space
	\begin{align}\label{eq:grading-fshat}
		&\fshat=\fshat_{-2}\oplus\fshat_{-1}\oplus\fshat_{0}
		&&\text{where}
		&\fshat_{-2} &= V,
		&\fshat_{-1} &= S,
	    &\fshat_0 &= \fso(V)\oplus \fr,
	    &\fshat_i &= 0 \text{ otherwise},
	\end{align}
	and brackets given by \eqref{eq:flat-model-bracket} and
	\begin{align}\label{eq:flat-r-symm-bracket}
	    \comm{a}{b} &= \comm{a}{b}_\fr,
	    &\comm{a}{A} &= 0,
	    &\comm{a}{v} &= 0,
	    &\comm{a}{\epsilon} &= a \epsilon,
	\end{align}
	where \(a,b\in\fr\), \(A\in\fso(V)\), \(v\in V\), \(\epsilon\in S\).
\end{definition}

Note that, by construction, \(\fs\) is a \(\ZZ\)-graded subalgebra of \(\fshat\) which differs only in the zeroth degree, and the Jacobi identities are satisfied precisely because those of \(\fs\) are and \(\fr\) acts by derivations on \(\fs\). The subspace \(\fshat_-=\fs_-:=V\oplus S\) is an ideal subalgebra of both \(\fs\) and \(\fshat\) which is sometimes known as the \emph{supertranslation ideal} \cite{Figueroa-OFarrill2016,Figueroa-OFarrill2017,Figueroa-OFarrill2017_1}, although we will not use this term in order to avoid confusion with the use of the term ``supertranslation'' in the literature on asymptotic symmetries.

\subsection{The spin-\(R\) group and spin-\(R\) structures}
\label{sec:spin-R}

As mentioned in the introduction, we seek to define Killing (super)algebras for spinors associated to a generalised spin structure (a spin-\(G\) structure), but the construction will require a Dirac current which is invariant under the action of the internal symmetry group \(G\). But such an action is nothing but a morphism of Lie groups \(G\to R\), so without loss of generality we need only consider the \(G=R\) case. 

As such, let us fix an inner product space \((V,\eta)\), a (possibly \(N\)-extended) real spinor module \(S\) of \(\Spin(V)\) and a Dirac current \(\kappa:\Odot^2S\to V\) or \(\kappa:\Wedge^2S\to V\) and denote the \(R\)-symmetry group by \(R\) and its Lie algebra by \(\fr\). We will adapt the relevant parts of the more general treatment found in \cite{Beckett2025_gen_spin} to the spin-\(R\) case. See also the references therein, particularly \cite{Nagase1995,Lazaroiu2019_1}.

\subsubsection{The spin-\(R\) group}
\label{sec:spin-R-group}

Recall that \(\Spin(V)\) has a central \(\ZZ_2\) subgroup which is the kernel of the canonical morphism \(\pi:\Spin(V)\to\SO(V)\):
\begin{equation}\label{eq:spin-cover-so}
\begin{tikzcd}
	1 \ar[r] &\ZZ_2 \ar[r] &\Spin(V) \ar[r,"\pi"] &\SO(V) \ar[r] &1.
\end{tikzcd}
\end{equation}
We will denote the non-trivial element of the \(\ZZ_2\) subgroup of \(\Spin(V)\) by \(-\1\). Denoting the spinor representation map by \(\sigma: \Spin(V)\to \GL(S)\), we have \(\sigma(\pm\1)=\pm\1\), where we also denote the identity in \(\GL(S)\) by \(\1\). The \(R\)-symmetry group consists of automorphisms of this representation, so \(\sigma(a)r=r\sigma(a)\) for all \(r\in R\), \(a\in \Spin(V)\), and we also have a central \(\ZZ_2\) subgroup \(\{\1,-\1\} \subseteq R\).

The subgroup \(\ZZ_2\cong\{(\pm\1,\pm\1)\}\subseteq \Spin(V)\times R\) is contained in the kernel of the action morphism \(\Spin(V)\times R\to \GL(S)\) given by \((a,r)\mapsto \sigma(a) r = r\sigma(a)\), so this action factors through the \emph{spin-\(R\) group} (of \((V,\eta,S,\kappa)\)), defined as
\begin{equation}
	\Spin^R(V,\eta) := \faktor{\Spin(V,\eta)\times R}{\{(\pm\1,\pm\1)\}}.
\end{equation}
We typically omit \(\eta\) in the notation where this is not ambiguous, and we also write \(\Spin^R(s,t):=\Spin^R(\RR^{s,t})\) for a spin-\(R\) group of a standard inner product space. This construction gives us the following commutative diagram:
\begin{equation}\label{eq:spin-R-reps-diagram}
\begin{tikzcd}
	\Spin(V) \ar[rd, hook] \ar[rrrd, bend left=15,"\sigma"] &&&\\
	& \Spin(V)\times R \ar[r, two heads, "\nu"] & \Spin^R(V) \ar[r,"\sigmahat"] & \GL(S)\\
	R \ar[ru, hook] \ar[rrru, bend right=15, hook] &&&
\end{tikzcd}.
\end{equation}
Here, the unlabelled maps \(\Spin(V)\hookrightarrow\Spin(V)\times R\), \(R\hookrightarrow\Spin(V)\times R\) and \(R\hookrightarrow\GL(S)\) denote natural inclusions, \(\nu:\Spin(V)\times R\to\Spin^R(V)\) denotes the natural quotient and \(\sigmahat: \Spin^R(V)\to \GL(S)\) denotes the action of \(\Spin^R(V)\) on \(S\) given by 
\begin{equation}
	[a,r]\cdot s = \sigmahat([a,r]) = \sigma(a)rs = r\sigma(a)s,
\end{equation}
where \([a,r]:=\nu(a,r)\in\Spin^R(V)\) denotes the image of a pair \((a,r)\in\Spin(V)\times R\). Considering an action of \(\ZZ_2\) by multiplication by \(\pm \1\) on the right of \(\Spin(V)\) and on the left of \(R\), we can also view \(\Spin^R(V)\) as the balanced product \(\Spin(V)\times_{\ZZ_2}R\). With either description, we thus have a short exact sequence
\begin{equation}\label{eq:spin-R-SES-defining}
\begin{tikzcd}
	1 \ar[r] &\ZZ_2 \ar[r] &\Spin(V)\times R \ar[r,"\nu"] &\Spin^R(V) \ar[r] &1,
\end{tikzcd}
\end{equation}
whence \(\nu\) is a double covering of \(\Spin^R(V)\). We have two natural injective Lie group morphisms
\begin{align}
	& \Spin(V) \xhookrightarrow{\quad} \Spin^R(V),	& a&\longmapsto [a,\1],\\
	& R \xhookrightarrow{\quad} \Spin^R(V),			& r&\longmapsto [\1,r],
\end{align}
and two surjective morphisms
\begin{align}
	\pihat&: \Spin^R(V) {\relbar\joinrel\twoheadrightarrow} \SO(V),		& [a,r]&\longmapsto \pi(a),\\
	\pi_R &: \Spin^R(V) {\relbar\joinrel\twoheadrightarrow} \Rbar = \faktor{R}{\ZZ_2},	& [a,r]&\longmapsto \bar{r},
\end{align}
where we note that we have defined \(\Rbar=R/\ZZ_2\) and \(\bar{r}=r\ZZ_2 \in \Rbar\). Together with the sequence \eqref{eq:spin-cover-so}, these maps fit into the following diagram with exact rows and columns:
\begin{equation}\label{eq:spin-R-SES-square}
\begin{tikzcd}
	 & 1 \ar[d] & 1\ar[d]	\\
	1 \ar[r] & \ZZ_2 \ar[r]\ar[d] & \Spin(V) \ar[r,"\pi"]\ar[d] & \SO(V) \ar[r]\ar[d,equal] & 1	\\
	1 \ar[r] & R \ar[r]\ar[d,"p_R"] & \Spin^R(V) \ar[r,"\pihat"]\ar[d,"\pi_R"] & \SO(V) \ar[r] & 1	\\
	 & \Rbar \ar[r,equal]\ar[d] & \Rbar \ar[d]	\\
	 & 1 & 1
\end{tikzcd}.
\end{equation}
The map \(\pihat\) will be particularly important in our discussion of spin-\(R\) structures. 

We also have the short exact sequence
\begin{equation}\label{eq:spin-R-SES-SO-RmodZ2}
\begin{tikzcd}
	1 \ar[r] &\ZZ_2 \ar[r] &\Spin^R(V) \ar[r,"\pihat\times\pi_R",pos=0.4] &\SO(V)\times\Rbar \ar[r] &1
\end{tikzcd}
\end{equation}
where the surjection is given by \([a,r]\mapsto(\pi(a),\bar{r})\) and the kernel \(\ZZ_2\) is generated by the element \([\1,-\1]=[-\1,\1]\in\Spin^R(V)\); in particular, \(\Spin^R(V)\) is a two-sheeted cover of \(\SO(V)\times\Rbar\). We then have a commutative diagram of group covering maps
\begin{equation}\label{eq:spin-R-double-covers}
\begin{tikzcd}
	\Spin(V)\times R \ar[rr, two heads, "\nu"']\ar[rrrr,bend left=12, two heads,"\pi\times p_R"] 
	&& \Spin^R(V) \ar[rr, two heads, "\pihat\times\pi_R"'] 
	&& \SO(V)\times\Rbar,
\end{tikzcd}
\end{equation}
so the Lie algebras of these three groups are naturally isomorphic:
\begin{equation}
	\Lie\Spin^R(V)\cong \Lie(\Spin(V)\times R) \cong \Lie\qty(\SO(V)\times \Rbar) \cong \fso(V)\oplus \fr.
\end{equation}
Denoting the adjoint action of any group \(G\) by \(\Ad^G\), we have
\begin{equation}\label{eq:spin-R-adjoint}
	\Ad^{\Spin^R(V)}_{[a,r]}(X,x) 
		= \qty(\Ad^{\Spin(V)}_a X, \Ad^R_r x)
		= \qty(\Ad^{\SO(V)}_{\pi(a)} X, \Ad^{\Rbar}_{\bar{r}} x)
\end{equation}
for all \(a\in \Spin(V)\), \(r\in R\), \(X\in\fso(V)\) and \(x\in \fr\).

\subsubsection{Spin-\(R\) structures}
\label{sec:spin-R-structures}

Let \((M,g)\) be an oriented pseudo-Riemannian manifold of signature \((s,t)\) and \(F_{SO}\to M\) its special orthonormal frame bundle, which has the structure of an \(\SO(s,t)\)-principal bundle. Recall that a spin structure on \((M,g)\) is a \(\Spin(s,t)\)-principal bundle \(P\to M\) with a \(\Spin(s,t)\)-equivariant bundle map \(\varpi:P\to F_{SO}\) where \(\Spin(s,t)\) acts on \(F_{SO}\) via the natural covering group morphism \(\pi:\Spin(s,t)\to\SO(s,t)\). Now note that \(\Spin^R(s,t)\) also acts on \(F_{SO}\) via the natural map \(\pihat:\Spin^R(s,t)\to\SO(s,t)\).

\begin{definition}[Spin-\(R\) structure]\label{def:spin-R-struc}
	A \emph{spin-\(R\) structure}\footnote{
		We note that our notion of spin-\(R\) is closely related to but distinct from the notion of \emph{generalised spin-\(r\)}, or \emph{spinorially twisted spin} structures, of e.g. \cite{Espinosa2016,Albanese2021,Artacho2023} which are spin-\(G\) structures with \(G=\Spin(r)\) for \(r\in\mathbb{N}\).
		} 
	on an oriented pseudo-Riemannian manifold \((M,g)\) of signature \((s,t)\) is a \(\Spin^R(s,t)\)-principal bundle \(\Phat\to M\) with a \(\Spin^R(s,t)\)-equivariant bundle map \(\varpihat:\Phat \to F_{SO}\).
	
	Two such structures \(\varpihat:\Phat\to M\), \(\varpihat':\Phat'\to M\) are \emph{equivalent} if there exists an isomorphism of \(\Spin^R(s,t)\)-principal bundles \(\Phi:\Phat\to\Phat'\) such that the following diagram commutes:
	\begin{equation}
	\begin{tikzcd}
		& \Phat\ar[rr,"\Phi"]\ar[dr,"\varpihat"']	&	& \Phat' \ar[dl,"\varpihat'"]\\
		&	& F_{SO} &
	\end{tikzcd}.
	\end{equation}
	We will call \((M,g)\) a \emph{spin-\(R\) manifold} if it is equipped with a spin-\(R\) structure.
\end{definition}

Unlike a spin structure, a spin-\(R\) structure is not in general a covering of \(F_{SO}\); instead, \(\varpihat: \Phat \to F_{SO}\) is an \(R\)-principal bundle on \(F_{SO}\), where \(R\) acts on \(\Phat\) via the natural homomorphism \(R\hookrightarrow \Spin^R(s,t)\). On the other hand \(\Phat\) can be viewed as a cover of a larger bundle as follows. The homomorphism \(\pi_R:\Spin^R(s,t)\to \Rbar\) induces a \(\pi_R\)-equivariant bundle map
\begin{equation}
	\varpi_R: \Phat \longrightarrow \Qbar := \Phat\times_{\pi_R}\Rbar,
\end{equation}
and we call \(\Qbar\to M\) the \emph{canonical \(\Rbar\)-bundle} -- note that we can consider \(\varpi_R\) as a \(\Spin(s,t)\)-principal bundle structure over \(\Qbar\). The map \(\pihat\times\pi_R:\Spin^R(s,t)\to\SO(s,t)\times \Rbar\) in the short exact sequence \eqref{eq:spin-R-SES-SO-RmodZ2} induces a two-sheeted covering map
\begin{equation}
	\varpihat\times\varpi_R:\Phat \longrightarrow F_{SO}\times_M\Qbar.
\end{equation} 
We can summarise the relationships between these principal bundles, each of which can be considered as an associated bundle to \(\Phat\), in the following diagram:
\begin{equation}
\begin{tikzcd}[sloped]
	\Phat \ar[dr,"\ZZ_2",dashed]
		\ar[drr,"R",bend left=12] 
		\ar[ddr,"{\Spin}"', bend right=12] 
		\ar[ddrr,"{\Spin^R}",controls={+(5,1) and +(1,1)},dotted] & &	\\
	& F_{SO}\times_M \Qbar \ar[r,"\Rbar"]\ar[d,"{\SO}"] 
		\ar[dr,"{\SO\times\Rbar}",dotted]	& F_{SO} \ar[d,"{\SO}"]	\\
	& \Qbar \ar[r,"\Rbar"]	& M	
\end{tikzcd}
\end{equation}
where each arrow represents a principal bundle over the target and is labelled by the structure group (omitting \((s,t)\) for clarity); equivalently, each represents a lift of the structure group along an epimorphism of Lie groups and is labelled with the kernel of that morphism. The solid and dashed arrows of the diagram form a commutative sub-diagram, and the dotted and dashed arrows form a separate commutative triangle. Each of these maps is also a \(\Spin^R(s,t)\)-equivariant bundle morphism.

If we have a \emph{spin} structure \(\varpi:P\to F_{SO}\) and a principal \(R\)-bundle \(Q\to M\), we can define a spin-\(R\) structure as follows. The fibred product \(P\times_M Q\to M\) is a \(\Spin(s,t)\times R\)-principal bundle on \(M\), and the canonical homomorphisms \(\nu:\Spin(s,t)\times R\to\Spin^R(s,t)\) and \(\pihat:\Spin^R(s,t)\to\SO(s,t)\) induce the spin-\(R\) structure
\begin{equation}\label{eq:spin-R-reducible-PQ}
	\varpihat_{PQ}:\Phat_{PQ}:= (P\times_M Q)\times_\nu\Spin^R(s,t) \longrightarrow F_{SO},
\end{equation}
where \(\varpihat_{PQ}([(p,q),[a,r]])=\varpi(p)\cdot\pihat([a,r])=\varpi(p\cdot a)\). In this case, there is a double cover \(Q\to\Qbar\) which is equivariant with respect to \(p_R:R\to\Rbar\), hence our notation for \(\Qbar\). Of course, not all spin-\(R\) structures arise this way. We call those that do \emph{reducible} \cite{Beckett2025_gen_spin}.

\subsubsection{Representations of the spin-\(R\) group and associated bundles}
\label{sec:spin-R-reps-bundles}

With our definition of spin-\(R\) structures in hand, we can now define an associated bundle of spinors. As with spin manifolds, many other natural bundles over a spin-\(R\) manifold can be described as associated bundles to the spin-\(R\) structure, essentially lifting their structure group from the special orthonormal frame bundle to the spin-\(R\)-bundle. Before discussing this in detail, we take a short detour to discuss the representation theory of the spin-\(R\) group. See \cite{Beckett2025_gen_spin} for more details.

Any representation of \(\Spin^R(V)\) pulls back along \(\nu:\Spin(V)\times R\to\Spin^R(V)\) and thus also pulls back to representations of \(\Spin(V)\) and of \(R\). Conversely, two representations \(\varrho_1: \Spin(V)\to \GL(W)\) and \(\varrho_2:R\to \GL(W)\) on the same space \(W\) assemble into a representation \(\varrho_1\varrho_2:\Spin(V)\times R\to\GL(W)\) given by \((\varrho_1\varrho_2)(a,r)=\varrho_1(a)\circ\varrho_2(r)\) if and only if \(\varrho_1(a)\circ\varrho_2(r)=\varrho_2(r)\circ \varrho_1(a)\) for all \(a\in\Spin(V)\) and \(r\in R\), and \(\varrho_1\varrho_2\) factors through a representation \(\varrho:\Spin^R(V)\to \GL(W)\) (such that \(\varrho\circ\nu=\varrho_1\varrho_2\)) if and only if \(\varrho_1(-\1)=\varrho_2(-\1)\).

\begin{definition}
	A representation \(\varrho:\Spin^R(V)\to \GL(W)\) with pull-backs \(\varrho_1:\Spin(V)\to\GL(W)\) and \(\varrho_2:R\to\GL(W)\) is \emph{even} if \(\varrho_1(-\1)=\varrho_2(-\1)=\Id_W\) and it is \emph{odd} if \(\varrho_1(-\1)=\varrho_2(-\1)=-\Id_W\).
\end{definition}

For example, the spinor representation \(\sigmahat:\Spin^R(V)\to\GL(S)\) is odd, while \(\pihat:\Spin^R(V)\to\SO(V)\) gives us an even representation on \(V\) in which \(R\) acts trivially, and the adjoint representation of \(\Spin^R(V)\) is also even. Since
\begin{align*}
	\varrho_1(-\1)= \Id_W &\iff \text{\(\varrho_1\) factors through a representation} \quad \overline{\varrho}_1:\SO(V) \to \GL(W),\\
	\varrho_2(-\1)= \Id_W &\iff \text{\(\varrho_2\) factors through a representation} \quad \overline{\varrho}_2: \Rbar = \faktor{R}{\ZZ_2} \to \GL(W),
\end{align*}
even representations of \(\Spin^R(V)\) are those which factor through a representation \(\overline\varrho=\overline\varrho_1\overline\varrho_2\) of \(\SO(V)\times \Rbar\) such that \(\overline\varrho\circ(\pihat\times\pi_R)=(\overline\varrho_1\circ\pihat)(\overline\varrho_2\circ\pi_R)=\varrho\). It can easily be shown that any representation of \(\Spin^R(V)\) is a direct sum of even and odd representations, and consequently that irreducible representations are either even or odd.

Now let \((M,g)\) be a pseudo-Riemannian manifold equipped with a spin-\(R\) structure \(\varpihat:\Phat\to F_{SO}\). We can also define parity for associated vector bundles of the principal bundle \(\Phat\to M\): the bundle \(\Phat\times_\varrho W \to M\) is \emph{even} or \emph{odd} if the representation \(\varrho:\Spin^R(s,t)\to\GL(W)\) that induces it is. 

The most important odd bundle is the one associated to \(\sigmahat:\Spin^R(s,t)\to\GL(S)\),
\begin{equation}
	\Shatbundle := \Phat\times_{\sigmahat} S \longrightarrow M,
\end{equation}
which we call the \emph{(twisted) spinor bundle}; we denote its space of sections by \(\fShat=\Gamma(\Shatbundle)\) and call this the space of \emph{(twisted) spinor fields}. 
We note that in the case that the spin-\(R\) structure is reducible to the product of a spin structure and an \(R\)-principal product as in \eqref{eq:spin-R-reducible-PQ}, this need not coincide with the ordinary spinor bundle \(\Sbundle:=P\times_\sigma S\to M\) unless \(Q\) is the trivial \(R\)-bundle, hence our use of the term ``twisted'' to distinguish the two.

The parity of associated bundles is important for gauge theory; if \(\varrho:\Spin^R(s,t)\to\GL(W)\) is an even representation then we have a natural isomorphism of vector bundles
\begin{equation}
	\Phat \times_\varrho W \cong \qty(F_{SO}\times_M \Qbar)\times_{\overline\varrho} W
\end{equation}
where we recall that \(\Qbar=\Phat\times_{\pi_R}\Rbar\to M\) is the canonical \(\Rbar\)-bundle, while if \(\varrho\) is odd, no such isomorphism exists. Note that if either \(R\) or \(\Spin(s,t)\) act trivially on \(W\), we have
\begin{align}
	& \Phat \times_\varrho W \cong F_{SO}\times_{\overline{\varrho}_1} W
	&& \text{or}
	&& \Phat \times_\varrho W \cong \Qbar\times_{\overline{\varrho}_2} W
\end{align}
respectively, as vector bundles. In particular, we can realise \(TM\) as an even bundle
\begin{equation}\label{eq:TM-assoc-spin-R}
	\Phat \times_{\pihat} \RR^{s,t} \cong F_{SO} \times_{\Id_{SO(s,t)}} \RR^{s,t} \cong  TM,
\end{equation}
by mapping \([p,v] \mapsto \qty[\varpihat(p)=(\varpihat(p)_i)_{i=1}^{s+t},v=(v^j)_{j=1}^{s+t}] \mapsto \sum_{i=1}^{s+t}v^i\varpihat(p)_i\) where we recall that \(\varpihat(p)=(\varpihat(p)_i)_{i=1}^{s+t}\) is nothing but an orthonormal basis for \(T_xM\) where \(p\) lies in the fibre \(\Phat_x\) above \(x\in M\). We can define similar isomorphisms for \(T^*M\), \(\Wedge^\bullet T^*M\), etc. Particularly important for connections on a spin-\(R\) structure, to be discussed soon, is the sequence of vector bundle isomorphisms
\begin{equation}
\begin{split}
	\Phat\times_{\Ad^{\Spin^R(s,t)}}\qty(\fso(s,t)\oplus\fr)
	& \cong \qty(\Phat\times_{\Ad^{\SO(s,t)}\circ\pihat} \fso(s,t))\oplus_M \qty(\Phat\times_{\Ad^{\Rbar}\circ\pi_R} \fr)	\\
	& \cong \qty(F_{SO}\times_{\Ad^{\SO(s,t)}} \fso(s,t))\oplus_M \qty(\Qbar\times_{\Ad^{\Rbar}} \fr)
\end{split}
\end{equation}
where we recall the relation between adjoint representations described in equation \eqref{eq:spin-R-adjoint}; that is, 
\begin{equation}\label{eq:adjoint-bundles-isom}
	\ad \Phat\cong \ad F_{SO} \oplus_M \ad \Qbar.
\end{equation}
Note that \(\ad F_{SO}\) and \(\ad \Qbar\) can also be considered as (even) associated bundles of \(\Phat\) by restricting the adjoint representation of \(\Spin^R(s,t)\). Any representation \(\varrho:\Spin^R(s,t)\to\GL(W)\) induces, via (left) conjugation on \(\End W\), a representation \(\End\varrho:\Spin^R(s,t)\to\GL(\End W)\) for which we have
\begin{equation}
	\Phat \times_{\End\varrho}{\End W} \cong \End\qty(\Phat \times_\varrho W),
\end{equation}
where the latter is the bundle of fibre-wise endomorphisms; in particular, the Lie algebra embedding \(\sigmahat_*:=\fso(s,t)\oplus\fr\hookrightarrow \End S\) induced by the twisted spinor representation \(\sigmahat\) induces a vector bundle embedding
\begin{equation}\label{eq:adPhat-embed-EndS}
	\ad\Phat\cong \ad F_{SO} \oplus_M \ad \Qbar \xhookrightarrow{\hspace{5mm}} \End\Shatbundle
\end{equation}
which allows us to define actions of sections of \(\ad F_{SO}\) and of \(\ad\Qbar\) on \(\End\Shatbundle\) fibre-wise via \(\sigmahat\).

\subsubsection{Connections and curvatures}
\label{sec:connections-rsymm}

Let \(\varpihat:\Phat\to F_{SO}\) be a spin-\(R\) structure on \((M,g)\). Recall that \(\Lie\Spin^R(s,t)\cong \fso(s,t)\oplus \fr\) and that its adjoint representation is given by equation \eqref{eq:spin-R-adjoint}. A (principal or Ehresmann) connection on \(\Phat\to M\) is a Lie algebra-valued 1-form \(\eA\in\Omega^1(\Phat;\fso(s,t)\oplus\fr)\) satisfying
\begin{equation}
\begin{aligned}
	R_{[a,r]}^*\eA &= \Ad^{\Spin^R(s,t)}_{[a,r]^{-1}}\circ\eA
	&& \text{and}
	& \eA(\widehat\xi_{(X,x)}) &= (X,x)
\end{aligned}
\end{equation}
for all \(a\in\Spin(s,t)\), \(r\in R\), \(X\in\fso(s,t)\) and \(x\in \fr\), where \(R_{[a,r]}:\Phat\to\Phat\) denotes the right action of \([a,r]=\pihat(a,r)\) and \(\widehat\xi_{(X,x)}\in\fX(\Phat)\) denotes the fundamental vector field generated by \((X,x)\). Letting \(\varpi_R:\Phat\to \Qbar=\Phat\times_{\pi_R}\Rbar\) denote projection to the canonical \(\Rbar\)-bundle, one can show that
\begin{equation}\label{eq:principal-connections-add}
	\eA = \varpihat^*\omega + \varpi_R^* \alpha,
\end{equation}
where \(\omega\in \Omega^1(F_{SO};\fso(s,t))\), \(\alpha\in \Omega^1(\Qbar;\fr)\) are principal connections, and conversely, any two such connections on \(F_{SO}\to M\) and \(\Qbar\to M\) induce a connection on \(\Phat\to M\). In what follows, we will demand that \(\omega\) is the Levi-Civita connection; equivalently, we can assign a torsion to \(\eA\) using the isomorphism of vector bundles in equation \eqref{eq:TM-assoc-spin-R} and demand that this torsion vanishes. This torsion is equal to that of the associated metric connection \(\omega\), so vanishes if and only if \(\omega\) is the Levi-Civita connection.

The curvature 2-forms \(\Omega_\eA\in\Omega^2(\Phat;\fso(s,t)\oplus\fr)\), \(\Omega_\omega\in\Omega^2(F_{SO};\fso(s,t))\), and \(\Omega_\alpha\in\Omega^2(\Qbar;\fr)\) satisfy
\begin{equation}\label{eq:principal-curvatures-add}
	\Omega_\eA = \varpihat^*\Omega_\omega + \varpi_R^*\Omega_\alpha,
\end{equation}
and if \(\omega\) is the Levi-Civita connection, \(\Omega_\omega\) is of course the Riemann 2-form. The curvature can of course be encoded in a different way; by standard theory of principal connections, these 2-forms are \emph{basic} (horizontal and invariant). whence they descend to sections \(\Rhat\in\Omega^2(M;\ad\Phat)\), \(R\in\Omega^2(M;\ad F_{SO})\), and \(F\in\Omega^2(M;\ad \Qbar)\) on \(M\) satisfying
\begin{equation}\label{eq:base-curvatures-add}
	\Rhat(X,Y) = R(X,Y) + F(X,Y)
\end{equation}
for all \(X,Y\in\fX(M)\), where we implicitly use the isomorphism of vector bundles in \eqref{eq:adjoint-bundles-isom}.

The principal connection \(\eA\) on \(\Phat\to M\) induces Koszul connections on associated bundles, in particular on the  twisted spinor bundle \(\Shatbundle\to M\), which we will denote by \(\nablahat\) and call the \emph{twisted covariant derivative}. On even vector bundles, which we have seen can be considered as associated bundles to \(F_{SO}\times_M\Qbar\), we can write \(\nablahat=\nabla + \alpha\) in local trivialisations where \(\nabla\) denotes the covariant derivative of the Levi-Civita connection; where the action of \(R\) on the defining representation is trivial, in particular on \(TM\), we have \(\nablahat=\nabla\). More generally, it is sometimes convenient to locally abuse notation and write \(\nablahat=\nabla + \alpha\) even on non-even associated bundles, where neither \(\nabla\) nor \(\alpha\) are globally defined.

We then have the following formulae for the curvatures above considered as operators on associated bundles; note that we could have equivalently used these equations to \emph{define} \(\Rhat\) as the curvature of the covariant derivative \(\nablahat\) on \(\Shatbundle\) (which is simply the representation of \(\Omega_\eA\) on \(\Shatbundle\)), \(R\) as the Riemann curvature and \(F\) as the field strength of the local \(\fr\)-valued gauge fields \(\alpha\): for \(X,Y\in\fX(M)\),
\begin{equation}\label{eq:curvatures-derivative-def}
\begin{gathered}
	\Rhat(X,Y) 	= \comm{\nablahat_X}{\nablahat_Y} - \nablahat_{\comm{X}{Y}},	\\
	R(X,Y) 		= \comm{\nabla_X}{\nabla_Y} - \nabla_{\comm{X}{Y}},			\\
	F(X,Y)		= (\nabla_X\alpha)(Y) - (\nabla_Y\alpha)(X) + \comm{\alpha(X)}{\alpha(Y)},
	\qquad \text{or} \qquad
	F = d\alpha + \tfrac{1}{2}[\alpha\wedge\alpha]. 
\end{gathered}
\end{equation}
Note that although the LHSs are defined globally, as is the RHS of the first equation, the RHS of the second equation is well-defined globally only on sections of even associated bundles or in local trivialisations, and the third is valid only in local trivialisations.

\subsubsection{Covariant Lie derivative and Cartan calculus}
\label{sec:lie-der-rsymm}

Recall that for Killing vectors \(X\in\fiso(M,g)\), the endomorphism \(A_X:=-\nabla X\) of the tangent bundle (given by \(Y\mapsto \nabla_Y X\)) lies in \(\fso(M,g)\), the Lie algebra of \(g\)-skew-symmetric sections of \(\End(TM)\), which can also be identified with the space of sections of \(\ad F_{SO}\to M\). By equation \eqref{eq:adPhat-embed-EndS}, as a section of \(\ad F_{SO}\), \(A_X\) naturally acts pointwise (via \(\sigmahat_*:\fso(s,t)\oplus\fr\to\End(S)\)) on \(\fShat\). We can thus define the \emph{covariant (spinorial) Lie derivative} \(\eLhat_X: \fShat\to \fShat\) along the Killing vector \(X\) by
\begin{equation}\label{eq:eLhat-nablahat}
	\eLhat_X\epsilon = \nablahat_X\epsilon + A_X\cdot\epsilon.
\end{equation}
for all \(\epsilon\in\fShat\). This is a direct generalisation of the Kosmann spinorial Lie derivative \cite{Kosmann1971}. Such a derivative is introduced via local expressions in \cite[Sec.7]{deMedeiros2018} in the special case of signature \((1,5)\) with, in the terminology of the present work, a \emph{reducible} spin-\(h\) (that is, spin-\(\Sp(1)\)) structure with trivial \(\Sp(1)\)-bundle \(Q\to M\) (in particular, the twisted spinor bundle \(\Shatbundle\) coincides with the untwisted one \(\Sbundle\)); we have simply re-stated the definition, in global geometric language, for general spin-\(R\) structures.

In fact, the formula \eqref{eq:eLhat-nablahat} defines a covariant Lie derivative on sections of any associated bundle to \(\Phat\to M\) and agrees with the ordinary Lie derivative on even bundles for which \(R\) acts trivially on the defining representation (in particular on \(TM\) and \(\ad F_{SO}\)); moreover, it satisfies the Leibniz rule with respect to any products of sections or actions of sections on other sections which we will encounter, as well as the following \cite{Beckett2025_gen_spin}.

\begin{proposition}\label{prop:eLhat-props}
	The covariant Lie derivative described above satisfies the following properties:
	\begin{align}
		\comm{\eLhat_X}{\nablahat_Y} &= \nablahat_{\comm{X}{Y}} + F(X,Y) \label{eq:eLhat-nablahat-comm}
	\end{align}
	for \(X\in\fiso(M,g)\) and \(Y\in\fX(M)\) (equivalently, \(\comm{\eLhat_X}{\nablahat} = \imath_X F\) for all \(X\in\fiso(M,g)\)) and
	\begin{equation}\label{eq:eLhat-comm}
		\comm{\eLhat_X}{\eLhat_Y} = \eLhat_{\comm{X}{Y}} + F(X,Y)	
	\end{equation}
	for \(X,Y\in\fiso(M,g)\).\footnote{Note that the \(F(X,Y)\) terms in these equations are to be understood as operators acting on sections of an associated bundle via the appropriate representation of \(\fr\).}
\end{proposition}

\begin{proof}
	We will make use of some identities originally due to Kostant \cite{Kostant1955} but distilled in our present notation in \cite[\S{}A.2.1]{Beckett2024_ksa}. Since \(X\) is a Killing vector, \(\nabla_YA_X = R(X,Y)\) \cite[eq.(61)]{Beckett2024_ksa}, so
	\begin{equation}
	\begin{split}
		\comm{\eLhat_X}{\nablahat_Y} 
			&= \comm{\nablahat_X}{\nablahat_Y} + \comm{A_X}{\nablahat_Y}	\\
			&= \comm{\nablahat_X}{\nablahat_Y} - \nabla_Y A_X	\\
			&= \nablahat_{\comm{X}{Y}} + \Rhat(X,Y) - R(X,Y)	\\
			&= \nablahat_{\comm{X}{Y}} + F(X,Y),
	\end{split}
	\end{equation}
	where we have used the curvature identities \eqref{eq:base-curvatures-add} and \eqref{eq:curvatures-derivative-def}; this demonstrates the first property. A Leibniz rule, a formula for the Lie derivative \cite[eq.(59)]{Beckett2024_ksa} and some identities for Killing vectors \cite[eq.(61),(62)]{Beckett2024_ksa} give us
	\begin{equation}
		\comm{\eLhat_X}{A_Y}=\eL_X A_Y = \nabla_X A_Y + \comm{A_X}{A_Y} = - R(X,Y) + \comm{A_X}{A_Y} = A_{\comm{X}{Y}},
	\end{equation}
	and using this as well as the first property gives us
	\begin{equation}
		\comm{\eLhat_X}{\eLhat_Y}
			= \comm{\eLhat_X}{\nablahat_Y}  + \comm{\eLhat_X}{A_Y} 
			= \nablahat_{\comm{X}{Y}} + F(X,Y) + A_{\comm{X}{Y}} 
			= \eLhat_{\comm{X}{Y}} + F(X,Y),
	\end{equation}
	which demonstrates the second property.
\end{proof}

Recall that the choice of connection \(\eA\) on \(\Phat\) induces not just a covariant derivative on all associated bundles, but a \emph{covariant exterior derivative}, which we will denote by \(\dhat\), on differential forms with values in an associated bundle. With this notation, we have Bianchi identities \(\dhat\Rhat =0\) and \(\dhat F = 0\).

The ordinary Lie derivative, exterior derivative and interior derivative (contraction) on differential forms satisfy a number of identities collectively known as \emph{Cartan calculus} which are extremely useful for computations. It turns out that the covariant derivatives satisfy a similar \emph{covariant Cartan calculus}, which we will take advantage of.

\begin{proposition}[{Covariant Cartan calculus \cite{Beckett2025_gen_spin}}]\label{prop:cov-cartan-calculus}
	Let \(\varrho:\Spin^R(s,t)\to \GL(W)\) be a representation with associated bundle \(\Wbundle = \Phat\times_\varrho W \to M\). Then we have the following Cartan formula for the covariant Lie derivative along any Killing vector \(X\) when acting on \(\Omega^\bullet(M;\Wbundle):=\Gamma\qty(\Wedge^\bullet T^*M\otimes \Wbundle)\):
	\begin{equation}
		\eLhat_X = \imath_X\dhat + \dhat\imath_X.
	\end{equation}
	We also have the following identities of operators on \(\Omega^\bullet(M;\Wbundle)\):\footnote{
		Here and elsewhere, the wedge of a form with values in \(\ad\Qbar\) (or locally in \(\fr\)) with a form with values in \(\Wbundle\) is understood to include the action of the former on the latter via the representation of \(\fr\) on \(W\).
			}
	\begin{equation} \label{eq:eLhat-dhat-comm}
		\comm{\eLhat_X}{\dhat} = \imath_X F \wedge,
	\end{equation}
	for all \(X\in\fiso(M,g)\);
	\begin{equation}\label{eq:eLhat-i-comm}
		\comm{\eLhat_X}{\imath_Y} = \imath_{\comm{X}{Y}},
	\end{equation}
	for all \(X\in\fiso(M,g)\) and \(Y\in\fX(M)\);
	\begin{equation}\label{eq:eLhat-comm-forms}
		\comm{\eLhat_X}{\eLhat_Y} = \eLhat_{\comm{X}{Y}} + F(X,Y),
	\end{equation}
	for all \(X,Y\in\fiso(M,g)\).
\end{proposition}

\subsubsection{Time orientability and Dirac currents on twisted spinor bundles}
\label{sec:dirac-current-bundle-rsymm}

One final structure which we will require on spin-\(R\) manifolds is the \emph{bundle Dirac current}, that is, a bundle map \(\kappabold:\Wedge^2 \Shatbundle\to TM\) or \(\kappabold:\Odot^2 \Shatbundle \to TM\) which is essentially the Dirac current \(\kappa\) (as in Definition~\ref{def:flat-model-alg}) used to define the \(R\)-symmetry when restricted to fibres over a point. In indefinite signature, there is an obstruction to defining such a map in that by definition \(\kappa\) is \(\fso(s,t)\)-invariant and \(R\)-equivariant but not \(\Spin(s,t)\)-equivariant, and thus not \(\Spin^R(s,t)\)-equivariant, in general. As such, it is necessary to assume that \((M,g)\) is time-orientable so that structure groups can be reduced.

\begin{lemma}[\cite{Beckett2025_gen_spin}]\label{lemma:time-orientable-spin-R}
	Let \((M,g)\) be an oriented pseudo-Riemannian manifold of indefinite signature with spin-\(R\) structure \(\varpihat:\Phat\to F_{SO}\). Then \((M,g)\) is time-orientable if and only if there exists a reduction \(\Phat_0\to M\) of the structure group of the \(\Spin^R(s,t)\)-principal bundle \(\Phat\to M\) to the index-2 subgroup\footnote{
			Note that although \(\Spin_0(s,t)\) is connected, \(\Spin^R_0(s,t)\) may not be since we do not assume that \(R\) is connected.
			}
	\(\Spin^R_0(s,t):=\Spin_0(s,t)\times_{\ZZ_2} R\).
\end{lemma}

We then have the following, which is a generalisation of the result for spin manifolds \cite[Lem.3.3]{Beckett2024_ksa}. In particular examples, it may be the case that, without assuming time-orientability, one can define the Dirac current at a point and show that it happens to be \(\nablahat\)-holonomy-invariant, hence one can still extend it globally by parallel transport (if \(M\) is connected), but we cannot expect this in general.

\begin{lemma}[Existence of bundle Dirac current]\label{lemma:bundle-dirac-current-exist-rsymm}
	Let \((M,g)\) be a connected, oriented and (in indefinite signature) time-oriented pseudo-Riemannian manifold with spin-\(R\) structure \(\varpihat:\Phat\to F_{SO}\) and connection \(\eA\) on \(\Phat\) lifting the Levi-Civita connection. Then there exists a \emph{bundle Dirac current} \(\kappabold:\Wedge^2 \Shatbundle\to TM\) or \(\kappabold:\Odot^2 \Shatbundle \to TM\) satisfying \(\nablahat\kappabold=0\) and \(\eLhat_X\kappabold=0\) for all Killing vectors \(X\) (where we use the covariant derivatives induced by \(\eA\)) and \(a\cdot\kappa=0\) for all \(a\in\Gamma(\ad\Qbar)\).
\end{lemma}

\begin{proof}
	For definite signature, set \(G:=\Spin^R(s,t)\) and \(\Phat_0:=\Phat\); for indefinite signature, reduce the structure group of the spin-\(R\) structure as described above and set \(G:=\Spin^R_0(s,t)\). We give the proof for symmetric Dirac current; the skew-symmetric case is entirely analogous.
	
	Fixing \(x\in M\) and \(p\in(\Phat_0)_x\), we can view the holonomy group at \(x\) of the connection \(\eA\) as a subgroup of \(G\), and we can define a linear isomorphism  \(\Hom(\Odot^2S,\RR^{s,t})\cong E_x\), equivariant under the action of the holonomy group, where \(E_x\) is the fibre of the bundle
	\begin{equation}
		E = \Hom\qty(\Odot^2\Shatbundle,TM) \xlongrightarrow{\cong} \Phat_0\times_G\Hom\qty(\Odot^2S,\RR^{s,t})
	\end{equation}
	 by sending \(\lambda\mapsto[p,\lambda]\). Under this map, \(\kappa\) is sent to \(\kappabold_x=[p,\kappa]\) which is holonomy-invariant since \(\kappa\) is invariant under \(G\). Thus by the holonomy principle, \(\kappabold_x\) extends by parallel transport to a \(\nablahat\)-parallel global section \(\kappabold\) of \(E\), as required, and this is independent of the choice of \(p\) by \(G\)-invariance of \(\kappa\). Infinitesimal \(G\)-invariance also ensures that \(\kappabold\) is invariant under the pointwise action of \(\ad\Phat\), in particular \(\eLhat_X\kappabold=\nablahat_X\kappabold + A_X\cdot\kappabold = 0\) for all \(X\in\fiso(M,g)\) and \(a\cdot\kappabold=0\) for all \(a\in\Gamma(\ad\Qbar)\).
\end{proof}

\section{Killing spinors and Killing (super)algebras on generalised spin manifolds}
\label{sec:killing-rsymm}

We now seek to find a notion of Killing spinors and Killing (super)algebras, and to describe the algebraic structure of the latter, in the geometric setting of a spin-\(R\) structure, generalising the treatment of \cite[Sec.3]{Beckett2024_ksa}. We will first introduce an analogue of the Lie algebra of Killing vectors on the spin-\(R\) structure itself.

Throughout, we let \(S\) be a (possibly \(N\)-extended) real spinor module of \(\Spin(s,t)\) equipped with a Dirac current \(\kappa: \Odot^2 S\to\RR^{p,q}\) or \(\kappa: \Wedge^2 S\to\RR^{p,q}\) and \(R\) be the corresponding \(R\)-symmetry group. We furthermore let \((M,g)\) be a connected, oriented and (in indefinite signature) time-oriented pseudo-Riemannian manifold of signature \((s,t)\) with special orthonormal frame bundle \(F_{SO}\to M\), spin-\(R\) structure \(\varpihat:\Phat\to F_{SO}\) (Def.~\ref{def:spin-R-struc}) and twisted spinor bundle \(\Shatbundle=\Phat\times_{\sigmahat} S\to M\). We define the space of twisted spinor fields \(\fShat=\Gamma(\Shatbundle)\), identify the space of sections of \(\ad F_{SO}\) with the space of \(g\)-skew-symmetric sections of \(\End(TM)\), denoted \(\fso(M,g)\), and denote the Lie algebra of Killing vectors by \(\fiso(M,g)\). We also fix a torsion-free connection \(\eA\) on \(\Phat\).

\subsection{The symmetry algebra of a spin-\(R\) structure}
\label{sec:symm-alg-spin-r}

The even part of a Killing (super)algebra of the type described in \cite[Sec.3]{Beckett2024_ksa} consists of a subalgebra of the algebra of \(\fiso(M,g)\). We now seek to find an appropriate analogue of \(\fiso(M,g)\) in the spin-\(R\) setting, and we will later define the even part of the Killing (super)algebra to be a subalgebra of this structure. We call this structure the \emph{symmetry algebra} of the spin-\(R\) structure (and the connection \(\eA\)). We note that this symmetry algebra can be described quite independently of the existence of any Killing (super)algebra, however, and we have introduced it for a more general audience elsewhere \cite{Beckett2025_gen_spin}.

We begin by identifying some Lie algebras of sections on \(M\) which (one might expect to) have natural actions on the space of spinor fields \(\fShat\). Two such algebras present themselves in our geometric setting: the algebra of Killing vectors \(\fiso(M,g)\) (via the covariant Lie derivative) and the space of sections of \(\ad\Qbar\to M\), the (infinite-dimensional) Lie algebra of infinitesimal gauge transformations on the canonical \(\Rbar\)-bundle \(\Qbar=\Phat\times_{\pi_R}\Rbar\to M\). Let us denote the latter by \(\fR=\Gamma(\ad\Qbar)\); we will refer to its elements as \emph{infinitesimal \(R\)-symmetries}, or \emph{\(\fr\)-symmetries} for short. In a local trivialisation, these sections are functions with values in \(\fr\).

One might naïvely expect that the direct product Lie algebra \(\fiso(M,g)\oplus\fR\) or some subalgebra thereof is the algebra we seek; however, we note that the curvature term in \eqref{eq:eLhat-comm} means that the action of \(\fiso(M,g)\) on \(\fShat\) via the covariant Lie derivative does \emph{not} define a representation, so we must choose a modified bracket on \(\fiso(M,g)\oplus\fR\). We will denote this bracket by \(\ccomm{-}{-}\) in order to distinguish it from the ordinary Lie bracket of vector fields and other Lie brackets and commutators.

Let us first consider a natural non-trivial bracket between \(\fiso(M,g)\) and \(\fR\) given by the covariant Lie derivative defined in equation~\eqref{eq:eLhat-nablahat}. Such a derivative exists because \(\ad \Qbar\) is an associated bundle to the spin-\(R\) structure: \(\ad\Qbar\cong \Phat\times_{\Ad^{\Rbar}\circ\pi_R}\fr\to M\), where the defining representation \(\Ad^{\Rbar}\circ\pi_R:\Spin^R(s,t)\to \fr\) can also be viewed as the restriction of \(\Ad^{\Spin^R(s,t)}\) to the submodule \(\fr\). Since \(\fso(V)\) acts trivially on \(\fr\), for all \(X\in\fiso(M,g)\) and \(a\in\fR\) we have \(A_X\cdot a = \comm{A_X}{a} = 0\) in \(\fR\) (where \(A_X:=-\nabla X\in\fso(M,g)\)), so
\begin{equation}
	\eLhat_X a = \nablahat_X a + A_X\cdot a = \nablahat_X a + \comm{A_X}{a} = \nablahat_X a.
\end{equation}
Thus we will choose to set
\begin{equation}\label{eq:kvf-rsymm-bracket}
	\ccomm{X}{a} := \eLhat_X a = \nablahat_X a
\end{equation}
for a Killing vector \(X\in\fiso(M,g)\) and \(\fr\)-symmetry \(a\in\fR\). For \(a,b\in\fR\) let us take
\begin{equation}\label{eq:rsymmm-rsymm-bracket}
	\ccomm{a}{b} := \comm{a}{b}
\end{equation}
where the latter is the natural bracket on \(\fR\); that is, the one given pointwise by the bracket on \(\fr\). Thus the Jacobi identity for three elements of \(\fR\) is automatically satisfied, and by a Leibniz rule for the covariant derivative we have
\begin{equation}
	\ccomm{X}{\ccomm{a}{b}} = \nablahat_X\comm{a}{b} = \comm{\nablahat_X{a}}{b} + \comm{a}{\nablahat_Xb} = \ccomm{\ccomm{X}{a}}{b} + \ccomm{a}{\ccomm{X}{b}},
\end{equation}
which verifies another Jacobi identity for \(X\in\fiso(M,g)\) and \(a,b\in\fR\).

Let us now use another component of the Jacobi identity to motivate our definition of \(\ccomm{X}{Y}\); using the identity \eqref{eq:eLhat-comm} we find
\begin{equation}
	\ccomm{X}{\ccomm{Y}{a}} - \ccomm{Y}{\ccomm{X}{a}}
		= \comm{\eLhat_X}{\eLhat_Y}a
		= \eLhat_{\comm{X}{Y}} a + \comm{F(X,Y)}{a},
\end{equation}
for all \(X,Y\in\fiso(M,g)\) and \(a\in\fR\), so the Jacobi identity for \(X,Y,a\) is satisfied if we define
\begin{equation}\label{eq:kvf-kvf-bracket}
	\ccomm{X}{Y} := \comm{X}{Y} + F(X,Y) = \eL_X Y + F(X,Y).
\end{equation}
Now we need to check the Jacobi identity for three Killing vectors; we find that
\begin{equation}
	\ccomm{X}{\ccomm{Y}{Z}} = \eL_X\eL_Y Z + F(X,\eL_YZ) + \eLhat_X(F(Y,Z)),
\end{equation}
and, noting that \(\eLhat_X(F(Y,Z))=\nablahat_X(F(Y,Z))\), upon taking cyclic permutations of the above we find that this component of the Jacobi identity is equivalent to
\begin{equation}
	\qty(\nablahat_XF)(Y,Z) + \qty(\nablahat_YF)(Z,X) + \qty(\nablahat_ZF)(X,Y) = 0
\end{equation}
for all Killing vector fields \(X,Y,Z\). But the above is nothing but a form of the Bianchi identity \(\dhat F = 0\) for \(F\), so it is identically satisfied. This gives us the first part of the following result, a slightly more general version of which appears in \cite{Beckett2025_gen_spin}.

\begin{proposition}\label{prop:symm-alg-spin-r}
	The space \(\fsymm(\varpihat,\eA) := \fiso(M,g)\oplus\fR\) equipped with the bracket \(\ccomm{-}{-}\) given by \eqref{eq:kvf-rsymm-bracket}, \eqref{eq:rsymmm-rsymm-bracket} and \eqref{eq:kvf-kvf-bracket} is a Lie algebra. The space of spinor fields \(\fShat\) is a module of this algebra, where elements of \(\fiso(M,g)\) act via \(\eLhat\) and those of \(\fR\) act via the pointwise action of \(\fr\) on \(S\).
\end{proposition}

\begin{proof}
	The first claim follows from the discussion above. Equation \eqref{eq:eLhat-comm} gives us
	\begin{equation}
		\ccomm{X}{Y}\cdot\epsilon 
			= \eLhat_{\comm{X}{Y}}\epsilon + F(X,Y)\epsilon
			= \eLhat_X\qty(\eLhat_Y\epsilon) - \eLhat_Y\qty(\eLhat_X\epsilon)
	\end{equation}
	for all \(X,Y\in\fiso(M,g)\) and \(\epsilon\in\fShat\),	while the Leibniz rule for \(\eLhat_X\) gives 
	\begin{equation}
		\ccomm{X}{a}\cdot\epsilon
		= \qty(\eLhat_X a)\epsilon 
		= \eLhat_X (a\epsilon) - a \qty(\eLhat_X \epsilon)
	\end{equation}
	for \(X\in\fiso(M,g)\), \(a\in\fR\) and \(\epsilon\in\fShat\). We also have
	\begin{equation}
		\ccomm{a}{b}\cdot\epsilon
		= \comm{a}{b} \epsilon
		= a(b\epsilon) - b(a\epsilon)
	\end{equation}
	for all \(a,b\in\fR\) and \(\epsilon\in\fShat\), since the second equality holds pointwise for the action of \(\fr\) on \(S\). This shows that \(\fShat\) is a module for \(\fsymm(\varpihat,\eA)\).
\end{proof}

It is manifest in the definition of the bracket \(\ccomm{-}{-}\) that \(\fsymm(\varpihat,\eA)\) is \emph{not} the direct sum of Lie algebras \(\fiso(M,g)\oplus\fR\) in general, but rather an extension of \(\fiso(M,g)\) by \(\fR\); we have a short exact sequence of Lie algebras
\begin{equation}\label{eq:symm-alg-SES}
\begin{tikzcd}
	0 \ar[r] & \fR \ar[r] & \fsymm(\varpihat,\eA) \ar[r] & \fiso(M,g) \ar[r] & 0
\end{tikzcd}
\end{equation}
where the map \(\fR\to\fsymm(\varpihat,\eA)\) is the inclusion and \(\fsymm(\varpihat,\eA)\to\fiso(M,g)\) is projection. From another point of view, \(\fsymm(\varpihat,\eA)\) is a Lie algebra deformation of \(\fiso(M,g)\oplus\fR\).

The notation emphasises the dependence of the algebra, through the definitions of \(\eLhat\) and \(\nablahat\), on the connection \(\eA\). In particular, note that if \(\imath_XF=0\) for all \(X\in\fiso(M,g)\) then \(\ccomm{-}{-}\) restricted to \(\fiso(M,g)\) is simply the Lie bracket, and we have \(\fsymm(\varpihat,\eA)=\fiso(M,g)\ltimes\fR\). 

\subsection{Construction of Killing (super)algebras}
\label{sec:ksa-construction-rsymm}

We now seek to construct a Lie superalgebra or \(\ZZ_2\)-graded Lie algebra whose odd part is the space of spinor fields which are parallel with respect to some connection \(D\) on the twisted spinor bundle \(\fShat\),
\begin{equation}
	\fShat_D = \qty{ \epsilon\in\fShat \,\middle|\, D\epsilon = 0 }.
\end{equation}
We begin with some results on such connections.

\subsubsection{Connections on twisted spinor bundles}

Let \(D\) be a connection on \(\Shatbundle\). The difference between two connections on a bundle is a 1-form with values in the endomorphisms of the bundle, so there exists some \(\betabold\in\Omega^1\qty(M;\End\Shatbundle)\) such that 
\begin{equation}
	D=\nablahat-\betabold.	
\end{equation}
We denote the contraction of \(\betabold\) with a vector field \(X\in\fX(M)\) by \(\betabold_X=\imath_X\betabold\).

\begin{lemma}
	For all \(X\in\fiso(M,g)\),
	\begin{equation}\label{eq:Lhat-D-comm}
		\comm{\eLhat_X}{D} = \imath_X F - \eLhat_X\betabold,
	\end{equation}
	where the bracket is the commutator of operators on \(\fShat\).
\end{lemma}

\begin{proof}
	Using Proposition~\ref{prop:eLhat-props}, we have
	\begin{equation}
		\comm{\eLhat_X}{D} 
			= \comm{\eLhat_X}{\nablahat} - \comm{\eLhat_X}{\betabold}
			= \imath_XF - \eLhat_X\betabold
	\end{equation}
	for all \(X\in\fiso(M,g)\).
\end{proof}

We denote the curvature 2-form \(R^D\in\Omega^2(M;\End\Shatbundle)\) of \(D\) by
\begin{equation}
		R^D(X,Y) = \comm{D_X}{D_Y} - D_{\comm{X}{Y}}
\end{equation}
for all \(X,Y \in \fX(M)\). Recalling now the embedding of adjoint bundles into \(\End\Shatbundle\) in equation \eqref{eq:adPhat-embed-EndS}, we have the following.

\begin{proposition}\label{prop:D-curvature-rsymm}
Considering each term as a 2-form with values in the bundle \(\End\Shatbundle\), the curvature \(R^D\) is given by
\begin{equation}
	R^D = R + F + \tfrac{1}{2}[\betabold\wedge\betabold] - \dhat\betabold
\end{equation}
where \(R\) is the Riemann curvature, \(F\) is the \(\fr\)-symmetry field strength (the \(\ad\Qbar\)-component of the curvature of the torsion-free connection \(\eA\) on the \(\Spin^R(s,t)\)-bundle \(\Phat\to M\)), and \(\dhat\) is the covariant exterior derivative.
\end{proposition}

\begin{proof}
For all \(X,Y\in\fX(M)\), using definitions and the Leibniz rule,
\begin{equation}
\begin{split}
	&R^D(X,Y)	\\
		&= \comm{\nablahat_X}{\nablahat_Y} + \comm{\betabold_X}{\betabold_Y} - \comm{\nablahat_X}{\betabold_Y} - \comm{\betabold_X}{\nablahat_Y} - \nablahat_{\comm{X}{Y}} + \betabold_{\comm{X}{Y}}	\\
		&= \Rhat(X,Y) + \comm{\betabold_X}{\betabold_Y} - \nablahat_X\qty(\betabold_Y) + \nablahat_Y\qty(\betabold_X) + \betabold_{\comm{X}{Y}}\\
		&= R(X,Y) + F(X,Y) + \comm{\betabold_X}{\betabold_Y} - (\dhat\betabold)(X,Y).
\end{split}
\end{equation}
The identity follows by abstracting the vector fields \(X,Y\).
\end{proof}

%The identity \(R^D\epsilon=0\) for \(\epsilon\in\fShat_D\) is an integrability condition for the existence of \(D\)-parallel sections of \(\Shatbundle\).

\subsubsection{Existence of the Killing (super)algebra}

The even part of our (super)algebra must be a Lie algebra which has \(\fShat_D\) as a module. Since we have already seen that the symmetry algebra \(\fsymm(\varpihat,\eA)\) of the spin-\(R\) structure discussed in \S\ref{sec:symm-alg-spin-r} acts on the space of spinor fields \(\fShat\) (Proposition~\ref{prop:symm-alg-spin-r}), the natural starting point is to find a subalgebra of \(\fsymm(\varpihat,\eA)\) which preserves \(\fShat_D\). Our point of departure is the identity \eqref{eq:Lhat-D-comm}.

Recall that \(\ad\Qbar\hookrightarrow\End\Shatbundle\) (see equation \eqref{eq:adPhat-embed-EndS}), so we can consider \(a\in\fR\) as a section of \(\End\Shatbundle\). The connection \(D\) induces a connection (also denoted \(D\)) on the endomorphism bundle via the Leibniz rule. Thus, to an \(\fr\)-symmetry \(a\in\fR\), we can assign the element \(Da\in\Omega^1(M;\End\Shatbundle)\) given by \(D_Xa:=\comm{D_X}{a} = \nablahat_X a - \comm{\betabold_X}{a}\) for \(X\in\fX(M)\) -- note that \(D_Xa\) need not be an \(\fr\)-symmetry.

\begin{lemma}\label{lemma:kvf-weak-condition}
	If \(X+a \in\fsymm(\varpihat,\eA)\) with \(X\in\fiso(M,g)\) and \(a\in\fR\) satisfies
	\begin{equation}\label{eq:kvf-preserve-ksf}
		\imath_X F - \eLhat_X\betabold - Da = 0
	\end{equation}
	then the action of that element on \(\fShat\) as described in Proposition~\ref{prop:symm-alg-spin-r} preserves \(\fShat_D\), and furthermore,
	\begin{equation}
		\eLhat_X R^D + \comm{a}{R^D} = 0.
	\end{equation}
	Moreover, if \(X+a, Y+b \in\fsymm(\varpihat,\eA)\) both satisfy equation~\eqref{eq:kvf-preserve-ksf}, then \(\ccomm{X+a}{Y+b}\) also satisfies that equation.
\end{lemma}

\begin{proof}
	The first claim follows essentially by construction: for all \(X\in\fiso(M,g)\), \(a\in\fR\) and \(\epsilon\in\fShat\),
	\begin{equation}
		\comm{\eLhat_X+a}{D}\epsilon 
		= \comm{\eLhat_X}{D}\epsilon + \comm{a}{D}\epsilon
		= (\imath_X F - \eLhat_X\betabold - Da)\epsilon,
	\end{equation}
	where we have used equation~\eqref{eq:Lhat-D-comm}; it immediately follows that if the quantity in brackets in the final expression vanishes, the operator \(\eLhat_X+a\) commutes with \(D\), whence it preserves \(\fShat_D\).

	Now note that \(\eLhat_X R = \eL_X R = 0\) identically for all \(X\in\fiso(M,g)\), so using the formula for \(R^D\) in Proposition~\ref{prop:D-curvature-rsymm} and the covariant Cartan calculus (Proposition~\ref{prop:cov-cartan-calculus}) as well as the Bianchi identity for \(F\), we have
	\begin{equation}
	\begin{split}
		\eLhat_X R^D
		&= \cancel{\eLhat_X R} + \eLhat_X F - \eLhat_X\dhat\betabold + [\eLhat_X\betabold\wedge\betabold]	\\
		&= \dhat\imath_X F + \cancel{\imath_X\dhat F} - \dhat\eLhat_X\betabold - [\imath_XF\wedge\betabold] + \qty[\eLhat_X\betabold\wedge\betabold]	\\
		&= \dhat\qty(\imath_XF - \eLhat_X\betabold) - \qty[\betabold\wedge \qty(\imath_XF - \eLhat_X\betabold)],
	\end{split}
	\end{equation}
	while on the other hand, noting that \(Da=\dhat a-\comm{\betabold}{a}\), we have
	\begin{equation}
	\begin{split}
		\comm{a}{R^D}
			&= \comm{a}{\cancel{R} + F + \tfrac{1}{2}[\betabold\wedge\betabold] - \dhat\betabold} \\
			&= \comm{a}{F} + [\comm{a}{\betabold}\wedge\betabold] - \dhat\comm{a}{\betabold} + \qty[(\dhat a)\wedge \betabold]	\\
			&= \comm{a}{F} - \dhat\comm{a}{\betabold} + \qty[(Da)\wedge\betabold]  \\
			&= \comm{a}{F} + \dhat\qty(\dhat a - Da) + \qty[(Da)\wedge\betabold]  \\
			&= - \dhat(Da) + \qty[\betabold\wedge (Da)],
	\end{split}
	\end{equation}
	where in the last line we have used that \(\dhat^2=F\wedge\) as an operator on forms with values in bundles associated to the spin-\(R\) structure. Thus if identity \eqref{eq:kvf-preserve-ksf} holds, we have \(\eLhat_X R^D + \comm{a}{R^D} = 0\).
	
	Using Cartan calculus and the Bianchi identity again, we can compute
	\begin{equation}
	\begin{split}
		\imath_{\comm{X}{Y}} F
			&= \eLhat_X\imath_Y F  - \imath_Y\eLhat_X F	\\
			&= \eLhat_X\imath_Y F  - \imath_Y\dhat\imath_X F - \cancel{\imath_Y\imath_X\dhat F}	\\
			&= \eLhat_X\imath_Y F - \eLhat_Y\imath_X F + \dhat\imath_Y\imath_X F,
	\end{split}
	\end{equation}
	for all Killing vectors \(X,Y\). 
	Then, using equation~\eqref{eq:eLhat-comm-forms}, we have
	\begin{equation}
	\begin{split}
		\imath_{\comm{X}{Y}} F - \eLhat_{\comm{X}{Y}}\betabold
			&= \eLhat_X\qty(\imath_Y F - \eLhat_Y\betabold) 
				- \eLhat_Y\qty(\imath_X F - \eLhat_X\betabold) 
				+ \dhat\imath_Y\imath_X F + \comm{F(X,Y)}{\betabold}	\\
			&= \eLhat_X\qty(\imath_Y F - \eLhat_Y\betabold) 
				- \eLhat_Y\qty(\imath_X F - \eLhat_X\betabold) 
				+ D(F(X,Y)).	\\
	\end{split}
	\end{equation}
	Then, assuming \(X+a\) and \(Y+b\) satisfy equation~\eqref{eq:kvf-preserve-ksf}, or equivalently \(\comm{\eLhat_X}{D}=Da\) and \(\comm{\eLhat_Y}{D}=Db\), we find
	\begin{equation}
	\begin{split}
		\imath_{\comm{X}{Y}} F - \eLhat_{\comm{X}{Y}}\betabold
			&= \eLhat_X\qty(Db) - \eLhat_Y\qty(Da) + D(F(X,Y))	\\
			&= D(\eLhat_Xb - \eLhat_Ya) + \comm{Da}{b} - \comm{Db}{a} + D(F(X,Y))	\\
			&= D(\eLhat_Xb - \eLhat_Ya + \comm{a}{b} + F(X,Y)).
	\end{split}
	\end{equation}
	This shows in particular that
	\begin{equation}
		\ccomm{X+a}{Y+b} = \comm{X}{Y} + (F(X,Y) + \eLhat_Xb - \eLhat_Ya + \comm{a}{b})
	\end{equation}
	satisfies \eqref{eq:kvf-preserve-ksf}, so the subspace of the symmetry algebra defined by that equation is closed under the bracket \(\ccomm{-}{-}\).
\end{proof}

We immediately have the following.

\begin{corollary}\label{coro:fVhat-D}
	The subspace of \(\fsymm(\varpihat,\eA)\) defined by equation~\eqref{eq:kvf-preserve-ksf}
	\begin{equation}
		\fVhat_D = \qty{ X+a\in\fsymm(\varpihat,\eA) \, \middle| \, X\in\fiso(M,g), a\in\fR \text{ s.t. } \imath_X F - \eLhat_X\betabold - Da = 0}
	\end{equation}
	is a Lie subalgebra which preserves \(\fShat_D\) as a subspace of \(\fShat\).
\end{corollary}

\begin{lemma}\label{lemma:fR-D}
	The space
	\begin{equation}
		\fR_D := \qty{a\in\fR \, | \, Da =0 }
	\end{equation}
	is a Lie subalgebra of \(\fR\) and an ideal in \(\fVhat_D\) with dimension at most \(\dim\fr\). Moreover, there is a short exact sequence of Lie algebras
	\begin{equation}\label{eq:fV-D-SES}
	\begin{tikzcd}
		0 \ar[r]	& \fR_D	\ar[r] 
					& \fVhat_D \ar[r]
					& \fV_D = \qty{ X\in\fiso(M,g) \, \middle| \, \exists a\in\fR \text{ s.t. } X+a\in\fVhat_D }	\ar[r] & 0
	\end{tikzcd}
	\end{equation}
	embedding into the sequence \eqref{eq:symm-alg-SES}. That is, \(\fVhat_D\) is an extension of a subalgebra of the algebra of Killing vector fields by \(\fR_D\). In particular, \(\fVhat_D\) is finite-dimensional with \(\dim\fVhat_D\leq\dim\fr+\dim\fiso(M,g)\).
\end{lemma}

\begin{proof}
	The Leibniz rule for \(D\) shows that \(\fR_D\) is a subalgebra of \(\fR\) and thus of \(\fsymm(\varpihat,\eA)\); setting \(X=0\) in equation~\eqref{eq:kvf-preserve-ksf} gives us \(Da=0\), showing that \(\fR_D\) is contained in the subalgebra \(\fVhat_D\), in particular \(\fR_D=\fVhat_D\cap\fR\). To see that it is an ideal in \(\fVhat_D\), let \(X+a\in\fVhat_D\), \(b\in\fR_D\) and note that \eqref{eq:Lhat-D-comm} holds when acting on \(\fR\subseteq\End(S)\) (where the RHS acts via the commutator), so we have
	\begin{equation}
		D \qty(\eLhat_X b + \comm{a}{b}) 
			= \eLhat_X D b - \comm{\imath_XF - \eLhat_X\betabold}{b} + D\comm{a}{b}
			= - \comm{Da}{b} + D\comm{a}{b}
			= \comm{a}{Db}
			= 0,
	\end{equation}
	whence \(\ccomm{X+a}{b}=\eLhat_X b + \comm{a}{b}\in\fR_D\) as required. Since elements of \(\fR_D\) are sections which are parallel with respect to a connection, they are determined by their value at a point; these values lie in a fibre which is isomorphic to \(\fr\), whence \(\dim\fR_D\leq \dim\fr\).
	
	The remaining claims follow by noting that \(\fV_D\) is the image of \(\fVhat_D\) under the restriction of the projection \(\fsymm(\varpihat,\eA)\to\fiso(M,g)\) and \(\fR_D=\fVhat_D\cap\fR\) is the kernel.
\end{proof}

\begin{remark}\label{rem:better-KSA-def}	 
	The author's PhD thesis \cite[Ch.4]{Beckett2024_thesis} and earlier drafts of this work took a different, perhaps more naïve approach to constructing the even part of the (super)algebra, essentially seeking subalgebras of \(\fiso(M,g)\) and \(\fR\) which separately preserve \(\fShat_D\) and whose sum is closed under the bracket \(\ccomm{-}{-}\).
	Starting from \eqref{eq:Lhat-D-comm}, one tries to impose the equation \(\imath_X F - \eLhat_X\betabold = 0\) (and thus \(Da=0\)) in the place of \eqref{eq:kvf-preserve-ksf}, but condition is not preserved by the Lie bracket of vector fields -- one must impose the stronger condition \(\imath_X F = \eLhat_X\betabold = 0\). While this does allow one to define a notion of a KSA generalising that in \cite{deMedeiros2018}, it has a number of pathologies which are avoided with the approach taken in the present work.
\end{remark}

We now attempt to extend the bracket \(\ccomm{-}{-}\) to include the space of parallel spinors \(\fShat_D\). We define a \(\ZZ_2\)-grading on \(\fVhat_D\oplus\fShat_D\) by declaring \(\fVhat_D\) to be the even subspace and \(\fShat_D\) the odd subspace, and we seek to extend \(\ccomm{-}{-}\) to a bracket on this space which preserves this grading. We will of course take
\begin{equation}
	\ccomm{X+a}{\epsilon} := \eLhat_X\epsilon + a\epsilon,
\end{equation}
for \(X\in\fiso(M,g)\) and \(a\in\fR\) such that \(X+a\in\fVhat_D\) and \(\epsilon\in\fShat_D\). The bracket of two spinors is less straightforward. While the bundle Dirac current \(\kappabold:\Odot^2\Shatbundle\to TM\) or \(\kappabold:\Wedge^2\Shatbundle\to TM\) constructed in Lemma~\ref{lemma:bundle-dirac-current-exist-rsymm} gives us a natural way of producing a vector field from a pair of spinor fields, we have no obvious way of obtaining an \(\fr\)-symmetry, so let us define
\begin{equation}\label{eq:spinor-spinor-bracket}
	\ccomm{\epsilon}{\zeta} := \kappabold(\epsilon,\zeta) + \rhobold(\epsilon,\zeta)
\end{equation}
for \(\epsilon,\zeta\in\fShat_D\), where \(\rhobold\) is some section of the bundle \(\Odot^2\Shatbundle^* \otimes_M \ad\Qbar\to M\) (if \(\kappabold\) is symmetric) or \(\Wedge^2\Shatbundle^* \otimes_M \ad\Qbar\to M\) (if \(\kappabold\) is skew-symmetric) which we leave undetermined for now.

The first thing to be checked is that the image of this bracket is in \(\fVhat_D\). In particular, \(\kappabold(\epsilon,\zeta)\) must be a Killing vector, and we must have \(\imath_{\kappabold(\epsilon,\zeta)} F - \eLhat_{\kappabold(\epsilon,\zeta)}\betabold - D(\rhobold(\epsilon,\zeta)) = 0\). Assuming for now that this is the case, we have thus defined a pre-Lie (super)algebra structure \(\ccomm{-}{-}\) on \(\fVhat_D\oplus\fShat_D\), that is, a bracket satisfying the axioms of a Lie bracket except for the Jacobi identity.
The components of the Jacobi identity with three entries from \(\fVhat_D\) or two from \(\fVhat_D\) and one from \(\fShat_D\) are satisfied by Corollary~\ref{coro:fVhat-D}.
This leaves the components with two or three entries from \(\fShat_D\) to be checked. 

First, for \(X+a \in\fVhat_D\) and \(\epsilon,\zeta\in\fShat_D\), we have
\begin{equation}\label{eq:RKSA-VSS-Jacobi}
\begin{split}
	\ccomm{X+a}{\ccomm{\epsilon}{\zeta}} - &\ccomm{\ccomm{X+a}{\epsilon}}{\zeta} - \ccomm{\epsilon}{\ccomm{X+a}{\zeta}}	\\
		& = \eL_X(\kappabold(\epsilon,\zeta)) + F(X,\kappabold(\epsilon,\zeta)) + \eLhat_X(\rhobold(\epsilon,\zeta))\\
		& \qquad 
		- \kappabold\qty(\eLhat_X\epsilon,\zeta) - \rhobold\qty(\eLhat_X\epsilon,\zeta)
		- \kappabold\qty(\epsilon,\eLhat_X\zeta) - \rhobold\qty(\epsilon,\eLhat_X\zeta)	\\
		& \qquad
		-\eLhat_{\kappabold(\epsilon,\zeta)} a + \comm{a}{\rhobold(\epsilon,\zeta)}
		- \kappabold(a\epsilon,\zeta) - \rhobold(a\epsilon,\zeta)
		- \kappabold(\epsilon,a\zeta) - \rhobold(\epsilon,a\zeta)	\\
		& = \qty(\eLhat_X\rhobold + a\cdot \rhobold)(\epsilon,\zeta) + F(X,\kappabold(\epsilon,\zeta)) -\eLhat_{\kappabold(\epsilon,\zeta)} a,
\end{split}
\end{equation}
where we have used \(\eLhat_X\kappabold = 0\) and the \(\fR\)-invariance of \(\kappabold\), and we have
\begin{equation}
	(a\cdot\rhobold)(\epsilon,\zeta) = \comm{a}{\rhobold(\epsilon,\zeta)} - \rhobold(a\epsilon,\zeta) - \rhobold(\epsilon,a\zeta)
\end{equation}
for all \(\epsilon,\zeta\in\fShat\). Now, using the fact that \(X+a\in\fVhat_D\), we can rewrite the last two terms as follows:
\begin{equation}
\begin{split}
	F(X,\kappabold(\epsilon,\zeta)) -\eLhat_{\kappabold(\epsilon,\zeta)} a
	= \imath_{\kappabold(\epsilon,\zeta)}(\imath_XF - \nablahat a)
	= \imath_{\kappabold(\epsilon,\zeta)}(\eLhat_X\betabold + Da - \nablahat a)
	= \imath_{\kappabold(\epsilon,\zeta)}(\eLhat_X\betabold - \comm{\betabold}{a}).
\end{split}
\end{equation}
For \(\epsilon,\zeta,\vartheta\in\fShat_D\),
\begin{multline}\label{eq:RKSA-SSS-Jacobi}
	\ccomm{\ccomm{\epsilon}{\zeta}}{\vartheta} 
		+ \ccomm{\ccomm{\zeta}{\vartheta}}{\epsilon} 
		+ \ccomm{\ccomm{\vartheta}{\epsilon}}{\zeta}\\
	= \eLhat_{\kappabold(\epsilon,\zeta)}\vartheta 
		+ \eLhat_{\kappabold(\zeta,\vartheta)}\epsilon 
		+ \eLhat_{\kappabold(\vartheta,\epsilon)}\zeta
		+ \rhobold(\epsilon,\zeta)\vartheta 
		+ \rhobold(\zeta,\vartheta)\epsilon 
		+ \rhobold(\vartheta,\epsilon)\zeta.
\end{multline}
Thus the Jacobi identity is satisfied if and only if
\begin{gather}
	\qty(\qty(\eLhat_X\betabold + \comm{a}{\betabold})\circ\kappabold + \qty(\eLhat_X\rhobold + a\cdot \rhobold))(\epsilon,\zeta) = 0,	\label{eq:jacobi-Xa-rho}\\
	\eLhat_{\kappabold(\epsilon,\zeta)}\vartheta 
		 + \eLhat_{\kappabold(\zeta,\vartheta)}\epsilon 
		 + \eLhat_{\kappabold(\vartheta,\epsilon)}\zeta
		 + \rhobold(\epsilon,\zeta)\vartheta 
		 + \rhobold(\zeta,\vartheta)\epsilon 
		 + \rhobold(\vartheta,\epsilon)\zeta = 0. \label{eq:Jacobi-3spinors-rsymm}
\end{gather}
for \(\epsilon,\zeta,\vartheta\in\fShat_D\), \(X+a\in\fVhat_D\). The first equation here suggests that we should further restrict our even subspace to those elements \(X+a\) that preserve \(\betabold\) and \(\rhobold\), i.e. those satisfying \(\eLhat_X\betabold + \comm{a}{\betabold} = 0\) and \(\eLhat_X\rhobold + a\cdot \rhobold = 0\). These conditions are natural if we consider \(\betabold\) and \(\rhobold\) to be part of the ``background data'' which we demand that our ``symmetries'', i.e. the even part of the Killing (super)algebra, must preserve. Under these conditions, we note that the defining equation~\eqref{eq:kvf-preserve-ksf} for \(\fVhat_D\) becomes
\begin{equation}
	\imath_X F - \nablahat a = 0.
\end{equation}

\begin{lemma}\label{lemma:fVhat-fR-D-rho}
	Corollary~\ref{coro:fVhat-D} and Lemma~\ref{lemma:fR-D} still hold if we replace \(\fVhat_D\), \(\fR_D\) and \(\fV_D\) by their respective subalgebras
	\begin{equation}\label{eq:fVhat-fR-D-rho}
	\begin{aligned}
		\fVhat_{(D,\rhobold)} 
			&= \qty{ X+a \in\fsymm(\varpihat,\eA)	\, \middle| \, 
				\begin{gathered}
					X\in\fiso(M,g), a\in\fR \text{ s.t. }\\
					\imath_X F = \nablahat a,\,
					\eLhat_X\betabold + \comm{a}{\betabold} = 0,\,
					\eLhat_X\rhobold + a\cdot\rhobold = 0,
				\end{gathered}
			},
			\\
		\fR_{(D,\rhobold)} 
			&= \qty{ a\in\fR	\, \middle| \, 
				\nablahat a = 0, \comm{a}{\betabold} = 0, a\cdot\rhobold = 0},
			\\
		\fV_{(D,\rhobold)} 
			&= \qty{ X\in\fiso(M,g)	\, \middle| \, 
				\exists a\in\fR \text{ s.t. } X+a\in\fVhat_{(D,\rhobold)} }.
	\end{aligned}
	\end{equation}
\end{lemma}

\begin{proof}
	By the above discussion, \(\fVhat_{(D,\rhobold)}\) is a subspace of \(\fVhat_D\); we need only show that the conditions \(\eLhat_X\betabold + \comm{a}{\betabold} = 0\) and \(\eLhat_X\rhobold + a\cdot \rhobold = 0\) are preserved by the Lie bracket. 
	Suppose we have \(X+a,Y+b\in\fsymm(\varpihat,\eA)\) such that \(\eLhat_X\betabold=-\comm{a}{\betabold}\) and \(\eLhat_Y\betabold=-\comm{b}{\betabold}\).
	Then, using equation~\eqref{eq:eLhat-comm-forms}, we have
	\begin{equation}
	\begin{split}
		\eLhat_{\comm{X}{Y}}\betabold
			&= \eLhat_X\eLhat_Y\betabold - \eLhat_Y\eLhat_X\betabold - \comm{F(X,Y)}{\betabold}	\\
			&= -\eLhat_X\comm{b}{\betabold} + \eLhat_Y\comm{a}{\betabold} - \comm{F(X,Y)}{\betabold}	\\
			&= -\comm{\eLhat_Xb}{\betabold} - \comm{b}{\eLhat_X\betabold} + \comm{\eLhat_Ya}{\betabold} + \comm{a}{\eLhat_Y\betabold} - \comm{F(X,Y)}{\betabold}	\\
			&= -\comm{F(X,Y) + \eLhat_Xb - \eLhat_Ya}{\betabold} + \comm{b}{\comm{a}{\betabold}} - \comm{a}{\comm{b}{\betabold}}	\\
			&= -\comm{F(X,Y) + \eLhat_Xb - \eLhat_Ya + \comm{a}{b}}{\betabold},
	\end{split}
	\end{equation}
	which shows that \(\ccomm{X+a}{Y+b}=\comm{X}{Y}+(F(X,Y) + \eLhat_Xb - \eLhat_Ya + \comm{a}{b})\) satisfies the desired equation for \(\betabold\). One can demonstrate this in a completely analogous way for the condition on \(\rhobold\), so \(\fVhat_{(D,\rhobold)}\) is indeed a subalgebra of \(\fVhat_D\). Now, \(\fR_{(D,\rhobold)}\) is a subspace of \(\fR_D\) since \(Da=\nablahat a - \comm{\betabold}{a}\) for all \(a\in\fR\) -- in fact, \(\fR_{(D,\rhobold)}=\fVhat_{(D,\rhobold)}\cap\fR_D\), so it is an ideal in \(\fVhat_{(D,\rhobold)}\) since \(\fR_D\) is an ideal in \(\fVhat_D\), and it is the kernel of the projection \(\fVhat_{(D,\rhobold)}\to\fV_{(D,\rhobold)}\).
\end{proof}

For the following result, we define \(\gammabold\in\Gamma\qty(\Hom\qty(\Otimes^2\Shatbundle,\End(TM)))\) by
\begin{equation}\label{eq:gamma-geom-def-rsymm}
	\gammabold(\epsilon,\zeta)X := - \kappabold(\betabold(X)\epsilon,\zeta) - \kappabold(\epsilon,\betabold(X)\zeta)
\end{equation}
for \(\epsilon,\zeta\in\fShat\) and \(X\in\fX(M)\). This map has the same symmetry as \(\kappabold\) in its spinorial arguments. Recalling that for \(X\in\fX(M)\) we define \(A_X\in\Gamma(\End(TM))\) by \(A_XY=-\nabla_Y X\), we note that
\begin{equation}
	\gammabold(\epsilon,\zeta)=A_{\kappabold(\epsilon,\zeta)}
\end{equation}
for all \(\epsilon,\zeta\in\fShat_D\) \cite[Lem.3.4]{Beckett2024_ksa}.

\begin{theorem}[Existence of (super)algebra associated to \((D,\rhobold)\)]\label{thm:killing-algebra-exist-rsymm}
	Let \(S\) be a (possibly \(N\)-extended) spinor module of \(\Spin(s,t)\), \(\kappa\) a symmetric (resp. skew-symmetric) Dirac current on \(S\) and \(R\) the corresponding \(R\)-symmetry group. Let \((M,g)\) be a connected, oriented and (in indefinite signature) time-oriented pseudo-Riemannian manifold of signature \((s,t)\) and suppose that it is equipped with a spin-\(R\) structure \(\varpihat:\Phat\to F_{SO}\). Recall that we have an associated (twisted) spinor bundle \(\Shatbundle=\Phat\times_{\sigmahat}S\to M\) and canonical \(\Rbar = R/\ZZ_2\)-bundle \(\Qbar=\Phat\times_{\pi_R}\Rbar\to M\). We let \(\fShat=\Gamma(\Shatbundle)\) and \(\fR=\Gamma(\ad\Qbar)\).
	
	Let \(\eA\) be a torsion-free principal connection on \(\Phat\to M\) and let \(\nablahat\) and \(\eLhat\) denote the covariant derivative and Lie derivative respectively. Let \(D\) be a connection on \(\Shatbundle\) and define the \(\End\Shatbundle\)-valued 1-form \(\betabold=\nablahat-D\).
	Denote the symmetric (resp. skew-symmetric) bundle Dirac current induced on \(\Shatbundle\) by \(\kappabold\) and suppose we have some section \(\rhobold \in \Gamma(\Odot^2\Shatbundle^* \otimes_M \ad\Qbar)\) (resp. \(\rhobold \in \Gamma(\Wedge^2\Shatbundle^* \otimes_M \ad\Qbar)\)). Then the \(\ZZ_2\)-graded vector space \(\fKhat_{(D,\rhobold)} = \fVhat_{(D,\rhobold)} \oplus\fShat_D\), where
	\begin{gather}
		\fKhat_{(D,\rhobold)\overline{0}} = \fVhat_{(D,\rhobold)} 
			= \qty{X+a \in \fX(M)\oplus\fR \,\middle| \, 
			\begin{gathered}
				\eL_Xg = 0,\,
				\imath_X F = \nablahat a,	\\
				\eLhat_X\betabold + \comm{a}{\betabold} = 0,\,
				\eLhat_X\rhobold + a\cdot \rhobold = 0,
			\end{gathered}
			},	\\
		\fKhat_{(D,\rhobold)\overline{1}} = \fShat_D
			= \qty{\epsilon\in\fShat \, \middle| \, D\epsilon = 0},
	\end{gather}
	 equipped with the bracket
	\begin{align*}
		\comm{X+a}{Y+b}
			&= \eL_X Y + F(X,Y) + \nablahat_Xb - \nablahat_Ya + \comm{a}{b}_{\fR}	\\ 
			&= \eL_X Y - F(X,Y) + \comm{a}{b}_{\fR},	\\
		\comm{X+a}{\epsilon} &= \eLhat_X \epsilon + a \epsilon,	\\
		\comm{\epsilon}{\zeta} &= \kappabold(\epsilon,\zeta) + \rhobold(\epsilon,\zeta),
	\end{align*}
	where \(X+a,Y+b\in\fVhat_{(D,\rhobold)} \) and \(\epsilon,\zeta\in\fShat_D\), is a Lie superalgebra (resp. \(\ZZ_2\)-graded Lie algebra) if and only if the following conditions are satisfied for all \(\epsilon,\zeta,\vartheta\in\fShat_D\):
	\begin{gather}
		\gammabold(\epsilon,\zeta) \in \fso(M,g),
			\label{eq:killing-algebra-exist-rsymm-1} \\
		\eLhat_{\kappabold(\epsilon,\zeta)}\betabold + \comm{\rhobold(\epsilon,\zeta)}{\betabold} = 0,
			\label{eq:killing-algebra-exist-rsymm-2beta}	\\
		\eLhat_{\kappabold(\epsilon,\zeta)}\rhobold + \rhobold(\epsilon,\zeta)\cdot\rhobold = 0,
			\label{eq:killing-algebra-exist-rsymm-2rho}	\\
		\imath_{\kappabold(\epsilon,\zeta)}F = \nablahat(\rhobold(\epsilon,\zeta)),	 
			\label{eq:killing-algebra-exist-rsymm-2F}	\\
		\betabold_{\kappabold(\epsilon,\zeta)}\vartheta + \gammabold(\epsilon,\zeta)\cdot\vartheta + \rhobold(\epsilon,\zeta)\vartheta + \text{cyclic perm's} = 0, \label{eq:killing-algebra-exist-rsymm-3}
	\end{gather}
	where \(\gammabold\) is the map defined in equation \eqref{eq:gamma-geom-def-rsymm}.
\end{theorem}

\begin{proof}
	By the preceding discussion, it remains only to check closure of the  \(\fKhat_{(D,\rhobold)\overline{1}}\otimes \fKhat_{(D,\rhobold)\overline{1}} \to \fKhat_{(D,\rhobold)\overline{0}}\) component of the bracket and the all-odd component of the Jacobi identity, equation~\eqref{eq:Jacobi-3spinors-rsymm}, which is the only component not satisfied by construction.
	For closure, let \(\epsilon,\zeta\in\fShat_D\); we must show that \(\kappabold(\epsilon,\zeta) + \rhobold(\epsilon,\zeta) \in\fVhat_{(D,\rhobold)}\).
	Recall that \(\kappabold(\epsilon,\zeta)\) is a Killing vector if and only if the endomorphism \(A_{\kappabold(\epsilon,\zeta)}=\gammabold(\epsilon,\zeta)\) is skew-symmetric with respect to the metric, which is just \eqref{eq:killing-algebra-exist-rsymm-1}. The remaining conditions are exactly \eqref{eq:killing-algebra-exist-rsymm-2beta}-\eqref{eq:killing-algebra-exist-rsymm-2F}.
	For the Jacobi identity, letting \(\epsilon,\zeta,\vartheta\in\fShat_D\) and using the fact that \(\nablahat_{\kappabold(\epsilon,\zeta)}\vartheta=\betabold_{\kappabold(\epsilon,\zeta)}\vartheta\) and \(A_{\kappabold(\epsilon,\zeta)}=\gammabold(\epsilon,\zeta)\), we can write
	\begin{equation}
		\eLhat_{\kappabold(\epsilon,\zeta)}\vartheta 
		= \nablahat_{\kappabold(\epsilon,\zeta)}\vartheta + A_{\kappabold(\epsilon,\zeta)}\cot\vartheta 
		= \betabold_{\kappabold(\epsilon,\zeta)}\vartheta + \gammabold(\epsilon,\zeta)\cdot\vartheta,
	\end{equation}
	so equation~\eqref{eq:Jacobi-3spinors-rsymm} is equivalent to \eqref{eq:killing-algebra-exist-rsymm-3}.
\end{proof}

This result motivates the following definition, which is a generalisation of the definition of admissible connections and Killing (super)algebras for spin manifolds \cite[Def.3.6]{Beckett2024_ksa}.

\begin{definition}[Admissible pairs, Killing spinors and Killing (super)algebras]\label{def:killing-spinor-rsymm}
	Using the notation of Theorem~\ref{thm:killing-algebra-exist-rsymm}, we say that \((D,\rhobold)\) is an \emph{admissible pair} if the following conditions are satisfied:
	\begin{enumerate}
		\item The map \(\gammabold\) defined by \eqref{eq:gamma-geom-def-rsymm} takes values \(\ad F_{SO}\), i.e. \(\gammabold(\epsilon,\zeta)\in\fso(M,g)\) for all \(\epsilon,\zeta\in\fShat\); \label{item:admiss-pair-killing}
		\item Equation~\eqref{eq:killing-algebra-exist-rsymm-3} holds for all \(\epsilon,\zeta,\vartheta\in\fShat\); \label{item:admiss-pair-jacobi}
		\item Equations~\eqref{eq:killing-algebra-exist-rsymm-2beta}-\eqref{eq:killing-algebra-exist-rsymm-2F} hold for all \(\epsilon,\zeta\in\fShat_D\). \label{item:admiss-pair-bg}
	\end{enumerate}
If \((D,\rhobold)\) is an admissible pair, we will refer to the differential equation \(D\epsilon=0\) (equivalently \(\nablahat\epsilon =\betabold\epsilon\)) as the \emph{Killing spinor equation}. We say that \(\fShat_D\) is the space of \emph{Killing spinors} and \(\fKhat_{(D,\rhobold)}\) is the \emph{Killing superalgebra} (resp. \emph{Killing algebra}) associated to the pair.
\end{definition}

\begin{remark}\label{rem:better-KSA-def-2}
	The definition presented here also generalises (in a number of directions) the one given in \cite[\S6.3, Sec.7]{deMedeiros2018}, which treats the special case of Lorentzian supergravity backgrounds in 6 dimensions with reducible spin-\(h\) structure with trivial \(\Sp(1)\)-bundle and flat \(\Sp(1)\) connection.
	Indeed, the cited sections of \textit{loc. cit.} provide two alternative constructions of the Killing (super)algebra: the first takes a global geometric perspective as we have here but essentially assumes from the start that the connection on the \(R\)-bundle is flat; the second takes a local perspective and does not assume flatness \textit{a priori} but nonetheless arrives at it after imposing further constraints, which amount to assuming that \(\fVhat_D\) is a sum of subalgebras of \(\fiso(M,g)\) and \(\fR\) (see Remark~\ref{rem:better-KSA-def}).
\end{remark}

\subsection{Algebraic structure of the Killing (super)algebras}
\label{sec:alg-structure-killing-rsymm}

Generalising the 11-dimensional Lorentzian case \cite{Figueroa-OFarrill2017_1}, our previous work \cite[Thm.3.10]{Beckett2024_ksa} established that Killing (super)algebras on spin manifolds are filtered subdeformations of a flat model (super)algebra \(\fs\). The goal of this section is to prove the obvious generalisation for spin-\(R\) manifolds.

\begin{theorem}[Structure of Killing (super)algebras with \(R\)-symmetry]\label{thm:killing-algebra-filtered-rsymm}
	Using the notation of Theorem~\ref{thm:killing-algebra-exist-rsymm} and Definitions \ref{def:flat-model-alg-rsymm} and \ref{def:killing-spinor-rsymm}, let \((D,\rhobold)\) be an admissible pair on \(\Shatbundle\to M\). Then the associated Killing (super)algebra \(\fKhat_{(D,\rhobold)}\) is a filtered subdeformation of the \(\fr\)-extended flat model (super)algebra \(\fshat\) associated to \((\RR^{s,t},S,\kappa)\).
\end{theorem}

\subsubsection{Localising at a point}
\label{sec:localise-ksa-at-point-rsymm}

Consider the underlying space of the \(\fr\)-extended flat model (super)algebra \(\fshat=\RR^{s,t}\oplus S\oplus(\fso(s,t)\oplus\fr)\) as a \(\Spin^R(s,t)\)-module and let us define the associated vector bundle \(\eEhat:=\Phat\times_{\pihat\oplus\sigmahat\oplus\Ad}\fshat\). Note that
\begin{equation}\label{eq:eE-def-rsymm}
	\eEhat \cong TM\oplus_M\Shatbundle\oplus_M\ad\Phat
		\cong TM\oplus_M\Shatbundle\oplus_M\ad F_{SO}\oplus_M\ad\Qbar,
\end{equation}
so its space of sections is \(\Gamma(\eEhat) = \fX(M)\oplus\fShat\oplus\fso(M,g)\oplus\fR\). We define a connection \(\eDhat\) on \(\eEhat\) as follows:
\begin{equation}\label{eq:eDhat-def}
	\eDhat_Y(X,\epsilon,A,a) := 
		\qty(\nabla_Y X + AY,\, D_Y\epsilon,\, \nabla_Y A + R(Y,X),\, \nablahat_Ya + F(Y,X))
\end{equation}
for \(X,Y\in\fX(M)\), \(\epsilon\in\fShat\), \(A\in\fso(M,g)\), \(a\in\fR\). This is a generalisation of the Killing transport connection used in \cite[\S3.2.1]{Beckett2024_ksa} and \cite{Figueroa-OFarrill2017_1} (which itself extends the classical notion of Kostant and Geroch \cite{Kostant1955,Geroch1969}). We will use it to express the structure of \(\fKhat_{(D,\rhobold)}\) in terms of data over a point \(x\in M\).\footnote{
	Recalling \eqref{eq:kvf-preserve-ksf}, one might be tempted to define the fourth component of \(\eDhat_Y(X,\epsilon,A,a)\) as \(D_Ya + \imath_Y\eLhat_X\betabold - F(X,Y)\) so that parallel sections correspond to elements of \(\fVhat_D\oplus\fShat_D\). However, this expression lies not in \(\fR\) but in \(\End\fShat\) in general, so \(\eDhat\) would not be a connection on \(\eEhat\).}

\begin{proposition}\label{prop:eDhat-parallel}
	The \(\eDhat\)-parallel sections of \(\eEhat\) are tuples \((X,\epsilon,A_X,a)\), where \(X\) is a Killing vector, \(A_X=-\nabla X\), \(\nablahat a = \imath_XF\) and \(\epsilon\in\fShat_D\).
\end{proposition}

\begin{proof}
	A section \((X,\epsilon,A,a)\) of \(\eEhat\) is \(\eDhat\)-parallel if and only if
	\begin{equation}
	\begin{aligned}
		\nabla_Y X + AY & = 0,
		& D_Y \epsilon &= 0,
		& \nabla_YA &= R(X,Y),
		& \nablahat_Ya &= F(X,Y),
	\end{aligned}
	\end{equation}
	for all \(Y\in\fX(M)\). The first equation tells us that \(A=-\nabla X = A_X\), whence \(X\) is a Killing vector since \(A\in\fso(M,g)\); the third equation is implied by the first (see e.g. \cite[Lem.2.2]{Kostant1955}). The second and fourth equations precisely say that \(\epsilon\in\fShat_D\) and \(\nablahat a = \imath_XF\). 
\end{proof}

We immediately see that any \(X+a+\epsilon\in\fKhat_{(D,\rhobold)}\) with \(X+a\in\fKhat_{(D,\rhobold)\overline{0}}=\fVhat_{(D,\rhobold)}\), \(\epsilon\in\fKhat_{(D,\rhobold)\overline{1}}=\fShat_D\) defines a \(\eDhat\)-parallel section \((X,\epsilon,A_X,a)\) of \(\eEhat\). By parallel transport, any such section is determined by its value \((X_x,\epsilon_x,(A_X)_x,a_x)\) at any point \(x\in M\), which we call the \emph{(Killing) transport data of \(X+a+\epsilon\) at \(x\)}. 
We can thus localise the brackets of \(\fKhat_{(D,\rhobold)}\) at \(x\) -- that is, we can express the transport data of the bracket of a pair of elements of \(\fKhat_{(D,\rhobold)}\) in terms of the transport data of the elements of the pair. We let \(\tau_x: \fKhat_{(D,\rhobold)}\to\eEhat_x\) be the \(\RR\)-linear map sending an element of \(\fKhat_{(D,\rhobold)}\) to its transport data at \(x\):
\begin{equation}
	\tau_x(X+a+\epsilon):=(X_x,\epsilon_x,(A_X)_x,a_x),
\end{equation}
and we can then compute the following for \(X+a,Y+b\in\fVhat_{(D,\rhobold)}\) and \(\epsilon,\zeta\in\fShat_D\):
\begin{align}
	\tau_x(\comm{X+a}{Y+b})
		&= \qty((A_X)_x Y_x - (A_Y)_x X_x,
			0,
			\begin{aligned}
				\comm{(A_X)_x}{(A_Y)_x}_{\fso(T_xM)} \quad \\ - R_x(X_x,Y_x)
			\end{aligned},
			\begin{aligned}
				\comm{a_x}{b_x}_{\ad\Qbar_x} \qquad \\ -F_x(X_x,Y_x)
			\end{aligned}
			),
		\label{eq:phihat-x-XY}	\\
	\tau_x(\comm{X+a}{\epsilon}) 
		&= (0,\betabold_x(X_x)\epsilon_x + (A_X)_x\cdot\epsilon_x + a_x\epsilon_x,0,0),
		\label{eq:phihat-x-Xepsilon}	\\
	\tau_x(\comm{\epsilon}{\zeta}) 
		&= (\kappabold_x(\epsilon_x,\zeta_x),0,\gammabold_x(\epsilon_x,\zeta_x) = - \kappabold_x(\betabold_x\epsilon_x,\zeta_x) - \kappabold_x(\epsilon_x,\betabold_x\zeta_x),\rhobold_x(\epsilon_x,\zeta_x)),
		\label{eq:phihat-x-epsilonzeta}
\end{align}
where \(R_x\in\Wedge^2T^*_xM\otimes\fso(T_xM,g_x)\) is the value of the Riemann tensor, \(F_x\in\Wedge^2T^*_xM\otimes\ad\Qbar_x\) is the value of the field strength tensor, and \(\betabold_x\in T^*_xM\otimes\End\Shatbundle_x\), \(\kappabold_x\in\Hom(\Odot^2S,T_xM)\) (resp. \(\kappabold_x\in\Hom(\Wedge^2S,T_xM)\)) and \(\rhobold_x\in\Hom(\Odot^2S,\ad\Qbar_x)\) (resp. \(\rhobold_x\in\Hom(\Wedge^2S,\ad\Qbar_x)\)) are the values of the maps \(\betabold,\kappabold\) and \(\rhobold\) at \(x\). Note that all of the operations on the right-hand side are defined entirely in terms of data at the point \(x\).

\subsubsection{Proof of the structure theorem}

\begin{proof}[Proof of Theorem \ref{thm:killing-algebra-filtered-rsymm}]
	We fix a basepoint \(x\in M\) and define \((V,\eta):=(T_xM,g_x)\) and, by a convenient abuse of notation, \(S=\Shatbundle_x\) and \(\fr=\ad\Qbar_x\). Indeed, recalling the isomorphism \eqref{eq:eE-def-rsymm}, the fibre \(\eEhat_x=V\oplus \Shatbundle_x\oplus(\fso(V)\oplus\ad\Qbar_x)\) is naturally an \(\fr\)-extended flat model algebra with Dirac current \(\kappabold_x\); in fact, by choosing a basepoint \(p\in\Phat_x\), it can be identified with the typical fibre \(\fshat\), resolving the notational ambiguity. Let us then define the ``evaluation'' maps
	\begin{equation}
	\begin{aligned}
		\evaluate^V_x&: 
			& \fVhat_{(D,\rhobold)}  &\longrightarrow V,
			& X+a &\longmapsto X_x, \\
		\evaluate^S_x&:
			& \fShat_D &\longrightarrow S,
			&\epsilon &\longmapsto \epsilon_x,
	\end{aligned}
	\end{equation}
	and the related spaces
	\begin{equation}
	\begin{aligned}
		\fhhat &= \qty{(A_X)_x + a_x \in\fso(V)\oplus\fr \,\middle |\, X+a\in\fVhat_{(D,\rhobold)}:\,X_x=0 },	\\
		V' &= \Im\evaluate^V_x = \qty{X_x\in V \,\middle |\, X+a\in\fVhat_{(D,\rhobold)} },	\\
		S' 	&= \Im\evaluate^S_x = \qty{\epsilon_x\in S \,\middle |\, \epsilon\in\fShat_D }.
	\end{aligned}
	\end{equation}
	There is natural linear map \(\ker\evaluate^V_x\to\fhhat\) given by \(X+a\mapsto (A_X)_x+a_x\); this is an isomorphism since \(\ker\evaluate^V_x\) is precisely the set of elements of \(\fVhat_{(D,\rhobold)}\) whose vectorial part vanishes at \(x\), so the inverse is given by \(\eDhat\)-transport.
	We can also identify \(S'\) with \(\fShat_D\) via \(D\)-transport.
	
	Now note that \(\fa:=V'\oplus S'\oplus\fhhat\) is a \(\ZZ\)-graded subspace of \(\eEhat_x\cong \fshat\), and one can see by inspecting transport data that \(\fa\) is preserved by the flat model bracket, hence it is in fact a graded subalgebra.\footnote{
		A similar argument with more detail is given in the proof of \cite[Thm.3.10]{Beckett2024_ksa}.} 
	 Indeed, \(\fhhat\) is a subalgebra of \(\fso(V)\) preserving \(V'\), which can be seen by inspecting \eqref{eq:phihat-x-XY}; it also preserves \(S\) by \eqref{eq:phihat-x-Xepsilon}. 
	 We also have \(\comm{S'}{S'}=\kappabold_x(S',S')\subseteq V'\) by \eqref{eq:phihat-x-epsilonzeta}, and \(\comm{V'}{V'}=0\) trivially.
	 
	 Finally, we show that \(\fKhat_{(D,\rhobold)}\) is filtered with associated graded (super)algebra \(\fa\). The filtration \(\fKhat_{(D,\rhobold)}^\bullet\) is as follows: \(\fKhat_{(D,\rhobold)}^i = \fKhat_{(D,\rhobold)}\) for \(i<-2\), \(\fKhat_{(D,\rhobold)}^i=0\) for \(i>0\), and
	 \begin{equation}
 	\begin{aligned}
 		&\fKhat_{(D,\rhobold)}^{-2} = \fKhat_{(D,\rhobold)} = \fVhat_{(D,\rhobold)} \oplus \fShat_D
 		&&\supseteq 
 		&&\fKhat_{(D,\rhobold)}^{-1} = \ker\evaluate^V_x\oplus \fShat_D
 		&&\supseteq
 		&&\fKhat_{(D,\rhobold)}^0 = \ker\evaluate^V_x,
 	\end{aligned}
	\end{equation}
	and it is easily checked that \(\comm{\fKhat_{(D,\rhobold)}^i}{\fKhat_{(D,\rhobold)}^j}\subseteq\fKhat_{(D,\rhobold)}^{i+j}\) for all \(i,j\in\ZZ\).
	Recall that the associated \(\ZZ\)-graded (super)algebra \(\Gr\fKhat_{(D,\rhobold)}^\bullet\) has components \(\Gr_i\fKhat_{(D,\rhobold)}^\bullet={\fKhat_{(D,\rhobold)}^i}/{\fKhat_{(D,\rhobold)}^{i+1}}\); we will denote the projection \(\fKhat_{(D,\rhobold)}^i\to\Gr_i\fKhat_{(D,\rhobold)}^\bullet\) by an overline. 
	The graded bracket \(\comm{-}{-}_{\Gr}\) on \(\Gr\fKhat_{(D,\rhobold)}^\bullet\) is defined such that the component 
	\(\Gr_i\fKhat_{(D,\rhobold)}^\bullet\otimes\Gr_j\fKhat_{(D,\rhobold)}^\bullet\to\Gr_{i+j}\fKhat_{(D,\rhobold)}^\bullet\) 
	is given by
	\(\comm{\overline{A}}{\overline{B}}_{\Gr} = \overline{\comm{A}{B}}\)
	 for \(A\in\fKhat_{(D,\rhobold)}^i\) and \(B\in\fKhat_{(D,\rhobold)}^j\).
	 We immediately have \(\Gr_{i}\fKhat_{(D,\rhobold)}^\bullet=0\) for \(i<-2\) or \(i>0\). For \(i=0\) we have an isomorphism of Lie algebras
	\begin{equation}
		\Gr_0 \fKhat_{(D,\rhobold)}^\bullet = \fKhat_{(D,\rhobold)}^0 = \ker\evaluate^V_x\cong \fa_0 = \fhhat,
	\end{equation}
	whence we also obtain 	
 	\begin{equation}
 	\begin{aligned}
 		\Gr_{-2}\fKhat_{(D,\rhobold)}^\bullet &\cong \faktor{\fVhat_{(D,\rhobold)}}{\ker\evaluate^V_x}		&&\cong \fa_{-2} = V'	&& \text{(using \(\evaluate^V_x\))},	\\
 		\Gr_{-1}\fKhat_{(D,\rhobold)}^\bullet &\cong \fShat_D	&& \cong \fa_{-1} = S'	&& \text{(using \(\evaluate^{S}_x\))},	\\
 	\end{aligned}
 	\end{equation}
 	as \(\Gr_0 \fKhat_{(D,\rhobold)}^\bullet\cong\fa_0\)-modules, thus \(\Gr\fKhat_{(D,\rhobold)}^\bullet\cong\fa\) as \(\ZZ\)-graded \(\fa_0\)-modules. To show that this is in fact an isomorphism of \(\ZZ\)-graded Lie (super)algebras, it remains only to check that brackets which do not involve \(\Gr_0 \fKhat_{(D,\rhobold)}^\bullet\) are preserved. 
 	The components \(\Gr_{-2}\fKhat_{(D,\rhobold)}^\bullet\otimes\Gr_{-2}\fKhat_{(D,\rhobold)}^\bullet\to\Gr_{-4}\fKhat_{(D,\rhobold)}^\bullet\)
 	and
 	\(\Gr_{-2}\fKhat_{(D,\rhobold)}^\bullet\otimes\Gr_{-1}\fKhat_{(D,\rhobold)}^\bullet\to\Gr_{-3}\fKhat_{(D,\rhobold)}^\bullet\)
 	are trivially zero, so only the
 	\(\Gr_{-1}\fKhat_{(D,\rhobold)}^\bullet\otimes\Gr_{-1}\fKhat_{(D,\rhobold)}^\bullet\to\Gr_{-2}\fKhat_{(D,\rhobold)}^\bullet\)
 	component now remains; for \(\epsilon,\zeta\in\fShat_D\), we have
	\begin{align}
		\comm{\overline{\epsilon}}{\overline{\zeta}}_{\Gr} = \overline{\comm{\epsilon}{\zeta}} 
			= \overline{\kappabold(\epsilon,\zeta)+\rhobold(\epsilon,\zeta)}
			\longmapsto \kappabold_x(\epsilon_x,\zeta_x) 
			= \comm{\epsilon_x}{\zeta_x}_{\fa},
	\end{align}
	under the \(\fa_0\)-module isomorphism \(\Gr_{-2}\fKhat_{(D,\rhobold)}^\bullet\xrightarrow{\cong}\fa_{-2}=V\) induced by \(\evaluate^V_x\), as required.
\end{proof}

Continuing to use the notation of the proof, we can explicitly describe the Lie (super)algebra structure of \(\fKhat_{(D,\rhobold)}\) as a deformation of the bracket on the graded subalgebra \(\fa=V'\oplus S'\oplus\fhhat\) of \(\eEhat_x\cong\fshat\) as follows. The map \(\evaluate^V_x\) gives us a short exact sequence of vector spaces
\begin{equation}
\begin{tikzcd}
	0 \ar[r] & \fhhat \ar[r] & \fVhat_{(D,\rhobold)} \ar[r,"\evaluate^V_x"] & V' \ar[r] & 0
\end{tikzcd}
\end{equation}
which can be split by a choice of linear map \(\Lambda:V'\to \fVhat_{(D,\rhobold)}\) with \(\Lambda_1(v)_x = v\), where we write \(\Lambda(v)=\Lambda_1(v)+\Lambda_2(v)\) as an element of \(\fiso(M,g)\oplus\fR\). 
This induces a linear isomorphism \(V'\oplus \fhhat\to \fVhat_{(D,\rhobold)}\) given by \(v+(A+a)\mapsto \Lambda(v) + \widetilde{A+a}\), where \(\widetilde{A+a}\) is the element of \(\fVhat_{(D,\rhobold)}\) with transport data \((0,0,A,a)\). 
The inverse of this map is \(X+a \mapsto X_x + ((A_{X-\Lambda_1(X_x)})_x + (a - \Lambda_2(X_x))_x)\). 
We also define the map \(\lambda: V'\to \fso(V)\oplus\fr\) as the sum of the maps \(\lambda_1:V'\to\fso(V)\), \(\lambda_2:V'\to\fr\) with \(\lambda_1(v)=(A_{\Lambda_1(v)})_x\), \(\lambda_2(v)=\Lambda_2(v)_x\).
This contains the same information as \(\Lambda\) since \(\Lambda(v)\) can be reconstructed from its transport data \((v,0,\lambda_1(v),\lambda_2(v))\). 
Note that \(\lambda\) (or \(\Lambda\)) is unique up to addition of a map \(V'\to\fhhat\). 
We then define a vector space isomorphism \(\Psihat:\fa\to\fKhat_{(D,\rhobold)}\) by
\begin{equation}
	\Psihat(v+\epsilon+(A+a)) = \Lambda(v) + \widetilde{A+a} + \widetilde{\epsilon},
\end{equation}
for \(v\in V'\), \(\epsilon\in S'\), and \(A\in\fso(V)\), \(a\in\fr\) such that \(A+a\in\fhhat\); we use \(\widetilde{\epsilon}\in\fShat_D\) to denote the spinor field with value \(\epsilon\) at \(x\); \(\widetilde\epsilon_x=\epsilon\). The inverse is given by
\begin{equation}
	\Psihat^{-1}(X + a + \epsilon) 
		= X_x + \epsilon_x + \qty((A_X)_x-\lambda_1(X_x) + a_x-\lambda_2(X_x)),
\end{equation}
where \(X\in\fiso(M,g)\), \(a\in\fR\) such that \(X+a\in\fVhat_{(D,\rhobold)}\), and \(\epsilon\in\fShat_D\).
We then define a new bracket \(\comm{-}{-}'\) on \(\fa\) by using \(\Psihat\) to pull back the bracket of \(\fKhat_{(D,\rhobold)}\):
\begin{equation}
\begin{split}
	\Psihat\qty(\comm{v+\epsilon+(A+a)}{w+\zeta+(B+b)}') 
		&= \comm{\Psihat(v+\epsilon+(A+a))}{\Psihat(w+\zeta+(B+b))} \\
		&= \comm{\Lambda(v) + \widetilde{A+a} + \widetilde\epsilon}{\Lambda(w) + \widetilde{B+b} + \widetilde\zeta},
\end{split}
\end{equation}
so that by construction, \(\Psihat\) is a Lie algebra isomorphism \((\fa,\comm{-}{-}')\cong\fKhat_{(D,\rhobold)}\); inverting it gives us
\begin{equation}\label{eq:KSA-deformed-brackets-rsymm}
\begin{gathered}
	\begin{aligned}
		\comm{A+a}{B+b}' &= \comm{A}{B} + \comm{a}{b},\\
		\comm{A+a}{v}' &= Av + \comm{A}{\lambda_1(v)}-\lambda_1(Av) +\comm{a}{\lambda_2(v)}-\lambda_2(Av),	\\
		\comm{A+a}{\epsilon}	&= A\cdot\epsilon + a\cdot \epsilon,	\\
		\comm{v}{w}' &= \lambda_1(v)w - \lambda_1(w)v + \theta_1(v,w) + \theta_2(v,w),	\\
		\comm{v}{\epsilon}' &= \betabold_x(v)\epsilon + \lambda(v)\cdot \epsilon,\\
		\comm{\epsilon}{\zeta}'
			&= \kappabold_x(\epsilon,\zeta)
			+ \gammabold_x(\epsilon,\zeta) + \rhobold_x(\epsilon,\zeta)
			-\lambda\qty(\kappabold_x(\epsilon,\zeta)),	\\
	\end{aligned}
\end{gathered}
\end{equation}
where the \(\theta_i\)s are the \(\fso(V)\)- and \(\fr\)-components of a map \(\theta:\Wedge^2V'\to \fhhat\) given by
\begin{align}
	\theta_1(v,w) 
		&= \comm{\lambda_1(v)}{\lambda_1(w)} - R_x(v,w) - \lambda_1(\lambda_1(v)w-\lambda_1(w)v),	
			\label{eq:KSA-deformed-brackets-theta1}\\	
	\theta_2(v,w) 
		&= \comm{\lambda_2(v)}{\lambda_2(w)} - F_x(v,w) - \lambda_2(\lambda_1(v)w-\lambda_1(w)v).
			\label{eq:KSA-deformed-brackets-theta2}
\end{align}
We thus have
\begin{equation}
	\comm{-}{-}' = \comm{-}{-} + \betabold_x + \gammabold_x + \rhobold_x + \partial\lambda + \theta,
\end{equation}
where \(\betabold_x\), \(\gammabold_x\), \(\rhobold_x\)  and \(\theta\) have been trivially extended to maps \(\Otimes^2\fshat\to\fshat\). The map \(\mu = \betabold_x + \gammabold_x + \rhobold_x + \partial\lambda\) has degree \(+2\) with respect to the grading on \(\fa\), while \(\theta\) has degree \(+4\). The deformation therefore has defining sequence \((\mu,\theta,0,\dots)\).

\section{Filtered subdeformations of the \(\fr\)-extended Poincaré superalgebra}
\label{sec:filtered-def-flat-model-rsymm}

Now that we have established that the Killing (super)algebras are filtered subdeformations of an \(\fr\)-extended flat model (super)algebra \(\fshat\), we will develop a theory of such deformations in much the same way as we did for filtered subdeformations of \(\fs\) in \cite[Sec.4]{Beckett2024_ksa}. We specialise to the case of \(\fr\)-extended flat model \emph{super}algebras -- that is, we now assume that the Dirac current \(\kappa\) is symmetric -- in order to simplify the exposition. In particular, we exploit the fact that, as explained in \cite[\S2.2.3]{Beckett2024_ksa}, symmetric multilinear maps \(\Psi:\Odot^n U\to W\) are determined by their ``polarised'' values \(\Psi(u,u,\dots,u)\) for \(u\in U\). From \S\ref{sec:highly-susy-subalg-rsymm} onward, we will also assume Lorentzian signature and place further conditions on the Dirac current in order to make use of the Homogeneity Theorem (Thm.~\ref{thm:homogeneity}) and better control the \(R\)-symmetry. We make extensive use of the terminology and results for filtered and graded Lie superalgebras and the Spencer cohomology of the latter presented in \cite[\S2.2,\S{}A.3]{Beckett2024_ksa} (see also the references therein, in particular \cite{Cheng1998}). The approach here follows that of \cite[Sec.4]{Beckett2024_ksa} and \cite{Figueroa-OFarrill2017_1} but requires some more careful analysis.

We consider graded subalgebras of the \(\fr\)-extended flat model superalgebra \(\fshat\)
\begin{equation}\label{eq:graded-subalg-form}
	\fa = V'\oplus S'\oplus \fhhat
\end{equation}
where \(\fa_{-1}=S'\) is a vector subspace of \(S\), \(\fa_{-2}=V'\) is a vector subspace of \(V\) such that \(\kappa_s\in V'\) for all \(s\in S'\), where here and in the sequel we write \(\kappa_s:=\kappa(s,s)\) for a polarised Dirac current, and \(\fa_0=\fhhat\) is a subalgebra of \(\fshat_0=\fso(V)\oplus\fr\) whose action preserves \(S'\) and \(V'\).

\begin{remark}\label{rem:better-KSA-def-3}
	In the author's PhD thesis \cite[Ch.4]{Beckett2024_thesis} and earlier drafts of this work, only those graded subalgebras with \(\fhhat=\fh\oplus\fr'\) where \(\fh\) is a subalgebra of \(\fso(V)\) and \(\fr'\) is a subalgebra of \(\fr\) were considered. A similar restriction was applied in \cite{deMedeiros2018}. However, this restriction is incompatible with our definition of the even part of the Killing superalgebra for spin-\(R\) manifolds (see Remarks~\ref{rem:better-KSA-def} and \ref{rem:better-KSA-def-2}), and is in any case somewhat unnatural, so we consider the general case here.
\end{remark}

While \(\fhhat\) may not be a sum of subalgebras of \(\fso(V)\) and \(\fr\), it is in fact an \emph{extension} of the former by the latter, reflecting the fact that \(\fshat_{\overline 0}\) is an extension of \(\fs_{\overline 0}\) by \(\fr\).

\begin{lemma}\label{lemma:hhat-ext}
	Let \(\fh=\{A\in\fso(V)\,|\,\exists a\in\fr \text{ s.t. } A+a\in\fhhat\}\). Then there are embeddings of short exact sequences of Lie algebras such that the following diagram with exact rows commutes:
	\begin{equation}
	\begin{tikzcd}
		& 0 \ar[r]	& \fhhat\cap\fr	\ar[r] \ar[d,equal]
					& \fhhat		\ar[r] \ar[d,hook]
					& \fh			\ar[r] \ar[d,hook]	& 0	\\
		& 0 \ar[r]	& \fhhat\cap\fr	\ar[r] \ar[d,hook]
					& \fa_{\overline 0} = \fhhat\ltimes V'	\ar[r] \ar[d,hook]
					& \fh\ltimes V'		\ar[r] \ar[d,hook]	& 0 \\
		& 0 \ar[r]	& \fr				\ar[r] 
					& \fshat_{\overline 0} = \fs_{\overline 0}\oplus\fr	\ar[r]
					& \fs_{\overline 0} = \fso(V)\ltimes V	\ar[r] & 0.
	\end{tikzcd}
	\end{equation}
\end{lemma}

\begin{proof}
	Clearly \(\fh\) is the image of \(\fhhat\) under the projection \(\fso(V)\oplus\fr\to\fso(V)\) (which is a Lie algebra morphism), and the kernel of the restriction of that map to \(\fhhat\) is \(\fhhat\cap\fr\). The rest of the diagram follows by considering the projection \(\fshat_{\overline 0} = \fs_{\overline 0}\oplus\fr \to \fs_{\overline 0}\) and its restrictions.
\end{proof}

Note that we could also express \(\fhhat\) as an extension of a subalgebra of \(\fr\) by one of \(\fso(V)\) by using the projection \(\fso(V)\oplus\fr\to\fr\), but this description is less useful because, unlike \(\fr\), \(\fso(V)\) is not an ideal in \(\fshat_{\overline 0}\), so the latter is not an extension by the former, and similarly, the subalgebra \(\fhhat\cap\fso(V)\) is not an ideal in \(\fa_{\overline 0}\) in general.
Likewise, the subalgebras \(\fr\subseteq\fshat\) and \(\fhhat\cap\fr\subseteq\fa\) are not ideals since \(\comm{\fr}{\fshat_{\overline 1}}\subseteq\fshat_{\overline 1}\), so the diagram in the lemma above does not extend to include the superalgebras.

\subsection{The Spencer (2,2)-cohomology}
\label{sec:poincare-spencer-22-cohomology-rysmm}

As we will discuss in more detail later in this section, following work by Cheng and Kac \cite{Cheng1998}, filtered deformations of graded subalgebras of \(\fshat\) are governed by the Spencer cohomology of those subalgebras, particularly the degree-2 part. Here, we will describe the degree-2 complex and give a useful characterisation of the cohomology group \(\ssH^{2,2}(\fshat_-;\fshat)\) of the full \(\fr\)-extended flat model algebra.

We recall that the Spencer complex for a \(\ZZ\)-graded Lie superalgebra \(\fa=\bigoplus_{i=-h}^\infty\fa_i\) (where the \emph{depth} \(h\in\ZZ_{\geq 0}\) is finite and \(\dim\fa_i<\infty\)) is the (bi-graded) Chevalley--Eilenberg complex \((\ssC^{\bullet,\bullet}(\fa_-;\fa),\partial)\) for the subalgebra \(\fa_-=\bigoplus_{i=-h}^{-1}\fa_i\) with values in \(\fa\).
In more detail: \(\fa_-\) acts on \(\fa\) via \(\ad^{\fa}|_{\fa_-}\); the cochain space of homological degree \(p\) inherits a \(\ZZ\)-grading \(\ssC^p(\fa_-;\fa)=\bigoplus_{d\in\ZZ}\ssC^{d,p}(\fa_-;\fa)\) where \(\ssC^{d,p}(\fa_-;\fa)=\Hom(\Wedge^p\fa_-;\fa)_d\) is the space of degree-\(d\) linear maps with respect to the \(\ZZ\)-grading on  \(\fa\); and the differential \(\partial\) respects this grading, so for each \(d\in\ZZ\) we get a sub-complex \((\ssC^{d,\bullet}(\fa_-;\fa),\partial)\).
Moreover, since the adjoint action of \(\fa_0\) preserves the grading on \(\fa\), so each \(\fa_i\) is an \(\fa_0\)-submodule, the cochain spaces inherit an action of \(\fa_0\) which preserves all gradings and the differential.
Thus \(\fa_0\) also preserves the spaces of cocycles \(\ssZ^{d,p}(\fa_-;\fa)\) and coboundaries \(\ssB^{d,p}(\fa_-;\fa)\), so the cohomology groups \(\ssH^{d,p}(\fa_-;\fa) := \ssZ^{d,p}(\fa_-;\fa)/\ssB^{d,p}(\fa_-;\fa)\) inherit in turn an action of \(\fa_0\).
We will be particularly interested in the spaces of invariants \(\ssH^{2,2}(\fa_-;\fa)^{\fa_0}\) and \(\ssH^{2,2}(\fshat_-;\fshat)^{\fa_0}\) for graded subalgebras \(\fa\) of \(\fshat\). For more details on generalities, including an explicit formula for the differential \(\partial\), see \cite{Cheng1998} or \cite[\S2.2]{Beckett2024_ksa}.

\subsubsection{The complex}

The degree-2 Spencer complex for a graded subalgebra \(\fa=V'\oplus S'\oplus\fhhat\) of \(\fshat\) is
\begin{equation}
\begin{split}
	0
	&\longrightarrow \ssC^{2,1}(\fa_-;\fa) = \Hom(V',\fhhat)\\
	&\longrightarrow \ssC^{2,2}(\fa_-;\fa) = \Hom(\Wedge^2V',V')\oplus\Hom(V'\otimes S',S')
		\oplus\Hom(\Odot^2S',\fhhat)\\
	&\longrightarrow \ssC^{2,3}(\fa_-;\fa) = \Hom(V'\otimes\Odot^2S',V')\oplus\Hom(\Odot^3S',S')\\
	&\longrightarrow 0,
\end{split}
\end{equation}
where we note that \(\fa_-=V'\oplus S'\). We denote the projections to the components of \(\ssC^{2,2}(\fa_-;\fa)\) by
\begin{equation}
\begin{aligned}
	\pi_1&: \ssC^{2,2}(\fa_-;\fa) \longrightarrow \Hom(\Wedge^2V',V'),\\
	\pi_2&: \ssC^{2,2}(\fa_-;\fa) \longrightarrow \Hom(V'\otimes S',S'),\\
	\pi_3&: \ssC^{2,2}(\fa_-;\fa) \longrightarrow \Hom(\Odot^2S',\fhhat).
\end{aligned}
\end{equation}

Let us recall the following result, which we will use in several places.

\begin{lemma}[{\cite[Lem.4.1,Coro.4.2]{Beckett2024_ksa}}]\label{lemma:V-soV-maps}
	If \(\lambda:V\to\fso(V)\) is a linear map such that \(\lambda(v)w - \lambda(w)v = 0\) for all \(v,w\in V\) then \(\lambda=0\). In particular, there is an isomorphism \(\Hom(V,\fso(V))\cong\Hom(\Wedge^2V,V)\) given by \(\lambda\mapsto\alpha_\lambda\), where \(\alpha_\lambda(v,w)=\lambda(v)w-\lambda(w)v\).
\end{lemma}

We then have the following technical lemma which generalises a much more straightforward analogue for graded subalgebras of \(\fs\) \cite[Coro.4.2]{Beckett2024_ksa}. It already hints that ``highly supersymmetric'' graded subalgebras \(\fa\) of \(\fshat\) (see \S\ref{sec:homogeneity-faithful-rsymm} and also \cite[\S4.3.1]{Beckett2024_ksa}) as well as those for which \(\fhhat\) acts faithfully on \(S'\) will be particularly amenable to our homological methods.

\begin{lemma}\label{lemma:proj-diff-injective-rsymm}
	Let \(\fa=V'\oplus S'\oplus\fhhat\) be a graded subalgebra of \(\fshat\).
	\begin{enumerate}
		\item \label{item:proj-diff-injective-1}
			If \(V'=V\) then the kernel of the map
			\begin{equation}
				\pi_1\circ\partial:\Hom(V,\fhhat)\longrightarrow\Hom(\Wedge^2V,V)
			\end{equation}
			is \(\Hom(V,\fhhat\cap\fr)\), so we have an injective map \(\Hom(V,\fh)\to\Hom(\Wedge^2V,V)\) where \(\fh\) is the subalgebra of \(\fso(V)\) appearing in Lemma~\ref{lemma:hhat-ext}. If \(\fh=\fso(V)\), the latter map is the isomorphism in Lemma~\ref{lemma:V-soV-maps}.
		\item \label{item:proj-diff-injective-2}
			If \(\fhhat\) acts faithfully on \(S'\) then
			\begin{equation}
				\pi_2\circ\partial:
					\Hom(V',\fhhat)\longrightarrow\Hom(V'\otimes S',S')
			\end{equation}
			is an injective map of \(\fhhat\)-modules.
		\item \label{item:proj-diff-injective-34}
			If \(V'=\kappa(\Odot^2S')\) then
			\begin{align}
				\pi_3\circ\partial:&
					\Hom(V',\fhhat)\longrightarrow\Hom(\Odot^2S',\fhhat)
			\end{align}
			(which is simply a pull-back by \(-\kappa\)) is an injective maps of \(\fhhat\)-modules.
		\item \label{item:proj-diff-injective-cohom}
			If \(\fhhat\) acts faithfully on \(S'\) or if \(V'=\kappa(\Odot^2S')\), then 
			\begin{equation}
				\ssH^{2,1}(\fa_-;\fa)=\ssZ^{2,1}(\fa_-;\fa)=0.
			\end{equation}
	\end{enumerate}
\end{lemma}

\begin{proof}
	Let \(\lambda\in\ssC^{2,1}(\fa_-;\fa)=\Hom(V',\fhhat)\) throughout. For part \ref{item:proj-diff-injective-1}, let us assume \(V'=V\). Then for \(v,w\in V\),
	\begin{equation}
		(\partial\lambda)(v,w) = \lambda_1(v)w-\lambda_1(w)v = 0,
	\end{equation}
	where \(\lambda_1:V\to\fso(V)\) is the \(\fso(V)\)-component of \(\lambda\), so \(\lambda_1=0\) by Lemma~\ref{lemma:V-soV-maps}, i.e. \(\lambda\) takes values in \(\fhhat\cap\fr\).
	Taking the quotient by the kernel gives the injective map claimed, where we recall that \(\fh\cong\fhhat/(\fhhat\cap\fr)\). A dimension count shows that this map is an isomorphism in the \(\fh=\fso(V)\) case, and it is easy to see that it is indeed the one appearing in Lemma~\ref{lemma:V-soV-maps}.
	
	Then if \(\fhhat\) acts faithfully on \(S'\) and
	\begin{equation}
		((\pi_2\circ\partial)\lambda)(v,s) = \lambda(v)\cdot s = 0
	\end{equation}
	for all \(v\in V'\) and \(s\in S'\), we must have \(\lambda=0\), which proves \ref{item:proj-diff-injective-2}. If \(V'=\kappa(\Odot^2S')\) and
	\begin{equation}
		((\pi_3\circ\partial)\lambda)(s,s) = -\lambda(\kappa_s) = 0
	\end{equation}
	for all \(s\in S'\), then we must have \(\lambda=0\), so we have proven \ref{item:proj-diff-injective-34}.
	Part \ref{item:proj-diff-injective-cohom} follows immediately.
\end{proof}

\subsubsection{The cohomology}

We will not attempt to calculate the Spencer cohomology of an arbitrary graded subalgebra \(\fa\) of \(\fshat\) nor describe its structure; however, we can derive a useful description of the Spencer cohomology group  \(\ssH^{2,2}(\fshat_-;\fshat)\) of \(\fshat\) itself, from which we will later bootstrap some homological information about certain classes of subalgebras.

We observe that since \(V\) is a simple \(\fso(V)\)-module and the Dirac current \(\kappa:\Odot^2S\to V\) is an \(\fshat_0=\fso(V)\oplus\fr\)-module morphism, \(\kappa\) is either zero or surjective; thus, if it is non-zero, we have an \(\fshat_0\)-module isomorphism \(V\cong\Odot^2S/\ker\kappa\). It will be useful for us to find an orthogonal complement to \(\ker\kappa\) in \(\Odot^2S\), giving us a splitting
\begin{equation}\label{eq:spinor-square-decomposition}
	\Odot^2S \cong V\oplus \ker\kappa.
\end{equation}
Such a splitting of \(\fso(V)\)-modules always exists since \(\fso(V)\) is simple; often \(\fr\) is semisimple (indeed, it is often also simple), whence \(\fshat_0\) is also semisimple, or the group \(R\) is compact, so this can be taken to be a splitting of \(\fshat_0\)-modules.
In practice, \(\kappa\) is typically constructed from a non-degenerate \(\fshat_0\)-invariant bilinear which can be employed to decompose \(\Odot^2 S\) (as well as \(\Wedge^2 S\)) into a product of factors of the form \(\Wedge^p V\otimes\Delta_p\), where \(\Delta_p\) is some auxiliary \(\fr\)-module, but we will not concern ourselves with such details here. 
Assuming we have a splitting of \(\fshat_0\)-modules, we denote the pull-backs of any (linear or \(\fshat_0\)-module) map \(\phi:\Odot^2 S\to W\) (where \(W\) is some other vector space or module) to the components of this splitting by \(\phi_V\) and \(\phi_{\ker\kappa}\) respectively.

\begin{lemma} \label{lemma:H22-full-flat-model-rsymm}
	Let \(\alpha+\beta+\gamma+\rho\) be a cocycle in \(\ssC^{2,2}(\fshat_-;\fshat)\) where
	\begin{equation}
	\begin{aligned}
		\alpha \in	&\Hom(\Wedge^2V,V),
		&\beta \in	&\Hom(V\otimes S,S),\\
		\gamma \in	&\Hom(\Odot^2S,\fso(V)),
		&\rho \in	&\Hom(\Odot^2S,\fr),
	\end{aligned}
	\end{equation}
	and suppose that we have a splitting \eqref{eq:spinor-square-decomposition} of \(\fshat_0\)-modules. 
	Then the cohomology class \([\alpha+\beta+\gamma+\rho]\in \ssH^{2,2}(\fshat_-;\fshat)\) has a unique representative \(\beta'+\gamma'+\rho'\) with zero \(\Hom(\Wedge^2 V,V)\) component and for which \(\rho'_V=0\). 
	We denote the space of such \emph{normalised cocycles} by \(\cHhat^{2,2}\):
	\begin{equation}
		\cHhat^{2,2} 
		= \qty{\beta+\gamma+\rho\in\ssZ^{2,2}(\fshat_-;\fshat)
			~ \middle| ~
			\begin{gathered}
				\beta\in\Hom(V\otimes S,S),\\
				\gamma\in\Hom(\Odot^2 S,\fso(V)),\\
				\rho\in\Hom(\Odot^2S,\fr),\,
				\rho_V=0
			\end{gathered}
		}
	\end{equation}
	and then we have \(\ssH^{2,2}(\fshat_-;\fshat) \cong \cHhat^{2,2}\) as \(\fshat_0\)-modules.
\end{lemma}

\begin{proof}
	By Lemma~\ref{lemma:V-soV-maps} there is a unique map \(\lambda_1: V\to \fso(V)\) such that \(\alpha=\partial\lambda_1|_{\Wedge^2 V}\). Moreover, using the decomposition of \(\Odot^2 S\) discussed above, we can write \(\rho_V=\partial\lambda_2|_{\Odot^2 S}\) for a unique map \(\lambda_2:V\to\fr\). Then we can define a cocycle \(\beta'+\gamma'+\rho'\in\ssZ^{2,2}(\fshat_-;\fshat)\) with 	\(\beta'\in\Hom(V\otimes S,S)\), \(\gamma'\in\Hom(\Odot^2S,\fso(V))\), \(\rho' \in\Hom(\Odot^2S,\fr)\), in particular with trivial \(\Hom(\Wedge^2 V,\fso(V))\)-component and \(\rho_V'=0\), by
	\begin{equation}
		\beta' + \gamma' + \rho' = \alpha + \beta + \gamma + \rho - \partial\lambda
	\end{equation}
	where \(\lambda=\lambda_1+\lambda_2\in\Hom(V,\fso(V)\oplus\fr)=\ssC^{2,1}(\fshat_-;\fshat)\), which by construction lies in \(\cHhat^{2,2}\). It follows that \(\ssZ^{2,2}(\fshat_-;\fshat)=\cHhat^{2,2}\oplus\ssB^{2,2}(\fshat_-;\fshat)\) as \(\fshat_0\)-modules, hence the final claim.
\end{proof}

The normalised \((2,2)\)-cocycles satisfy the following pair of \emph{normalised Spencer cocycle conditions}:
\begin{gather}
	2\kappa(s,\beta(v,s)) + \gamma(s,s)v = 0, \label{eq:cocycle-1-rsymm}	\\
	\beta(\kappa_s,s) + \gamma(s,s)\cdot s + \rho(s,s)s = 0, \label{eq:cocycle-2-rsymm}
\end{gather}
for \(v\in V\), \(s\in S\) (again writing \(\kappa_s=\kappa(s,s)\)); we must also separately impose the \(\rho_V=0\) condition. Note that condition \eqref{eq:cocycle-1-rsymm} fully determines the component \(\gamma\) if \(\beta\) is already known.

A normalised cocycle \(\beta+\gamma+\rho\) is invariant under the action of a subalgebra \(\fhhat\subseteq\fso(V)\oplus\fr\) if and only if \(\beta,\gamma\) and \(\rho\) are all separately invariant. Moreover, by the condition \eqref{eq:cocycle-1-rsymm}, invariance of \(\beta\) implies invariance of \(\gamma\), whence invariance of the cocycle is equivalent to invariance of \(\beta\) and \(\rho\).

We also note that if we set \(\rho=0\), the equations above are precisely the normalised cocycle conditions for the unextended flat model superalgebra \(\fs\) (compare to \cite[eq.(20),(21)]{Beckett2024_ksa}). Thus we have \(\cH^{2,2}=\cHhat^{2,2}\cap\{\rho=0\}\), where \(\cH^{2,2}\) is the space of normalised cocycles for \(\fs\) defined in \cite[Lem.4.4]{Beckett2024_ksa}.

\subsection{General filtered deformations and cohomology}
\label{sec:filtered-def-cohom-rsymm}

As shown in \cite{Cheng1998} and recapitulated in \cite[\S2.2.2]{Beckett2024_ksa}, a filtered deformation of a \(\ZZ\)-graded Lie superalgebra can be described (though non-canonically) as a deformed super-Lie bracket on the same underlying space, where the deformation is the sum of a ``defining sequence'' of maps of positive, even degree. In the case at hand, a filtered deformation \(\fatilde\) of a graded subalgebra \(\fa=V'\oplus S'\oplus\fhhat\) of \(\fshat\) has the brackets
\begin{equation}\label{eq:def-brackets-general-rsymm}
\begin{aligned}
	\comm{A+a}{B+b} &= AB-BA + ab-ba,	
	&\comm{v}{w} &= \alpha(v,w) + \theta(v,w),	\\
	\comm{A+a}{v} &= Av + \delta(A+a,v),
	&\comm{s}{s} &= \kappa_s + \gamma(s,s) + \rho(s,s),	\\
	\comm{A+a}{s} &= A\cdot s + as,
	&\comm{v}{s} &= \beta(v,s),	\\
\end{aligned}
\end{equation}
where \(A+a,B+b\in\fhhat\), \(v,w\in V'\), \(s\in S'\). The deformation maps of degree 2 are
\begin{align}
	\alpha&: \Wedge^2V' \longrightarrow V',
	&\beta&: V'\otimes S' \longrightarrow S',
	&\gammahat&: \Odot^2S'\longrightarrow \fhhat,
\end{align}
where we note that we have written \(\gammahat\) in terms of the maps\footnote{The discrepancy in notation for \(\gammahat\) and its components in comparison to other \(\fhhat\)-valued maps is for consistency with notation in other works \cite{deMedeiros2018,Beckett2024_thesis}, and because it is more convenient for some calculations.}
\begin{align}
	\gamma&: \Odot^2 S' \longrightarrow \fso(V),
	&\rho&: \Odot^2 S' \longrightarrow \fr
\end{align}
with \(\gammahat(s,s)=\gamma(s,s)+\rho(s,s)\in\fhhat\) for all \(s\in S'\), as well as a map
\begin{align}\label{eq:delta-comps-def}
	\delta&: \fhhat\otimes V' \to \fhhat
	&&\text{with components}
	&\delta_1&: \fhhat\otimes V'\longrightarrow\fso(V),
	&\delta_2&: \fhhat\otimes V'\longrightarrow\fr,
\end{align}
such that \(\delta(A+a,v)=\delta_1(A+a,v)+\delta_2(A+a,v)\in\fhhat\) for all \(A+a\in\fhhat\), \(v\in V'\). The degree-4 deformation map is
\begin{align}
	\theta&: \Wedge^2 V' \longrightarrow \fhhat
	&& \text{with components}
	&\theta_1&: \Wedge^2 V' \longrightarrow \fso(V),
	&\theta_2&: \Wedge^2 V' \longrightarrow \fr,
\end{align}
such that \(\theta(v,w)=\theta_1(v,w)+\theta_2(v,w)\in\fhhat\) for all \(v,w\in V\). We denote the full degree-2 deformation by
\begin{equation}
	\mu = \alpha + \beta + \gamma + \rho + \delta: \fa\otimes\fa\longrightarrow\fa,
\end{equation}
so the defining sequence for the deformation is \((\mu,\theta,0,\dots)\). 

The filtration \(\fatilde^\bullet\) is as follows: as vector spaces, \(\fatilde^i= \fa\) for \(i<-2\), \(\fatilde^i=0\) for \(i>0\), and
\begin{equation}\label{eq:atilde-filtration}
	\begin{aligned}
		&\fatilde^{-2} = \fa_{-2}\oplus\fa_{-1}\oplus\fa_0 = V'\oplus S'\oplus\fhhat
		&&\supseteq 
		&&\fatilde^{-1} = \fa_{-1}\oplus\fa_0 = S'\oplus\fhhat
		&&\supseteq
		&&\fatilde^0 = \fa_0 = \fhhat.
	\end{aligned}
\end{equation}
In fact, the equalities in the zeroth part of the filtration are also equalities of Lie algebras. 
Note that the adjoint actions of \(\fhhat\) on \(\fa\) and on \(\fatilde\) preserve the graded and filtered subspaces respectively, and that the natural maps \(\fatilde^i\to\Gr_i\fatilde^\bullet=\fatilde^{i}/\fatilde_{i+1}=\fa_i\) are \(\fhhat\)-equivariant.

Just as \(\fhhat\) and \(\fa_{\overline0}\) are extensions by \(\fhhat\cap\fr\) of subalgebras of \(\fs_0=\fso(V)\) and \(\fs_{\overline0}\) respectively (Lemma~\ref{lemma:hhat-ext}), so too is \(\fatilde_{\overline 0}\) an extension of a filtered subdeformation \(\overline\fa_{\overline 0}\) of \(\fs_{\overline0}\) by \(\fhhat\cap\fr\); letting \(\overline\fa_{\overline 0}=V'\oplus S'\oplus\fh\) as a vector space and endowing it with brackets \(\comm{A}{B}=AB-BA\), \(\comm{v}{w}=\alpha(v,w)+\theta_1(v,w)\) and \(\comm{A}{v}=Av+\delta_1(A)\) for \(A,B\in\fh\) and \(v,w\in V'\), there is a short exact sequence
\begin{equation}
\begin{tikzcd}
	0 \ar[r]	& \fhhat\cap\fr				\ar[r] 
				& \fatilde_{\overline 0}	\ar[r]
				& \overline\fa_{\overline 0} \ar[r] & 0
\end{tikzcd}
\end{equation}
obtained by projecting out \(\fr\)-components. 
Compare this to the similar extension diagrams \eqref{eq:symm-alg-SES} and \eqref{eq:fV-D-SES}, the latter of which is essentially an example of this one.

Following \cite{Cheng1998}, the deformation maps correspond to some homological data for \(\fa\). The relevant results from \textit{loc. cit.} are summarised in \cite[Prop.2.3]{Beckett2024_ksa}; for the case at hand, we obtain the following (see also the first part of Proposition~\ref{prop:highly-susy-filtered-def-rsymm}):
\begin{align}
	\mu &\in \ssZ^{2,2}(\fa;\fa),	\label{eq:mu-CE-cocycle-rsymm}	\\
	\mu|_{\fa_-\otimes\fa_-} &\in \ssZ^{2,2}(\fa_-;\fa),	\label{eq:mu-spencer-cocycle-rsymm}	\\
	\qty[\mu|_{\fa_-\otimes\fa_-}] &\in \ssH^{2,2}(\fa_-;\fa)^{\fa_0}. \label{eq:mu-invt-class-rsymm}
\end{align}
The first condition says that \(\mu\) is a degree-2 (with respect to the \(\ZZ\)-grading on \(\fa\)) Chevalley--Eilenberg 2-cocycle of \(\fa\) with values in the adjoint representation, and the second and third say that restricting \(\mu\) to \(\fa_-=V'\oplus S'\) gives one a Spencer cocycle with \(\fa_0\)-invariant cohomology class. We will examine these conditions explicitly below and see that the first condition actually implies the other two. In the case without \(\fr\)-extension, there was a converse to this implication for subalgebras with \(V'=V\) \cite[Lem.4.6]{Beckett2024_ksa}; in the \(\fr\)-extended case, the situation is more complicated, but we will nonetheless derive an analogous result (Lemma~\ref{lemma:invt-spencer-cocycle-rsymm}).

\subsubsection{Unpacking the cocycle conditions}

The condition \eqref{eq:mu-CE-cocycle-rsymm} is equivalent to the following system of equations:
\begin{align}
	\alpha(\kappa_s,v) + 2\kappa(s,\beta(v,s)) + \gamma(s,s)v &= 0,\label{eq:22-spencer-cocycle-vss-rsymm}\\
	\beta(\kappa_s,s) + \gamma(s,s)\cdot s + \rho(s,s)s &= 0,\label{eq:22-spencer-cocycle-sss-rsymm}
\end{align}
for all \(v\in V'\), \(s\in S'\);
\begin{align}
	A\alpha(v,w) - \alpha(Av,w) - \alpha(v,Aw) 
		- \delta_1(A+a,v)w + \delta_1(A+a,w)v 
		&= 0,\label{eq:alpha-delta1-rsymm}\\
	A\cdot(\beta(v,s)) - \beta(Av,s) - \beta(v,A\cdot s) 
		+ a(\beta(v,s)) - \beta(v,a s)
		- \delta_1(A+a,v)\cdot s - \delta_2(A+a,v)s 
		&= 0,\label{eq:beta-delta12-rsymm}\\
	\comm{A}{\gamma(s,s)} - 2\gamma(A\cdot s,s) 
		- 2\gamma(as,s)
		+ \delta_1(A+a,\kappa_s) 
		&= 0,\label{eq:gamma-delta1-rsymm}\\
	- 2\rho(A\cdot s,s)
		+ \comm{a}{\rho(s,s)} - 2\rho(as,s)
		+ \delta_2(A+a,\kappa_s) 
		&= 0,\label{eq:rho-delta2-rsymm}
\end{align}
for all \(A+a\in\fhhat\), \(v,w\in V'\), \(s\in S'\);
\begin{align}
	\delta_1(\comm{A}{B}+\comm{a}{b},v) 
		- \comm{A}{\delta_1(B+b,v)} + \comm{B}{\delta_1(A+a,v)}
		- \delta_1(A+a,Bv) + \delta_1(B+b,Av)
		&= 0, \label{eq:delta1-cocycle-cond-rsymm}\\
	\delta_2(\comm{A}{B}+\comm{a}{b},v) 
	- \comm{a}{\delta_2(B+b,v)} + \comm{b}{\delta_2(A+a,v)}
	- \delta_2(A+a,Bv) + \delta_2(B+b,Av) &= 0, \label{eq:delta2-cocycle-cond-rsymm}
\end{align}
for all \(A+a,B+b\in\fhhat\), \(v\in V'\). We note that all brackets here are commutators in \(\fso(V)\oplus\fr\).

Equations \eqref{eq:22-spencer-cocycle-vss-rsymm} and \eqref{eq:22-spencer-cocycle-sss-rsymm} are equivalent to condition \eqref{eq:mu-spencer-cocycle-rsymm} and as such are referred to as the \emph{Spencer cocycle conditions} for \(\fa\); compare these to the \emph{normalised} Spencer cocycle conditions \eqref{eq:cocycle-1-rsymm} and \eqref{eq:cocycle-2-rsymm} for \(\fshat\) itself.
Equations~\eqref{eq:alpha-delta1-rsymm}-\eqref{eq:rho-delta2-rsymm} are equivalent to 
\begin{equation}
	(A+a)\cdot(\alpha+\beta+\gamma+\rho) = \partial\imath_{A+a}(\delta_1+\delta_2),
\end{equation}
for all \(A+a\in\fhhat\), which implies the condition \eqref{eq:mu-invt-class-rsymm}. 
The remaining equations \eqref{eq:delta1-cocycle-cond-rsymm} and \eqref{eq:delta2-cocycle-cond-rsymm} can be interpreted as yet another cocycle condition (using the isomorphism of \(\fhhat\)-modules \(\Hom(\fhhat\otimes V',\fhhat)\cong\Hom(\fhhat,(V')^*\otimes\fhhat)\)):
\begin{equation}
	\delta\in\ssZ^1\qty(\fhhat;(V')^*\otimes\fhhat).
\end{equation}

The following analogue of a result for graded subalgebras of \(\fshat\) \cite[Lem.4.6]{Beckett2024_ksa}, but needs some additional regularity conditions so that we can make use of Lemma~\ref{lemma:proj-diff-injective-rsymm}. 
The upshot is that, under these conditions, a Spencer \((2,2)\)-cocycle with \(\fa_0\)-invariant cohomology class uniquely determines a cocycle in \(\ssZ^{2,2}(\fa;\fa)\).

\begin{lemma}\label{lemma:invt-spencer-cocycle-rsymm}
	If \(\mu\in\ssZ^{2,2}(\fa;\fa)\) then \(\mu|_{\fa_-\otimes\fa_-}\in\ssZ^{2,2}(\fa_-;\fa)\) and \([\mu|_{\fa_-\otimes\fa_-}]\in\ssH^{2,2}(\fa_-;\fa)^{\fa_0}\). If \(V=V'\) and either \(V=\kappa(\Odot^2S')\) or \(\fhhat\) acts faithfully on \(S'\) then, conversely, if \(\mu_-\in\ssZ^{2,2}(\fa_-;\fa)\) such that \([\mu_-]\in\ssH^{2,2}(\fa_-;\fa)^{\fa_0}\), then there exists a unique \(\delta\in\ssZ^1(\fhhat;V^*\otimes\fhhat)\) such that \(\mu=\mu_-+\delta\in\ssZ^{2,2}(\fa;\fa)\).
\end{lemma}

\begin{proof}
	The first claim is immediate from the discussion above. For the converse statement, we first note that by part~\ref{item:proj-diff-injective-cohom} of Lemma~\ref{lemma:proj-diff-injective-rsymm}, we have \(\ssZ^{2,1}(\fa_-;\fa)=0\); that is, \(\partial:\ssC^{2,1}(\fa_-;\fa)\to\ssC^{2,2}(\fa_-;\fa)\) is injective.	
	Now, invariance of the Spencer cohomology class means that for all \(X\in\fhhat\), 
	\begin{equation}\label{eq:mu-partial-chi-invariance}
		X\cdot\mu_- = \partial\chi_X
	\end{equation}
	for some \(\chi_X\in\ssC^{2,1}(\fa_-;\fa)=\Hom(V;\fhhat)\), but by injectivity of \(\partial\), \(\chi_X\) is unique and the assignment \(X\mapsto\chi_X\) linear, thus we can uniquely define a linear map \(\delta:\fhhat\otimes V\to\fhhat\) by \(\imath_X\delta=\chi_X\), i.e. \(\delta(X,v)=\chi_X(v)\). By construction, the components of \(\delta\) satisfy the equations \eqref{eq:alpha-delta1-rsymm}-\eqref{eq:rho-delta2-rsymm} (where of course \(\mu=\alpha+\beta+\gamma+\rho+\delta\) as usual); it remains only to check that \eqref{eq:delta1-cocycle-cond-rsymm} and \eqref{eq:delta2-cocycle-cond-rsymm} are satisfied, i.e. that \(\delta\in\ssZ^1(\fhhat;V^*\otimes\fhhat)\). Let us define a map \(\delta_{X,Y}\in\ssC^{2,1}(\fa_-;\fa)\) by
	\begin{equation}
		\delta_{X,Y}
			= (\partial\delta)(X,Y)
			= X\cdot\imath_Y\delta - Y\cdot\imath_X\delta - \imath_{\comm{X}{Y}}\delta
			= X\cdot\chi_Y - Y\cdot\chi_X - \chi_{\comm{X}{Y}}
	\end{equation}
	for all \(X,Y\in\fhhat\), where \(\partial\) denotes the differential of \(\ssC^\bullet(\fhhat;V^*\otimes\fhhat)\). Then, using \(\fa_0=\fhhat\)-invariance of the differential \(\partial:\ssC^{2,1}(\fa_-;\fa)\to\ssC^{2,2}(\fa_-;\fa)\) and the defining relation \eqref{eq:mu-partial-chi-invariance} of \(\chi_X\),
	\begin{equation}
		\partial\delta_{X,Y}
			= X\cdot\partial\chi_Y - Y\cdot\partial\chi_X - \partial\chi_{\comm{X}{Y}}
			= X\cdot(Y\cdot\mu_-) - Y\cdot(X\cdot\mu_-) - \comm{X}{Y}\cdot\mu_-
			= 0,
	\end{equation}
	so \(\delta_{X,Y}=0\) for all \(X,Y\in\fhhat\) by injectivity of \(\partial\), whence \(\delta\) is a cocycle as claimed.
\end{proof}

Of course, we can be a little more explicit about the form of the solution.

\begin{lemma}\label{lemma:alpha-delta-soln}
	Let \(V=V'\) and suppose that \(\mu_-=\alpha+\beta+\gamma+\rho\in\ssZ^{2,2}(\fa_-;\fa)\) with \([\mu_-]\in\ssH^{2,2}(\fa_-;\fa)^{\fa_0}\). Then we have the following:
	\begin{enumerate}
	\item \label{item:alpha-delta1-soln}
		The general solution \((\alpha:\Wedge^2V\to V,\,\delta_1:\fhhat\otimes V\to\fso(V))\) to equation~\eqref{eq:alpha-delta1-rsymm} is
		\begin{align}
			\alpha(v,w) &= \lambda_1(v)w - \lambda_1(w)v, \label{eq:alpha-soln}\\
			\delta_1(A+a,v) &= \comm{A}{\lambda_1(v)} - \lambda_1(Av), \label{eq:delta1-soln}
		\end{align}
		for all \(A+a\in\fhhat\) with \(A\in\fso(V)\) and \(a\in\fr\) and all \(v,w\in V\), where \(\lambda_1:V\to \fso(V)\) is some linear map. Equation~\eqref{eq:delta1-cocycle-cond-rsymm} is identically satisfied by \(\delta_1\) in this form, i.e. \(\delta_1\in\ssZ^1(\fhhat,V^*\otimes\fso(V))\).
	\item \label{item:delta-beta-soln}
		If \(\fhhat\) acts faithfully on \(S'\) then equation \eqref{eq:beta-delta12-rsymm} defines \(\delta\) uniquely in terms of \(\beta\), and both cocycle conditions \eqref{eq:delta1-cocycle-cond-rsymm} and \eqref{eq:delta2-cocycle-cond-rsymm} are then identically satisfied; that is to say, the equation
		\begin{equation}\label{eq:delta-beta-soln}
			\delta(X,v)\cdot s = (X\cdot\beta)(v,s)
		\end{equation}
		for \(X\in\fhhat\), \(v\in V\), \(s\in S'\), defines a cocycle \(\delta\in\ssZ^1(\fhhat;V^*\otimes\fhhat)\).
	\item \label{item:delta-gammahat-soln}
		If \(V=\kappa(\Odot^2S')\) then equation \eqref{eq:gamma-delta1-rsymm} defines a map \(\delta_1:\fhhat\otimes V\to\fso(V)\) uniquely in terms of \(\gamma\) which identically satisfies equation~\eqref{eq:delta1-cocycle-cond-rsymm}, and equation~\eqref{eq:rho-delta2-rsymm} defines a map \(\delta_2:\fhhat\otimes V\to\fr\) uniquely in terms of \(\rho\) which identically satisfies \eqref{eq:delta2-cocycle-cond-rsymm}; that is to say, the equation
		\begin{equation}\label{eq:delta-gammahat-soln}
			\delta(X,\kappa(s,s)) = -\qty(X\cdot\gammahat)(s,s)
		\end{equation}
		for \(X\in\fhhat\), \(v\in V\), defines a cocycle \(\delta\in\ssZ^1(\fhhat;V^*\otimes\fhhat)\).
	\end{enumerate}
\end{lemma}

\begin{proof}
	For \ref{item:alpha-delta1-soln}, the expression \eqref{eq:alpha-soln} for \(\alpha\) in terms of a unique map \(\lambda_1:V\to\fso(V)\) follows from Lemma~\ref{lemma:V-soV-maps}; substituting this expression into equation~\eqref{eq:alpha-delta1-rsymm} gives us
	\begin{equation}
		\deltatilde_1(A+a,v)w - \deltatilde_1(A+a,w)v = 0
	\end{equation}
	where \(\deltatilde_1(A+a,v)w:=\delta_1(A+a,v)w - \comm{A}{\lambda_1(v)} + \lambda_1(Av)\), but then by Lemma~\ref{lemma:V-soV-maps} we have \(\deltatilde_1=0\), giving us \eqref{eq:delta1-soln}. One can then check that \eqref{eq:delta1-cocycle-cond-rsymm} is satisfied for any map \(\delta_1\) of this form.
	
	For \ref{item:delta-beta-soln}, it is clear that equation \eqref{eq:beta-delta12-rsymm} is equivalent to \eqref{eq:delta-beta-soln} and that this defines \(\delta\) uniquely if \(\fhhat\) acts faithfully on \(S'\).
	We can combine \eqref{eq:delta1-cocycle-cond-rsymm} and \eqref{eq:delta2-cocycle-cond-rsymm} to obtain the cocycle condition for \(\ssC^1(\fhhat;V^*\otimes\fhhat)\), which is explicitly
	\begin{equation}\label{eq:delta-cocycle-condi-combined}
		\delta(\comm{X}{Y},v) - \comm{X}{\delta(Y,v)} + \comm{Y}{\delta(X,v)} - \delta(X,Y\cdot v) + \delta(Y,X\cdot v) = 0,
	\end{equation}
	and this is the equation we must verify. Acting the left-hand side of this equation on \(s\in S'\), we apply \eqref{eq:beta-delta12-rsymm} to find
	\begin{equation}
	\begin{split}
		&\delta(\comm{X}{Y},v)\cdot s
			- \comm{X}{\delta(Y,v)}\cdot s + \comm{Y}{\delta(X,v)}\cdot s 
			- \delta(X,Y\cdot v)\cdot s + \delta(Y,X\cdot v)\cdot s	\\
		& = (\comm{X}{Y}\cdot\beta)(v,s)
			- X\cdot(Y\cdot \beta)(v,s) + (Y\cdot\beta)(v,X\cdot s)
			+ Y\cdot(X\cdot \beta)(v,s) - (X\cdot\beta)(v,Y\cdot s)	\\
			&\qquad - (X\cdot \beta)(Y\cdot v,s) + (Y\cdot \beta)(X\cdot v,s)	\\
		& = (\comm{X}{Y}\cdot\beta)(v,s) - (X\cdot (Y\cdot \beta))(v,s) - (Y\cdot (X\cdot \beta))(v,s)	\\
		& = 0,
	\end{split}
	\end{equation}
	whence the cocycle condition is satisfied by faithfulness of the action of \(\fhhat\) on \(S'\).
	
	For \ref{item:delta-gammahat-soln}, if \(V=\kappa(\Odot^2S')\) then clearly \(\delta_1\) and \(\delta_2\) are determined by equations \eqref{eq:gamma-delta1-rsymm} and \eqref{eq:rho-delta2-rsymm} respectively; these can be rendered as \(\delta_1(X,\kappa_s) = -\qty(X\cdot\gamma)(s,s)\) and \(\delta_2(X,\kappa_s) = -\qty(X\cdot\rho)(s,s)\) for \(X\in\fhhat\), \(s\in S'\), which sum to \eqref{eq:delta-gammahat-soln} (and can be recovered as the \(\fso(V)\)- and \(\fr\)-components of that equation).
	For brevity, we will simply show that \eqref{eq:delta-gammahat-soln}  implies the combined cocycle equation \eqref{eq:delta-cocycle-condi-combined}; the manipulations for the conditions on \(\delta_1,\delta_2\) separately are formally identical.
	We will verify the cocycle condition for \(v=\kappa_s=\kappa(s,s)\) with \(s\in S'\), since then by the hypothesis it must hold for all \(v\in V\). We can depolarise \eqref{eq:delta-gammahat-soln} to find
	\begin{equation}
		\delta(X,\kappa(s_1,s_2)) = 
			\gammahat(X\cdot s_1,s_2) + \gammahat(s_1,X\cdot s_2)
			- \comm{X}{\gammahat(s_1,s_2)}
	\end{equation}
	for \(s_1,s_2\in S\) and \(X\in\fhhat\), so using \(\fhhat\)-equivariance of \(\kappa\), we have
	\begin{equation}
	\begin{split}
		&\delta(\comm{X}{Y},\kappa_s) 
			- \comm{X}{\delta(Y,\kappa_s)} + \comm{Y}{\delta(X,\kappa_s)}
			- \delta(X,Y\kappa_s) + \delta(Y,X\kappa_s)\\
		&= \delta(\comm{X}{Y},\kappa_s)
			- \comm{X}{\delta(Y,\kappa_s)} + \comm{Y}{\delta(X,\kappa_s)}
			- 2\delta(X,\kappa(Y\cdot s,s)) + 2\delta(Y,\kappa(X\cdot s,s))	\\
		&= 2\gammahat(\comm{X}{Y}s,s) 
			- \comm{\comm{X}{Y}}{\gammahat(s,s)} \\
			&\qquad\qquad - \cancel{2\comm{X}{\gammahat(Y\cdot s,s)}} + \comm{X}{\comm{Y}{\gammahat(s,s)}}
			+ \cancel{2\comm{Y}{\gammahat(X\cdot s,s)}} - \comm{Y}{\comm{X}{\gammahat(s,s)}}	\\
			&\qquad\qquad - 2\gammahat(X\cdot (Y\cdot s),s) - \cancel{2\gammahat(Y\cdot s,X\cdot s)} + \cancel{2\comm{X}{\gammahat(Y\cdot s,s)}}	\\
			&\qquad\qquad + 2\gammahat(Y\cdot (X\cdot s),s) + \cancel{2\gammahat(X\cdot s,Y\cdot s)} - \cancel{2\comm{Y}{\gammahat(X\cdot s,s)}}	\\
		&= 2\gammahat\qty(\comm{X}{Y}\cdot s - X\cdot (Y\cdot s) + Y\cdot (X\cdot s),s)
			- \comm{\comm{X}{Y}}{\gammahat(s,s)} + \comm{X}{\comm{Y}{\gammahat(s,s)}} - \comm{Y}{\comm{X}{\gammahat(s,s)}}
		\\
		&= 0,
	\end{split}
	\end{equation}
	as required.
\end{proof}

The hypothesis that \(\mu_-\) is a Spencer cocycle with \(\fa_0\)-invariant cohomology class is necessary in order for \(\delta\) to be well-defined. Indeed, let us relax this hypothesis.
The map \(\delta\) is not automatically well-defined by equation~\eqref{eq:delta-beta-soln} for arbitrary \(\beta:V\otimes S'\to S'\); the map \(X\cdot\beta\) must take values in the image of \(\fhhat\) in \(\End(S')\) for all \(X\in\fhhat\).
Similarly, equation~\eqref{eq:delta-gammahat-soln} is only a valid definition for maps \(\gammahat:\Odot^2S'\to\fshat\) such that \(X\cdot\gammahat\) factors through \(\kappa\).
Finally, observe that the different definitions of \(\delta\) need not agree in general, even where valid.
However, by Lemma~\ref{lemma:invt-spencer-cocycle-rsymm}, the hypothesis guarantees that a unique solution \(\delta\) exists, whence these obstructions are automatically resolved.

\subsubsection{Jacobi identities}
\label{sec:jacobi-rsymm}

There are three types of conditions on the deformation maps which arise from the Jacobi identity for \(\fatilde\): those which are independent of the deformation maps and are identically satisfied since they are Jacobi identities for \(\fa\); those which are linear in the degree-2 deformation maps, independent of the degree-4 maps and correspond to one of the cocycle conditions  \eqref{eq:22-spencer-cocycle-vss-rsymm}-\eqref{eq:delta2-cocycle-cond-rsymm} for the Chevalley--Eilenberg cocycle \(\mu=\alpha+\beta+\gamma+\rho+\delta\); and finally those which are quadratic in the degree-2 maps and linear in the degree-4 map \(\theta\). For the second type, we further distinguish between the \emph{Spencer cocycle conditions} \eqref{eq:22-spencer-cocycle-vss-rsymm} and \eqref{eq:22-spencer-cocycle-sss-rsymm} which involve only \(\mu_-=\alpha+\beta+\gamma+\rho\), the \emph{\(\delta\)-cocycle conditions}  \eqref{eq:delta1-cocycle-cond-rsymm} and \eqref{eq:delta2-cocycle-cond-rsymm} which involve only \(\delta\), and the \emph{mixed cocycle conditions} \eqref{eq:alpha-delta1-rsymm}-\eqref{eq:rho-delta2-rsymm} which involve both. 
Denoting by \([W_1,W_2,W_3]\) the component of the Jacobi identity for \(\fatilde\) on \(W_1\otimes W_2\otimes W_3\), where each \(W_i\) is one of the subspaces \(V'\), \(S'\) or \(\fhhat\), the conditions arising from each component are as follows.
\begin{itemize}
	\item
	\([\fhhat,\fhhat,\fhhat]\): This is the Jacobi identity for the Lie algebra \(\fhhat\), thus is identically satisfied.
	\item
	\([\fhhat,\fhhat,S']\): Identically satisfied; follows from \(S'\) being a representation of \(\fhhat\).
	\item
	\([\fhhat,\fhhat,V']\): Equivalent to the  \(\delta\)-cocycle conditions \eqref{eq:delta1-cocycle-cond-rsymm} and \eqref{eq:delta2-cocycle-cond-rsymm}.
	\item
	\([\fhhat,S',S']\): Using the \(\fso(V)\oplus\fr\)-equivariance of \(\kappa\), this is equivalent to the mixed cocycle conditions \eqref{eq:gamma-delta1-rsymm} and \eqref{eq:rho-delta2-rsymm}.
	\item
	\([\fhhat,S',V']\): Equivalent to mixed cocycle condition~\eqref{eq:beta-delta12-rsymm}.
	\item
	\([\fhhat,V',V']\): This has components in \(V'\), \(\fso(V)\) and \(\fr\); the first is equivalent to the mixed cocycle condition \eqref{eq:alpha-delta1-rsymm}, and the other two, which are of the quadratic type, are as follows:
	\begin{align}
		\begin{split}\label{eq:jacobi-022a-rsymm}
		(A\cdot\theta_1)(v,w) &= \delta_1(\delta_1(A+a,v)+\delta_2(A+a,v),w) \\
			&\quad - \delta_1(\delta_1(A+a,w)+\delta_2(A+a,w),v) - \delta_1(A+a,\alpha(v,w)), 
		\end{split}\\
		\begin{split}\label{eq:jacobi-022b-rsymm}
		(A\cdot\theta_2)(v,w) + \comm{a}{\theta_2(v,w)} &= \delta_2(\delta_1(A+a,v)+\delta_2(A+a,v),w) \\
			&\quad - \delta_2(\delta_1(A+a,w)+\delta_2(A+a,w),v) - \delta_2(A+a,\alpha(v,w)), 
		\end{split}
	\end{align}
	for all \(A+a\in\fhhat\) and \(v,w\in V'\).
	\item
	\([S',S',S']\): Equivalent to the second Spencer cocycle condition \eqref{eq:22-spencer-cocycle-sss-rsymm}.
	\item
	\([S',S',V']\): This has components in \(V'\), \(\fso(V)\) and \(\fr\). The first is the first Spencer cocycle condition \eqref{eq:22-spencer-cocycle-vss-rsymm}, and the others are the quadratic conditions
	\begin{align}
		\theta_1(\kappa_s,v) + \delta_1(\gamma(s,s)+\rho(s,s),v) + 2\gamma(s,\beta(v,s)) &= 0, \label{eq:jacobi-112a-rsymm}\\
		\theta_2(\kappa_s,v) + \delta_2(\gamma(s,s)+\rho(s,s),v) + 2\rho(s,\beta(v,s)) &= 0, \label{eq:jacobi-112b-rsymm}
	\end{align}
	for all \(s\in S'\), \(v\in V'\) (where \(\kappa_s=\kappa(s,s)\)).
	\item
	\([S',V',V']\): For \(s\in S'\), \(v,w\in V'\), we have the quadratic condition
	\begin{equation}\label{eq:jacobi-122-rsymm}
	\begin{split}
		&\theta_1(v,w)\cdot s + \theta_2(v,w) s \\
		& \qquad + \beta(\alpha(v,w),s) - \beta(v,\beta(w,s)) + \beta(w,\beta(v,s)) =0.
	\end{split}
	\end{equation}
	\item
	\([V',V',V']\): This has components in \(V'\), \(\fso(V)\) and \(\fr\), which are all quadratic:
	\begin{align}
		\theta_1(u,v)w + \alpha(\alpha(u,v),w) + \dots &= 0, \label{eq:jacobi-222a-rsymm}\\
		\theta_1(\alpha(u,v),w) + \delta_1(\theta_1(u,v)+\theta_2(u,v),w) + \dots &=0, \label{eq:jacobi-222b-rsymm}\\
		\theta_2(\alpha(u,v),w) + \delta_2(\theta_1(u,v)+\theta_2(u,v),w) + \dots &=0, \label{eq:jacobi-222c-rsymm}
	\end{align}
	where \(\dots\) denote the omission of cyclic permutations of the explicitly given terms in \(u,v,w\in V\).
\end{itemize}

Determining a general filtered deformation of a graded subalgebra \(\fa\subseteq\fshat\) therefore consists of solving the cocycle conditions for \(\mu\) and then checking whether a solution \(\theta\) for the quadratic equations \eqref{eq:jacobi-022a-rsymm}-\eqref{eq:jacobi-222c-rsymm} exists.

\subsection{Highly supersymmetric and transitive graded subalgebras}
\label{sec:highly-susy-subalg-rsymm}

We have already seen in Lemmas~\ref{lemma:proj-diff-injective-rsymm}, \ref{lemma:invt-spencer-cocycle-rsymm} and \ref{lemma:alpha-delta-soln} that we acquire significant simplifications when \(V'=V=\kappa(\Odot^2S')\) and \(\fhhat\) acts faithfully on \(S'\). We will argue that when \((V,\eta)\) is Lorentzian and \(\kappa\) is chosen appropriately, we can arrange for this to be the case so long as \(S'\) is sufficiently large. Then, assuming that these choices have been made, we will describe a class of filtered subdeformations of \(\fshat\) which can be identified with (subalgebras of) Killing superalgebras.

\subsubsection{The Homogeneity Theorem and faithful actions}
\label{sec:homogeneity-faithful-rsymm}

The following result is quoted from \cite[\S4.3]{Beckett2024_ksa}. We say that a symmetric Dirac current \(\kappa\) is \emph{causal} if \(\kappa_s=\kappa(s,s)\) is either time-like or null for all \(s\in S\).

\begin{theorem}[Homogeneity Theorem \cite{Figueroa-OFarrill2012,Hustler2016}]\label{thm:homogeneity}
	Let \((V,\eta)\) be a Lorentzian vector space of dimension \(\dim V>2\) and \(S\) a (possibly \(N\)-extended) spinor module of \(\fso(V)\) with a symmetric, causal Dirac current \(\kappa:\Odot^2 S \to V\). If \(S'\) is a vector subspace of \(S\) with \(\dim S'>\tfrac{1}{2}\dim S\), then \(\kappa|_{\Odot^2S'}\) is surjective onto \(V\).
\end{theorem}

With this in mind, we will from now on assume that our Dirac current \(\kappa\) is causal. We then call a graded subalgebra \(\fa=V\oplus S'\oplus\fhhat\) of \(\fshat\) with \(\dim S'> \tfrac{1}{2}\dim S\) \emph{highly supersymmetric}, and it is \emph{maximally supersymmetric} if \(S'=S\), and likewise for their filtered deformations. Theorem~\ref{thm:homogeneity} tells us that we must have \(\fa_{-2}=V=\kappa(\Odot^2S')\) for such subalgebras. The following corollary is a generalisation of \cite[Coro.4.9]{Beckett2024_ksa}, which says that \(\fann_{\fso(V)}(S')=0\) when \(\dim S'>\tfrac{1}{2}\dim S'\). We recall that the ideal \(\fhhat\cap\fr\) of \(\fhhat\) also appeared in Lemma~\ref{lemma:hhat-ext}.

\begin{corollary}\label{coro:annihilators}
	If the Dirac current \(\kappa\) is causal and \(S'\) is a subspace of \(S\) with \(\dim S'>\tfrac{1}{2}\dim S\) then \(\fann_{\fso(V)\oplus\fr}(S')=\fann_\fr(S')\). In particular, if \(\fhhat\) is a subalgebra of \(\fso(V)\oplus\fr\) such that \(\fhhat\cdot S'\subseteq S'\), then \(\fann_{\fhhat}(S')\subseteq\fann_{\fhhat\cap\fr}(S')\), so \(\fhhat\) acts faithfully on \(S'\) if and only if \(\fhhat\cap\fr\) does.
\end{corollary}

\begin{proof}
	Clearly \(\fann_\fr(S')\subseteq \fann_{\fso(V)\oplus\fr}(S')\). On the other hand, given \(A\in\fso(V)\), \(a\in\fr\) such that \((A+a)\cdot S' = 0\), by \(\fso(V)\oplus\fr\)-invariance of \(\kappa\) we have
	\begin{equation}
		A\kappa(s,s)=2\kappa(A\cdot s + as,s)=0
	\end{equation}
	 for all \(s\in\ S'\). But then \(A=0\) by the Homogeneity Theorem since \(\fso(V)\) acts faithfully on \(V\), thus \(A+a=a\in\fann_{\fr}(S')\). The second claim follows immediately.
\end{proof}

The fact that the action of \(\fhhat\) on \(S'\) is not automatically faithful (i.e. we may have \(\fann_{\fhhat}(S')\neq 0\)) is a complication to extending results on the maximally supersymmetric filtered subdeformations of \(\fs\) to those of \(\fshat\). However, the following lemmas show that simply assuming faithfulness is fairly mild.

\begin{lemma}\label{lemma:hhat-faithful-wlog}
	If the Dirac current \(\kappa\) is causal, \(R\) is compact and \(\dim S'>\tfrac{1}{2}\dim S\) then \(\fhhat\) decomposes as a sum of ideals \(\fhhat=\fhhat'\oplus\fann_{\fhhat\cap\fr}(S')\) where \(\fhhat'\) acts faithfully on \(S'\).
\end{lemma}

\begin{proof}
	If \(R\) is compact then \(\fr\) admits an \(\ad\)-invariant positive-definite inner product. Since \(\fso(V)\) is simple but not compact, its Killing form is non-degenerate but indefinite. These bilinears induce a non-degenerate, indefinite, \(\ad\)-invariant bilinear form on the sum \(\fso(V)\oplus\fr\) with respect to which \(\fso(V)\) and \(\fr\) are mutually orthogonal. We denote the restriction of this form to \(\fhhat\), which may be degenerate, by \(\pair{-}{-}\). 
	Now define \(\fhhat'\) as the orthogonal complement of the annihilator of \(S'\), \(\fhhat':=\fann_{\fhhat\cap\fr}(S')^\perp\), and note that since \(\pair{-}{-}\) is non-degenerate on \(\fann_{\fhhat\cap\fr}(S')\subseteq\fr\), we have \(\fhhat'\cap\fann_{\fhhat\cap\fr}(S')=0\) and \(\fhhat'+\fann_{\fhhat\cap\fr}(S')=\fhhat\).
	A standard Lie theory computation then shows that \(\fhhat'\) is an ideal such that \(\fhhat=\fhhat'\oplus\fann_{\fhhat\cap\fr}(S')\), and it acts faithfully on \(S'\) by Corollary~\ref{coro:annihilators}.
\end{proof}

Recall from \S\ref{sec:R-symm} that in Lorentzian signature it is always possible to choose a symmetric Dirac current \(\kappa\) such that the \(R\)-symmetry group is compact \cite{VanProeyen1999,Gall2021}. The upshot of the lemma above is that if we have so chosen \(\kappa\), and \(\fhhat\) does not act faithfully, we have \(\fa=\fa'\oplus\fann_{\fhhat}(S')\) (as Lie superalgebras) where \(\fa'=V\oplus S'\oplus\fhhat'\) is a graded subalgebra of \(\fa\) (thus of \(\fshat\)) with \(\fhhat'\) acting faithfully on \(S'\), and we can simply choose to work with \(\fa'\) instead of \(\fa\).

A different hypothesis ensures not only that the action of \emph{any} choice of \(\fhhat\) faithful, but also that the corresponding group action is effective.
We will require this in \S\ref{sec:towards-classification-rsymm} and \S\ref{sec:high-susy-spin-mfld-homogeneous-rsymm}.

\begin{lemma}\label{lemma:hhat-faithful-stronger}
	If \(\kappa\) is causal, \(\dim S'>\tfrac{1}{2}\dim S\) and \(\Spin^R(V)\cdot S'=S\) then the only element of \(\Spin_0^R(V)\) which acts trivially on \(S'\) is the identity. In particular, if the action of some subgroup \(\widehat{H}\) of \(\Spin_0^R(V)\) preserves \(S'\), then that action is effective and the action of \(\fhhat=\Lie\widehat{H}\) is faithful.
\end{lemma}

\begin{proof}	
	Since \(\kappa\) is \(\fso(V)\)-equivariant, it is \(\Spin_0(V)\) equivariant; a similar argument to the one in the proof of Corollary~\ref{coro:annihilators} shows that it is sufficient to show that no non-identity element of \(R\) acts trivially on \(S'\).
	But if \(r\in R\) such that \(rs=s\) for all \(s\in S'\), then for all \(g\in\Spin(V)\) we have \(r(g\cdot s) = g\cdot(rs) = g\cdot s\), which shows that \(r\) acts trivially on \(S=\Spin(V)\cdot S'\), whence \(r=\1\).
	The remaining claims follow immediately.
\end{proof}

\begin{remark}\label{rem:smaller-poincaré-superalgebra}
	The hypothesis \(\Spin(V)\cdot S'=S\) essentially says that \(\fshat\) is the smallest \(N\)-extended Poincaré superalgebra containing \(S'\). Indeed, in the more general case take \(\Sch:=\Spin(V)\cdot S'\). Then \(\Sch\) is a (possibly extended) spinor module, and it is preserved by \(\fhhat\).
	Noting that the Dirac current on \(S\) restricts to one on \(\Sch\) and writing the corresponding \(R\)-symmetry algebra for \(\Sch\) as \(\frch\), we can define a smaller Poincaré superalgebra \(\fsch=V\oplus\Sch\oplus(\fso(V)\oplus\frch)\) and a sequence of embeddings \(\fa\hookrightarrow\fsch\hookrightarrow\fshat\).
	Now, considered as a subdeformation of \(\fsch\) instead of \(\fshat\), \(\fa\) satisfies the hypothesis, so we can assume it without loss of generality.
	We will not lay out the details of of this here since they are somewhat tangential to the present discussion, but we include them in Appendix~\ref{appendix}.
\end{remark}

In light of the above, we now freely assume wherever necessary that the action of \(\fhhat\) on \(S'\) is faithful.

\subsubsection{Homological properties}
\label{sec:homological-props-high-susy-rsymm}

The following pair of results are analogues of those for subalgebras of \(\fs\) \cite[Lem.4.12,Prop.4.13]{Beckett2024_ksa}, and the proofs are similar, though we must use the assumption that \(\fhhat\) acts faithfully on \(S'\) in order for \(\fa\) to be transitive; in the sequel, we will freely use ``\(\fa\) is transitive'' as a synonym for ``\(\fa_0=\fhhat\) acts faithfully on \(\fa_{-1}=S'\).''
The definitions of the named properties are taken from \cite{Cheng1998} and also stated in \cite[Def.2.4]{Beckett2024_ksa}.

\begin{lemma}\label{lemma:high-susy-subalg-coho-rsymm}
	Let \(\fa=V\oplus S'\oplus\fhhat\) be a highly supersymmetric graded subalgebra of \(\fshat\). Then
	\begin{itemize}
	\item \(\fa\) is fundamental (\(\fa_-\) is generated by \(\fa_{-1}\)),
	\item \(\fa\) is transitive (\(\forall X\in\fa_k\) with \(k\geq 0\), \(\comm{X}{\fa_-}=0\implies X=0\)) if and only if \(\fhhat\) acts faithfully on \(S'\),
	\item \(\fa\) is a full prolongation of degree 2 (\(\forall d\geq 2\), \(\ssH^{d,1}(\fa_-;\fa)=0\)),
	\item \(\ssH^{d,2}(\fa_-;\fa)=0\) for all even \(d>2\).
	\end{itemize}
\end{lemma}

\begin{proof}
	Fundamentality follows immediately from the Homogeneity Theorem (Thm.~\ref{thm:homogeneity}).
	Since \(\fa_k=0\) for \(k>0\), transitivity is equivalent to the adjoint action of \(\fa_0=\fhhat\) on \(\fa_-=V\oplus S'\) being faithful. Note that an element of \(\fhhat\) annihilates \(\fa_-\) if and only if it annihilates both \(V\) and \(S'\). But by Corollary~\ref{coro:annihilators} the annihilator of \(S'\) lies in \(\fr\), and the latter already annihilates \(V\), whence \(\fann_{\fhhat}(\fa_-)=\fann_{\fhhat}(S')\). Thus \(\fhhat\) acts faithfully on \(\fa_-\) if and only if \(\fhhat\), acts faithfully on \(S'\).
	
	We have \(\ssH^{2,1}(\fa_-;\fa)=0\) by the part~\ref{item:proj-diff-injective-cohom} of	Lemma~\ref{lemma:proj-diff-injective-rsymm}, and \(\ssH^{d,1}(\fa_-;\fa)=0\) trivially for \(d>2\) since \(\ssC^{d,1}(\fa_-;\fa)=0\), so \(\fa\) is a full prolongation of degree 2. For the final part, we can write the beginning of the degree-4 Spencer complex as
		\begin{equation}\label{eq:spencer-4}
		\begin{split}
			0	&\longrightarrow \Hom(\Wedge^2V,\fhhat)	\\
				&\longrightarrow \Hom(\Wedge^3 V,V)
					\oplus \Hom(\Wedge^2V\otimes S',S')
					\oplus \Hom(V\otimes \Odot^2 S', \fhhat) 
				\longrightarrow \cdots,
		\end{split}
		\end{equation}
		and for \(\psi\in \ssC^{4,2}(\fa_-;\fa)=\Hom(\Wedge^2V,\fhhat)\), we have
		\begin{equation}
			\partial\psi (v,s,s) = \psi(v,\kappa_s)
		\end{equation}
		for all \(v\in V\), \(s\in S'\), so using the Homogeneity Theorem again we have \(\partial\psi=0\) if and only if \(\psi=0\). This shows that \(\ssH^{4,2}(\fa_-;\fa)=0\), and we also have \(\ssH^{d,2}(\fa_-;\fa)=0\) for all \(d>4\) since \(\ssC^{d,2}(\fa_-;\fa)=0\).
\end{proof}

\begin{proposition}\label{prop:highly-susy-filtered-def-rsymm}
	Let \(\fa=V'\oplus S'\oplus\fhhat\) be a graded subalgebra of \(\fshat\) and \(\fatilde\) a filtered deformation of \(\fa\) with defining sequence \((\mu,\theta,0,\dots)\). Then
	\begin{enumerate}
		\item \(\mu|_{\fa_-\otimes\fa_-}\) is a cocycle in \(\ssC^{2,2}(\fa_-;\fa)\), and
			\begin{equation}
				\qty[\mu|_{\fa_-\otimes\fa_-}] \in \ssH^{2,2}(\fa_-;\fa)^{\fa_0}.
			\end{equation}
		Furthermore, \(\mu\) is a cocycle in \(\ssC^{2,2}(\fa;\fa)\).
		\item Suppose that \(\fa\) is highly supersymmetric and transitive (\(\dim S'> \tfrac{1}{2}\dim S\) so \(V'=V=\kappa(\Odot^2 S')\), and \(\fhhat\) acts faithfully on \(S'\)). Then if \(\fatilde'\) is another filtered deformation of \(\fa\) with degree-2 deformation map \(\mu'\) such that \([\mu'|_{\fa_-\otimes\fa_-}]=[\mu|_{\fa_-\otimes\fa_-}]\) then \(\fatilde\cong\fatilde'\) as filtered Lie superalgebras.
	\end{enumerate}
\end{proposition}

\begin{proof}
	See the proof of the analogous result \cite[Prop.4.13]{Beckett2024_ksa}.
\end{proof}

This result allows us to consider the \(\fa_0\)-invariants in \(\ssH^{2,2}(\fa_-;\fa)\) as ``infinitesimal deformations" of a highly supersymmetric and transitive graded subalgebra \(\fa\). That is, to determine the possible filtered deformations of \(\fa\) up to isomorphism, one must compute this space and then determine which of its elements \([\mu_-]\) ``integrate'' to a full deformation, in the sense that there exists a degree-4 map \(\theta\) satisfying the quadratic Jacobi conditions determined in \S\ref{sec:jacobi-rsymm} (since Lemma~\ref{lemma:invt-spencer-cocycle-rsymm} already ensures existence and uniqueness of \(\delta\) for some choice of representative \(\mu_-\)).

\subsubsection{Maps in cohomology}
\label{sec:maps-in-cohomology-rsymm}

Let us now recall Lemma~\ref{lemma:H22-full-flat-model-rsymm}, which established the existence of the space of normalised cocycles \(\cHhat^{2,2}\) of \(\fshat\) and showed that it was isomorphic to the cohomology group \(\ssH^{2,2}(\fshat_-;\fshat)\). It will be useful to relate \(\ssH^{2,2}(\fa_-;\fa)\) to this and some other cohomology groups.

Let \(\fa\) be a graded subalgebra of \(\fshat\). The following maps of cochains are induced by the inclusion \(i:\fa\hookrightarrow\fshat\):
\begin{equation}\label{eq:maps-in-cohomology-def}
\begin{aligned}
	i_*&: \ssC^{\bullet,\bullet}(\fa_-;\fa) \longrightarrow \ssC^{\bullet,\bullet}(\fa_-;\fshat),
		&& \phi\longmapsto i_*\phi = i\circ\phi,\\
	i^*&: \ssC^{\bullet,\bullet}(\fshat_-;\fshat) \longrightarrow \ssC^{\bullet,\bullet}(\fa_-;\fshat),
		&& \phi\longmapsto i^*\phi = \phi\circ i,
\end{aligned}
\end{equation}
where \(\fshat\) is viewed as representation of \(\fa_-\) via the restriction of the adjoint representation of \(\fshat\). The induced maps in cohomology are also denoted by \(i_*\) and \(i^*\). 

\begin{lemma}\label{lemma:maps-in-cohomology-rsymm}
	Let \(\fa\) be a highly supersymmetric graded subalgebra of \(\fshat\) and suppose we have a splitting \(\Odot^2S\cong V\oplus \ker\kappa\) as in equation~\eqref{eq:spinor-square-decomposition} so that the space of normalised cocycles \(\ssH^{2,2}(\fshat_-;\fshat)\cong\cHhat^{2,2}\) can be defined as in Lemma~\ref{lemma:H22-full-flat-model-rsymm}. Then we have a diagram with exact rows and columns
	\begin{equation}
	\begin{tikzcd}
		& & & 0 \arrow[d]\\
		& & & \ssH^{2,2}(\fa_-;\fa) \arrow[d,"i_*"] 
		& \\
		0 \arrow[r]
		& \cKhat^{2,2}(\fa_-) \arrow[r]
		& \ssH^{2,2}(\fshat_-;\fshat) \cong \cHhat^{2,2} \arrow[r, "i^*"]
		& \ssH^{2,2}(\fa_-;\fshat) \arrow[r]
		&0
	\end{tikzcd}
	\end{equation}
	where 
	\begin{equation}\label{eq:K22}
		\cKhat^{2,2}(\fa_-) = \qty{
			\beta+\gamma+\rho \in \cHhat^{2,2} \,\middle |\,
				\beta|_{V\otimes S'}=0,\,\rho|_{\Odot^2S'}=0	}.
	\end{equation}
\end{lemma}

\begin{proof}
	For \(\beta:V\otimes S\to S\), \(\gamma:\Odot^2S\to\fso(V)\), \(\rho:\Odot^2\to\fr\) such that \(\beta+\gamma+\rho\in\cHhat^{2,2}\) is a normalised cocycle in \(\ssC^{2,2}(\fshat_-;\fshat)\), \(i^*[\beta+\gamma+\rho]=0\) if and only if \(\beta|_{V\otimes S'}+\gamma|_{\Odot^2S'}+\rho|_{\Odot^2S'}=\partial\lambda_1+\partial\lambda_2\), where \(\lambda_1:V\to\fso(V)\) and \(\lambda_2:V\to\fr\). 
	We see that \(\partial\lambda_1|_{\Wedge^2 V}=0\), so \(\lambda_1=0\) by Lemma~\ref{lemma:V-soV-maps}. Taking the \(\fr\)-component, we find that \(\rho|_{\Odot^2S'}=\partial\lambda_2|_{\Odot^2S'}\). But by the normalisation of \(\rho\), this shows that \(\rho|_{\Odot^2S'}=0\) and \(\partial\lambda_2|_{\Odot^2S'}=0\). The latter says that \(\lambda_2(\kappa_s)=0\) for all \(s\in S'\), but then by Theorem~\ref{thm:homogeneity} we have \(\lambda_2=0\). Looking at the remaining components, we find that \(\beta|_{V\otimes S'}=0\) and \(\gamma|_{\Odot^2S'}=0\), and since the former implies the latter by the cocycle condition for \(\beta+\gamma+\rho\), we find that \(\ker i^* \cong \cKhat^{2,2}(\fa_-)\).
	
	For \(\alpha:\Wedge^2V\to V\), \(\beta:V\otimes S'\to S'\), \(\gammahat:\Odot^2S'\to\fhhat\) such that \(\alpha+\beta+\gammahat\) is a cocycle in \(\ssC^{2,2}(\fa_-;\fa)\), \(i_*[\alpha+\beta+\gammahat]=0\) if and only if \(\alpha+\beta+\gammahat=\partial\lambda\) for some map \(\lambda:V\to\fshat_0=\fso(V)\oplus\fr\). Then \(\gammahat(s,s)=-\lambda(\kappa_s)\), in particular \(\lambda(\kappa_s)\subseteq\fhhat\), for all \(s\in S'\). Thus by Theorem~\ref{thm:homogeneity}, \(\lambda\) takes values in \(\fhhat\); in particular it lies in \(\ssC^{1,2}(\fa_-;\fa)\), hence \([\alpha+\beta+\gammahat]=0\). This shows that \(i_*:\ssH^{2,2}(\fa_-;\fa)\to \ssH^{2,2}(\fa_-;\fshat)\) is injective.
\end{proof}

\begin{remark}\label{rem:dirac-current-splitting}
Unless \(\fa_-=\fshat_-\) (i.e. the subalgebra in question is maximally supersymmetric), there is no direct analogue in \(\ssZ^{2,2}(\fa_-;\fshat)\) of the space of normalised cocycles \(\cH^{2,2}(\fa_-)\subseteq \ssZ^{2,2}(\fa_-;\fs)\) defined in \cite[\S4.3.3]{Beckett2024_ksa} for flat model superalgebras without \(\fr\)-extension. Cocycles in \(\ssZ^{2,2}(\fa_-;\fshat)\) can be \emph{partially} normalised by setting the \(\alpha\) component to zero -- the corresponding subspace is invariant under \(\fhhat\) and surjects onto \(\ssH^{2,2}(\fa_-;\fshat)\). Recall, though, that the normalisation condition for the \(\rho\) component in \(\cHhat^{2,2}\) relied on the \(\fshat_0\)-module splitting \(\Odot^2S\cong V\oplus \ker\kappa\). An analogous normalisation of elements of \(\ssZ^{2,2}(\fa_-;\fshat)\) would require a splitting of the short exact sequence of \(\fhhat\)-modules
\begin{equation}\label{eq:dirac-current-SES}
\begin{tikzcd}
	0 \ar[r] & \fD \ar[r] & \Odot^2S' \ar[r,"\kappa|_{\Odot^2S'}"] & V \ar[r] & 0,
\end{tikzcd}
\end{equation}
where \(\fD=\ker\kappa|_{\Odot^2S'}\), the \emph{Dirac kernel} of \(S'\); equivalently, we must find a section of the Dirac current, i.e. a map \(\Sigma:V\to \Odot^2S'\) satisfying \(\kappa\circ\Sigma=\Id_V\). With such a section, we could fully normalise our cocycles with respect to \(\Sigma\) by demanding that \(\rho(\Sigma(V))=0\) -- every cohomology class has a unique representative satisfying both this and \(\alpha=0\). The problem is that there need not exist such a section which is \(\fhhat\)-invariant. If one does not demand invariance, normalised cocycles can still be defined, but the space of such cocycles would also fail to be \(\fhhat\)-invariant.
\end{remark}

\subsection{Realisable subdeformations}
\label{sec:realisable-subdef-rsymm}

We will not discuss re-parametrising filtered deformations of general (or even transitive) highly supersymmetric graded subalgebras \(\fa\) of \(\fshat\) as we did for those of \(\fs\) in \cite[\S4.3.3]{Beckett2024_ksa} since, although they are determined (at least in the transitive case) by an \(\fa_0\)-invariant Spencer cohomology class \([\alpha+\beta+\gamma+\rho]\in \ssH^{2,2}(\fa_-;\fa)^{\fa_0}\), the more complicated form of \(\delta:\fa_0\otimes\fa_{-2}\to\fa_0\) and the fact that general cocycles in \(\ssZ^{2,2}(\fa_-;\fshat)\) cannot be invariantly normalised (see Remark~\ref{rem:dirac-current-splitting}) means the re-parametrisation is not particularly enlightening. We will instead restrict our attention to a class of deformations corresponding to a restricted type of cocycle, which we will eventually see Killing superalgebras of highly supersymmetric backgrounds as a subclass of. We once again recall from Lemma~\ref{lemma:H22-full-flat-model-rsymm} that \(\cHhat^{2,2}\) is the space of normalised Spencer \((2,2)\)-cocycles of \(\fshat\), isomorphic to \(\ssH^{2,2}(\fshat_-;\fshat)\).

\subsubsection{Admissibility}

\begin{definition}\label{def:admissible-class-rsymm}
	A cohomology class \([\alpha+\beta+\gamma+\rho]\in \ssH^{2,2}(\fa_-;\fa)\) for a transitive highly supersymmetric graded subalgebra \(\fa\) of \(\fshat\) is \emph{admissible} if there exists an \(\fa_0=\fhhat\)-invariant element \(\betatilde+\gammatilde+\rhotilde\in \cHhat^{2,2}\) such that \(i_*[\alpha+\beta+\gamma+\rho]=i^*\qty[\betatilde+\gammatilde+\rhotilde]\).
	
	We call a cocycle in \(\ssZ^{2,2}(\fa_-;\fa)\) \emph{admissible} if its cohomology class is admissible.
\end{definition}

Fixing an admissible cohomology class \([\alpha+\beta+\gamma+\rho]\), the corresponding element \(\betatilde+\gammatilde+\rhotilde\) in \((\cHhat^{2,2})^{\fa_0}\) is unique up to elements of \(\cKhat^{2,2}(\fa_-)^{\fa_0}\) (see equation~\eqref{eq:K22}). We have
\begin{equation}
	i_*(\alpha + \beta + \gamma + \delta) = i^*\qty(\betatilde + \gammatilde + \rhotilde) + \partial\lambda,
\end{equation}
where \(\lambda\in \ssC^{2,1}(\fa_-;\fshat)=\Hom(V;\fso(V)\oplus\fr)\) is uniquely determined by \(\alpha+\beta+\gamma+\rho\) and is determined up to elements of \(\ssC^{2,1}(\fa_-;\fa)=\Hom(V;\fhhat)\) by the cohomology class \([\alpha+\beta+\gamma+\rho]\). Denoting the components of \(\lambda\) by \(\lambda_1:V\to\fso(V)\) and \(\lambda_2:V\to\fr\), we can describe this even more explicitly:
\begin{equation}
\begin{aligned}
	\alpha(v,w)		&= \lambda_1(v)w - \lambda_1(w)v,
	&&\qquad
	&\beta(v,s) 	&=  \betatilde(v,s) + \lambda_1(v)\cdot s + \lambda_2(v)s,	\\
	\gamma(s,s) 	&= \gammatilde(s,s) - \lambda_1(\kappa_s),
	&&\qquad
	&\rho(s,s) 		&= \rhotilde(s,s) - \lambda_2(\kappa_s),
\end{aligned}
\end{equation}
for \(v,w\in V\), \(s\in S'\). Note that the \(\alpha\) equation uniquely determines \(\lambda_1\) by Lemma~\ref{lemma:V-soV-maps}. By Theorem~\ref{thm:homogeneity}, the \(\gamma\) equation also determines \(\lambda_1\) and the \(\rho\) equation determines \(\lambda_2\). Since \(i_*\) is injective by Lemma~\ref{lemma:maps-in-cohomology-rsymm}, the \(\fa_0\)-invariance of \(\betatilde+\gammatilde+\rhotilde\) implies that the cohomology class \([\alpha+\beta+\gamma+\rho]\) is also invariant, whence it is an infinitesimal filtered deformation (i.e. it lies in \(\ssH^{2,2}(\fa_-;\fa)^{\fa_0}\)). Explicitly, we have
\begin{equation}
\begin{aligned}
	(A\cdot\beta)(v,s) + (a\cdot\beta)(v,s) 
		&= \comm{A}{\lambda_1(v)}\cdot s - \lambda_1(Av)\cdot s - \lambda_2(Av)s + \comm{a}{\lambda_2(v)}s,\\
	(A\cdot\gamma)(s,s) + (a\cdot\gamma)(s,s)
		&= -\comm{A}{\lambda_1(\kappa_s)} + \lambda_1(A\kappa_s),	\\
	(A\cdot\rho)(s,s) + (a\cdot\rho)(s,s) &= -\lambda_2(A\kappa_s) - \comm{a}{\lambda_2(\kappa_s)},
\end{aligned}
\end{equation}
for all \(v\in V\), \(s\in S'\), \(A+a\in\fhhat\) with \(A\in\fso(V)\) and \(a\in\fr\); thus \(X\cdot(\beta+\gamma+\rho)=\partial(X\cdot\lambda)\) for \(X\in\fhhat\). The cocycle \(\delta\in \ssZ^1(\fhhat;V^*\otimes\fhhat)\) making \(\alpha+\beta+\gamma+\rho+\delta\) into a cocycle in \(\ssZ^{2,2}(\fa;\fa)\) (which we expect from Lemma~\ref{lemma:invt-spencer-cocycle-rsymm}) is given by \(i_*\delta=\partial\lambda\), where we note that \(\lambda\in\ssC^0(\fhhat;V^*\otimes(\fso(V)\oplus\fr))\), or
\begin{equation}
\begin{aligned}
	\delta_1(A+a,v) 	&= \comm{A}{\lambda_1(v)} - \lambda_1(Av),
	&&\qquad
	&\delta_2(A+a,v) 	&= \comm{a}{\lambda_2(v)} - \lambda_2(Av).
\end{aligned}
\end{equation}

The following lemma follows directly from the definition, but we make note of it nonetheless as it will be useful when working with the maps \(\betatilde\), \(\gammatilde\) and \(\rhotilde\).

\begin{lemma}\label{lemma:admiss-props-rsymm}
	We have the following for any admissible cocycle \(\alpha+\beta+\gamma+\rho\):
	\begin{itemize}
		\item For all \(v\in V\), \(s\in S'\), \(\betatilde(v,s)+\lambda_1(v)\cdot s+\lambda_2(v)s \in S'\);
		\item For all \(s\in S'\), \(\gammatilde(s,s)-\lambda_1(\kappa_s) + \rhotilde(s,s)-\lambda_2(\kappa_s)\in\fhhat\);
		\item For all \(A+a\in\fhhat\), \(v\in V\), \(\comm{A+a}{\lambda(v)}-\lambda(Av) = \comm{A}{\lambda_1(v)}-\lambda_1(Av)+\comm{a}{\lambda_2(v)}-\lambda_2(Av) \in\fhhat\). 
	\end{itemize}
\end{lemma}

\subsubsection{Integrability}
\label{sec:integrability-rsymm}

We now consider obstructions to the integration of an admissible cocycle to a filtered deformation. What is required is the degree-4 deformation map \(\theta:\Wedge^2 V\to\fhhat\) satisfying the quadratic conditions \eqref{eq:jacobi-022a-rsymm}-\eqref{eq:jacobi-222c-rsymm} imposed by the Jacobi identities -- note that the linear identities, i.e. the cocycle conditions of various types discussed in \S\ref{sec:jacobi-rsymm}, are automatically satisfied by the discussion above. We first re-parametrise the deformation, defining \(\thetatilde:\Wedge^2V\to\fshat_0\) with components \(\thetatilde_1:\Wedge^2 V\to\fso(V)\) and \(\thetatilde_2:\Wedge^2 V\to\fr\) by
\begin{align}
	\thetatilde_1(v,w) &:= \theta_1(v,w) + \lambda_1(\lambda_1(v)w - \lambda_1(w)v) - \comm{\lambda_1(v)}{\lambda_1(w)},\\
	\thetatilde_2(v,w) &:= \theta_2(v,w) + \lambda_2(\lambda_1(v)w - \lambda_1(w)v) - \comm{\lambda_2(v)}{\lambda_2(w)},
\end{align}
for \(v,w\in V\). The quadratic Jacobi identities then take a much simpler form.
\begin{itemize}
	\item Equations~\eqref{eq:jacobi-022a-rsymm} and \eqref{eq:jacobi-022b-rsymm} become
	\begin{equation}\label{eq:A-thetatilde}
		X\cdot \thetatilde_1 = 0, \qquad X\cdot \thetatilde_2 = 0,
	\end{equation}
	for all \(X\in\fhhat\); that is, \(\thetatilde\) is \(\fhhat\)-invariant.
	\item Equations~\eqref{eq:jacobi-112a-rsymm} and \eqref{eq:jacobi-112b-rsymm} become
	\begin{align}
		\thetatilde_1(v,\kappa_s) &= 2\gammatilde\qty(s,\betatilde(v,s)) - \qty(\lambda_1(v)\cdot\gammatilde)(s,s) - \qty(\lambda_2(v)\cdot\gammatilde)(s,s), \label{eq:tilde-theta1-def-rsymm}\\
		\thetatilde_2(v,\kappa_s) &= 2\rhotilde\qty(s,\betatilde(v,s)) - \qty(\lambda_1(v)\cdot\rhotilde)(s,s) - \qty(\lambda_2(v)\cdot\rhotilde)(s,s), \label{eq:tilde-theta2-def-rsymm}
	\end{align}
	and by the Homogeneity Theorem (Theorem~\ref{thm:homogeneity}) fully determine \(\thetatilde_1\) and \(\thetatilde_2\) in terms of the maps \(\betatilde\), \(\gammatilde\), \(\rhotilde\)  and \(\lambda\). By \(\fhhat\)-invariance of \(\betatilde\) and \(\gammatilde\), changing \(\lambda\) by addition of a map \(V\to\fhhat\) does not change \(\thetatilde_1\) or \(\thetatilde_2\), and by a slight rearrangement one can show that changing \(\betatilde+\gammatilde+\rhotilde\) by an element of \(\cK^{2,2}(\fa_-)\) also does not change \(\thetatilde\). Thus these equations determine \(\thetatilde\) uniquely for any admissible \emph{cohomology class}, not just for a particular representative cocycle.
	\item Equation \eqref{eq:jacobi-122-rsymm} becomes
	\begin{equation}\label{eq:tilde-theta-act-s-rsymm}
		\begin{split}
			\thetatilde_1(v,w)\cdot s + \thetatilde_2(v,w)s
			&= \betatilde\qty(v,\betatilde(w,s)) - \betatilde\qty(w,\betatilde(v,s)) \\
			&\qquad + \qty(\lambda_1(v)\cdot\betatilde)(w,s) - \qty(\lambda_1(w)\cdot\betatilde)(v,s) \\
			&\qquad + \qty(\lambda_2(v)\cdot\betatilde)(w,s) - \qty(\lambda_2(w)\cdot\betatilde)(v,s).
		\end{split}
	\end{equation}
	\item The remaining identities \eqref{eq:jacobi-222a-rsymm}, \eqref{eq:jacobi-222b-rsymm}, \eqref{eq:jacobi-222c-rsymm} are, respectively,
	\begin{gather}
		\thetatilde_1(u,v)w + \thetatilde_1(v,w)u + \thetatilde_1(w,u)v =0,\label{eq:tilde-theta1-bianchi-rsymm}\\
		\qty(\lambda_1(u)\cdot\thetatilde_1)(v,w) + \qty(\lambda_1(v)\cdot\thetatilde_1)(w,u) +	\qty(\lambda_1(w)\cdot\thetatilde_1)(u,v) = 0,\label{eq:tilde-theta1-lambda-bianchi-rsymm}\\
		\qty(\lambda_1(u)\cdot\thetatilde_2)(v,w) + \qty(\lambda_2(u)\cdot\thetatilde_2)(v,w) + \dots = 0, \label{eq:tilde-theta2-lambda-bianchi-rsymm}
	\end{gather}
	for all \(u,v,w\in V\), where the omitted terms in \(\dots\) are of course cyclic permutations.
\end{itemize}

Thus to integrate an admissible cocycle, we must show that the maps \(\thetatilde_1\), \(\thetatilde_2\) determined by \eqref{eq:tilde-theta1-def-rsymm}, \eqref{eq:tilde-theta2-def-rsymm} are well-defined, \(\fhhat\)-invariant maps such that
\begin{equation}\label{eq:theta-lambda-in-hhat}
	\thetatilde(v,w) - \lambda(\lambda(v)w-\lambda(w)v) + \comm{\lambda(v)}{\lambda(w)} \in \fhhat,
\end{equation}
and that they satisfy the remaining equations \eqref{eq:tilde-theta-act-s-rsymm}-\eqref{eq:tilde-theta2-lambda-bianchi-rsymm}. If any of these fail, the cocycle does not define a filtered deformation; it is a non-integrable infinitesimal deformation.

We begin with a number of results which will reduce the number of equations we have to check.

\begin{lemma}\label{lemma:theta-lambda-in-a0}
	Any map \(\thetatilde:\Wedge^2V\to \fso(V)\oplus\fr\) satisfying equations~\eqref{eq:tilde-theta1-def-rsymm} and \eqref{eq:tilde-theta2-def-rsymm} also satisfies the condition \eqref{eq:theta-lambda-in-hhat}.
\end{lemma}

\begin{proof}
	Let us write \(\widecheck\gamma=\gammatilde+\rhotilde:\Odot^2\to\fso(V)\oplus\fr\). Then Lemma~\ref{lemma:admiss-props-rsymm} then says that \(\widecheck\gamma(s,s)-\lambda(\kappa_s)\in\fhhat\) for all \(s\in S'\), and we can write \eqref{eq:tilde-theta1-def-rsymm} and \eqref{eq:tilde-theta2-def-rsymm} as \(\thetatilde(v,\kappa_s)=2\widecheck\gamma(s,\betatilde(v,s))-(\lambda(v)\cdot\widecheck\gamma)(s,s)\). The proof is now formally identical to that of \cite[Lem.4.19]{Beckett2024_ksa}. We must show that \(\thetatilde(v,w)\in\fhhat\) for all \(v,w\in V\), but by the Homogeneity Theorem, it is sufficient to show this for \(w=\kappa_s\) for some \(s\in S'\). We have
	\begin{equation}
	\begin{split}
		&\thetatilde(v,\kappa_s) - \lambda(\lambda(v)\kappa_s - \lambda(\kappa_s)v) + \comm{\lambda(v)}{\lambda(\kappa_s)}	\\
		& = 2\widecheck\gamma\qty(s,\betatilde(v,s)+\lambda(v)\cdot s)
			- \comm{\lambda(v)}{\widecheck\gamma(s,s)-\lambda(\kappa_s)}
			+ \lambda(\lambda(\kappa_s)v - \lambda(v)\kappa_s)	\\
		& = 2\qty(\widecheck\gamma-\lambda\circ\kappa)\qty(s,\betatilde(v,s)+\lambda(v)\cdot s)
			- \comm{\lambda(v)}{\widecheck\gamma(s,s)-\lambda(\kappa_s)}	\\
			&\qquad\qquad  + \lambda\qty(2\kappa\qty(s,\betatilde(v,s)+\lambda(v)\cdot s)
			+ \lambda(\kappa_s)v - \lambda(v)\kappa_s)	\\
		& = 2\qty(\widecheck\gamma-\lambda\circ\kappa)(s,\betatilde(v,s)+\lambda(v)\cdot s) + \qty(\comm{\widecheck\gamma(s,s)-\lambda(\kappa_s)}{\lambda(v)}
			- \lambda\qty(\widecheck\gamma(s,s)v - \lambda(\kappa_s)v)),
	\end{split}
	\end{equation}
	where in the last line we have made use of the Spencer cocycle condition \eqref{eq:cocycle-1-rsymm} as well as \(\fso(V)\oplus\fr\)-invariance of \(\kappa\). By Lemma~\ref{lemma:admiss-props-rsymm}, both terms in the last line lies in \(\fhhat\).
\end{proof}

The proof of the following is nearly identical to that of the analogous result \cite[Lem.4.20]{Beckett2024_ksa}. One applies \(\kappa(s,-)\) to \eqref{eq:tilde-theta-act-s-rsymm}, uses \(\fso(V)\oplus\fr\)-equivariance of \(\kappa\) to write \(\kappa(\thetatilde_1(v,w)\cdot s+\thetatilde_2(v,w)s,s)=\thetatilde(v,w)\kappa_s\) and then applies a cocycle condition.

\begin{lemma}\label{lemma:tilde-theta-def-2-rsymm}
	Suppose that we have maps \(\thetatilde_1:\Wedge^2 V\to \fso(V)\) and \(\thetatilde_2:\Wedge^2 V\to \fr\) satisfying equation~\eqref{eq:tilde-theta-act-s-rsymm}. Then for all \(v,w\in V\) and \(s\in S'\),
	\begin{equation}\label{eq:tilde-theta1-def2-rsymm}
	\begin{split}
		\thetatilde_1(v,w)\kappa_s
			&= 2\gammatilde\qty(\betatilde(v,s),s)w - 2\gammatilde\qty(\betatilde(w,s),s)v\\
				&\qquad - \qty(\lambda_1(v)\cdot\gammatilde)(s,s)w + \qty(\lambda_1(w)\cdot\gammatilde)(s,s)v	\\
				&\qquad - \qty(\lambda_2(v)\cdot\gammatilde)(s,s)w + \qty(\lambda_2(w)\cdot\gammatilde)(s,s)v.
	\end{split}
	\end{equation}
	This equation determines \(\thetatilde_1\) uniquely.
\end{lemma}

\begin{lemma}\label{lemma:quadr-jacobi-dependency-rsymm}
	Let \(\alpha+\beta+\gamma+\rho\) be an admissible Spencer (2,2)-cocycle for a highly supersymmetric graded subalgebra \(\fa\) of \(\fshat\), and assume that we have some maps \(\thetatilde_1:\Wedge^2 V\to\fso(V)\) and \(\thetatilde_2:\Wedge^2 V\to\fr\). Then
	\begin{itemize}
	\item Equations~\eqref{eq:tilde-theta1-def-rsymm} and \eqref{eq:tilde-theta2-def-rsymm} imply that \(\thetatilde\) is \(\fhhat\)-invariant (that is, they imply equation \eqref{eq:A-thetatilde});
	\item Equation~\eqref{eq:tilde-theta-act-s-rsymm} implies that \(\thetatilde_1+\thetatilde_2:\Wedge^2 V\to\fhhat\) is \(\fhhat\)-invariant; 
	\item Assuming equation~\eqref{eq:tilde-theta1-def-rsymm}, equation~\eqref{eq:tilde-theta1-def2-rsymm} is equivalent to the algebraic Bianchi identity~\eqref{eq:tilde-theta1-bianchi-rsymm} for \(\thetatilde_1\); in particular, equations \eqref{eq:tilde-theta1-def-rsymm} and \eqref{eq:tilde-theta-act-s-rsymm} together imply \eqref{eq:tilde-theta1-bianchi-rsymm};
	\item Equations~\eqref{eq:tilde-theta1-def-rsymm} and \eqref{eq:tilde-theta-act-s-rsymm} together imply \eqref{eq:tilde-theta1-lambda-bianchi-rsymm};
	\item Equations~\eqref{eq:tilde-theta2-def-rsymm} and \eqref{eq:tilde-theta-act-s-rsymm} together imply \eqref{eq:tilde-theta2-lambda-bianchi-rsymm}.
	\end{itemize}
\end{lemma}

\begin{proof}
	Note that equations \eqref{eq:tilde-theta1-bianchi-rsymm}-\eqref{eq:tilde-theta2-lambda-bianchi-rsymm} can be slightly more compactly presented as
	\begin{gather}
		\thetatilde(u,v)w + \thetatilde(v,w)u + \thetatilde(w,u)v =0,\\
		\qty(\lambda(u)\cdot\thetatilde)(v,w) + \qty(\lambda(v)\cdot\thetatilde)(w,u) +	\qty(\lambda(w)\cdot\thetatilde)(u,v) = 0,
	\end{gather}
	where we re-collect \(\lambda_i\) into a single map \(\lambda:V\to\fhhat\), and similarly for \(\thetatilde\); other expressions can be re-written similarly. The results then follow using \(\fso(V)\oplus\fr\)-invariance of \(\kappa\), the Spencer cocycle conditions \eqref{eq:cocycle-1-rsymm} and \eqref{eq:cocycle-2-rsymm}, \(\fhhat\)-invariance of \(\betatilde\), \(\gammatilde\) and \(\rhotilde\) and the properties of Lemma~\ref{lemma:admiss-props-rsymm}; the manipulations are formally identical to those in the proof of \cite[Lem.4.21]{Beckett2024_ksa}, so we will not repeat them here.
\end{proof}

We have thus reduced the set of conditions to be checked to the equations~\eqref{eq:tilde-theta1-def-rsymm}, \eqref{eq:tilde-theta2-def-rsymm} and \eqref{eq:tilde-theta-act-s-rsymm}. As previously noted, the first two equations determine \(\thetatilde_1\) and \(\thetatilde_2\) uniquely, so let us examine whether these are good definitions. Let us define the maps \(\Theta_1:V\otimes\Odot^2S'\to\fso(V)\) and \(\Theta_2:V\otimes\Odot^2S'\to\fr\) by
\begin{align}
	\Theta_1(v,s,s) &= 2\gammatilde\qty(s,\betatilde(v,s)) - \qty(\lambda_1(v)\cdot\gammatilde)(s,s) - \qty(\lambda_2(v)\cdot\gammatilde)(s,s), \label{eq:Theta1-def}\\
	\Theta_2(v,s,s) &= 2\rhotilde\qty(s,\betatilde(v,s)) - \qty(\lambda_1(v)\cdot\rhotilde)(s,s) - \qty(\lambda_2(v)\cdot\rhotilde)(s,s), \label{eq:Theta2-def}
\end{align}
for \(v\in V\), \(s\in S'\), so that the defining equations~\eqref{eq:tilde-theta1-def-rsymm} and \eqref{eq:tilde-theta2-def-rsymm} become \(\thetatilde_i(v,\kappa_s)=\Theta_i(v,s,s)\). 
Then the maps \(\thetatilde_i\) are well-defined if and only the \(\Theta_i\)s \emph{annihilate the Dirac kernel} \(\fD=\ker\kappa|_{\Odot^2S'}\) (see the short exact sequence \eqref{eq:dirac-current-SES}) -- that is, if \(\Theta_i(V\otimes\fD)=0\) -- and thus factor through \(V\otimes V\) as follows:
\begin{equation}\label{eq:Theta-factor-diagram-rsymm}
\begin{tikzcd}
	V\otimes \Odot^2 S' \ar[rr,"\Theta_1"]\ar[rd,"\Id\otimes \kappa"'] & & \fso(V)\\
		& V \otimes  V \ar[ru,"\thetatilde_1"'] &
\end{tikzcd} ,
\begin{tikzcd}
	V\otimes \Odot^2 S' \ar[rr,"\Theta_2"]\ar[rd,"\Id\otimes \kappa"'] & & \fr\\
		& V \otimes  V \ar[ru,"\thetatilde_2"'] &
\end{tikzcd},
\end{equation}
and we find that if this is the case, they are in fact alternating, as in the following lemma.

\begin{lemma}\label{lemma:theta-well-def-rsymm}
	If the maps \(\Theta_1\) and \(\Theta_2\) annihilate the Dirac kernel and thus factor as above, equations \eqref{eq:tilde-theta1-def-rsymm} and \eqref{eq:tilde-theta2-def-rsymm} define a maps \(\thetatilde_1:\Wedge^2V\to \fso(V)\) and \(\thetatilde_2:\Wedge^2V\to \fr\).
\end{lemma}

\begin{proof}
	Once again, the proof of the claim for \(\thetatilde_1\) is nearly identical to that of the analogous \cite[Lem.4.22]{Beckett2024_ksa}. For \(\thetatilde_2\), using a cocycle condition and \(\fhhat\)-invariance of \(\rhotilde\), we have
	\begin{equation}
	\begin{split}
		\thetatilde_2(\kappa_s,\kappa_s) 
		&= 2\rhotilde\qty(s,\betatilde(\kappa_s,s)+\lambda_1(\kappa_s)\cdot s+\lambda_2(\kappa_s) s) - \comm{\lambda_2(\kappa_s)}{\rhotilde(s,s)} \\
		&= -2\rhotilde\qty(s,\qty(\gammatilde(s,s)-\lambda_1(\kappa_s))\cdot s + \qty(\rhotilde(s,s)-\lambda_2(\kappa_s))s) - \comm{\lambda_2(\kappa_s)}{\rhotilde(s,s)} \\
		&= -\comm{\rho(s,s)-\lambda_2(\kappa_s)+\lambda_2(\kappa_s)}{\rho(s,s)}\\
		&= 0
	\end{split}
	\end{equation}
	so by the Homogeneity Theorem (Thm.~\ref{thm:homogeneity}), \(\thetatilde_2\) is alternating.
\end{proof}

\begin{definition}\label{def:integrable-cocycle-rsymm}
	Let \([\alpha+\beta+\gamma+\rho]\) be an admissible cohomology class and let \(\Theta_1:V\otimes\Odot^2S'\to\fso(V) \), \(\Theta_2:V\otimes\Odot^2S'\to\fr \) be the maps defined by equations~\eqref{eq:Theta1-def}, \eqref{eq:Theta2-def}. We say that \([\alpha+\beta+\gamma+\rho]\) is \emph{integrable} if
	\begin{enumerate}
	\item \(\Theta_1\) and \(\Theta_2\) annihilate \(\fD=\ker\kappa|_{\Odot^2S'}\), the Dirac kernel of \(S'\),
	\item The maps \(\thetatilde_1:\Wedge^2 V\to \fso(V)\), \(\thetatilde_2:\Wedge^2 V\to \fr\) defined by \(\thetatilde_i\circ\kappa=\Theta_i\) satisfy \eqref{eq:tilde-theta-act-s-rsymm}.
	\end{enumerate}
	An admissible cocycle is \emph{integrable} if its cohomology class is.
\end{definition}

\begin{theorem}[Integration of admissible cocycles]\label{thm:admiss-cocycle-integr-rsymm}
	Let \([\alpha+\beta+\gamma+\rho]\) be an admissible, integrable cohomology class for a transitive highly supersymmetric graded subalgebra \(\fa=V\oplus S'\oplus\fhhat\) of \(\fshat\) with \(i_*(\alpha+\beta+\gamma+\rho)=i^*(\betatilde+\gammatilde+\rhotilde)+\partial\lambda\) for \(\betatilde+\gammatilde+\rhotilde\in(\cHhat^{2.2})^{\fa_0}\) and some map \(\lambda:V\to\fso(V)\oplus\fr\), and let \(\thetatilde_1\), \(\thetatilde_2\) be the maps defined by equations \eqref{eq:tilde-theta1-def-rsymm} and \eqref{eq:tilde-theta2-def-rsymm}. Then the following brackets define a filtered deformation \(\fatilde\) of \(\fa\):
	\begin{equation}\label{eq:admiss-cocycle-brackets-rsymm}
	\begin{gathered}
	\begin{aligned}
		\comm{A+a}{B+b} &= \underbrace{(AB-BA) + (ab-ba)}_{\fhhat},
		\\
		\comm{A+a}{v} &= \underbrace{Av}_V + \underbrace{\comm{A}{\lambda_1(v)} - \lambda_1(Av)+\comm{a}{\lambda_2(v)}-\lambda_2(Av)}_{\fhhat},
		\\
		\comm{A+a}{s} &= \underbrace{A\cdot s + as}_{S'},
		\\
		\comm{s}{s} &= \underbrace{\kappa_s}_V + \underbrace{\gammatilde(s,s) - \lambda_1(\kappa_s) + \rhotilde(s,s) - \lambda_2(\kappa_s)}_{\fhhat},
		\\
		\comm{v}{s} &= \underbrace{\betatilde(v,s) + \lambda_1(v)\cdot s + \lambda_2(v)s}_{S'},
		\\
		\comm{v}{w} &= \underbrace{\lambda_1(v)w-\lambda_1(w)v}_V
			+ \underbrace{\theta_1(v,w) + \theta_2(v,w)}_{\fhhat},
	\end{aligned}
	\end{gathered}
	\end{equation}
	where
	\begin{equation}\label{eq:theta-thetatilde-reln}
		\theta_i(v,w)=\thetatilde_i(v,w) - \lambda_i(\lambda_1(v)w-\lambda_1(w)v) + \comm{\lambda_i(v)}{\lambda_i(w)},
	\end{equation}
	and where \(v,w\in V\), \(s\in S'\) and \(A+a,B+b\in\fhhat\) with \(A,B\in\fso(V)\), \(a,b\in\fr\).
\end{theorem}

	The brackets of the Killing superalgebra localised at a point given in equation~\eqref{eq:KSA-deformed-brackets-rsymm} are precisely of the form above, with \(-\thetatilde_1\) being the value of the Riemann tensor at that point and \(-\thetatilde_2\) being that of the \(\fr\)-symmetry field strength tensor.
	In view of this, we make the following definition; it remains to be shown that \emph{any} subdeformation of this form is (a subalgebra of) a Killing superalgebra, fully justifying the terminology, but we leave full discussion of this point this to \S\ref{sec:high-susy-spin-mfld-homogeneous-rsymm}.

\begin{definition}[Geometric realisability]\label{def:geom-real-rsymm}
	A highly supersymmetric filtered subdeformation of \(\fshat\) is \emph{geometrically realisable} if the associated graded subalgebra \(\fa\subseteq \fshat\) is transitive, and the associated Spencer cohomology class is admissible (and thus necessarily integrable).
\end{definition}

\subsection{Towards a classification scheme}
\label{sec:towards-classification-rsymm}

In \cite[\S4.5]{Beckett2024_ksa}, following the 11-dimensional Lorentzian case \cite{Figueroa-OFarrill2017_1,Santi2022}, a classification scheme for \emph{odd-generated} realisable subdeformations of the flat model \(\fs\) without \(\fr\)-extension was developed. Eliding some technicalities, a correspondence was established between isomorphism classes of such subdeformations and appropriately-defined equivalence classes of (integrable) \emph{Lie pairs} \((S',\beta+\gamma)\) consisting of a subspace \(S'\subseteq S\) with \(\dim S'>\frac{1}{2}\dim S\) and a normalised cocycle \(\beta+\gamma\in\cH^{2,2}\) of \(\fs\) such that the action of elements of the \emph{envelope} \(\fh_{(S',\beta+\gamma)}=\gamma(\fD) \subseteq \fso(V)\) (where \(\fD=\ker\kappa|_{\Odot^2S'}\) is the Dirac kernel for \(S'\)) leaves both \(S'\) and \(\beta+\gamma\) invariant \cite[Def.4.27]{Beckett2024_ksa}. Let us make some brief comments about generalising this treatment for subdeformations of \(\fshat\).

Fix a subspace \(S'\subseteq S\) with \(\dim S'>\tfrac{1}{2}\dim S\) and a normalised cocycle \(\beta+\gamma+\rho\in\cHhat^{2,2}\). We define the \emph{envelope} of such a pair to be
\begin{equation}
	\fhhat_{(S',\beta+\gamma+\rho)} := \qty(\gamma+\rho)(\fD).
\end{equation}
We call \((S',\beta+\gamma+\rho)\) a \emph{Lie pair} if both \(\beta+\gamma+\rho\) and \(S'\) are preserved by the action of elements in the envelope, in which case the latter is a Lie subalgebra of \(\fshat_0=\fso(V)\oplus\fr\) and  
\begin{equation}
	\fa_{(S',\beta+\gamma+\rho)}
	 = V \oplus S'\oplus \fhhat_{(S',\beta+\gamma+\rho)}
\end{equation}
is a graded subalgebra of \(\fshat\).

If there exists a (linear) section \(\Sigma:V\to\Odot^2S'\) of \(\kappa|_{\Odot^2S'}\) such that \(\imath_v(\beta+(\gamma+\rho)\circ\Sigma)(S')\subseteq S'\) for all \(v\in V\), one can use \(\beta+\gamma+\rho\) to define an admissible cocycle \(\mu_- = \beta + \gamma + \rho + \partial((\gamma+\rho)\circ\Sigma)\) for \(\fa_{(S',\beta+\gamma+\rho)}\) -- the argument is analogous to \cite[Lem.4.26,Prop.4.28]{Beckett2024_ksa}.
Now, by construction, if \(\mu_-\) integrates to filtered deformation \(\fatilde_{(S',\beta+\gamma+\rho)}\) then \(\fatilde_{(S',\beta+\gamma+\rho)}\) is odd-generated, and it is geometrically realisable if and only if the action of \(\fhhat_{(S',\beta+\gamma+\rho)}\cap\fr=\rho(\ker\gamma|_\fD)\) on \(S'\) is faithful (by Lemma~\ref{lemma:high-susy-subalg-coho-rsymm}).
We cannot use Lemma~\ref{lemma:hhat-faithful-wlog} to ensure faithfulness in this situation, since splitting \(\fhhat_{(S',\beta+\gamma+\rho)}\) into a kernel and a faithfully-acting ideal and passing to the latter would require us to re-define \(\rho\) in such a way that the cocycle condition may not be preserved. On the other hand, by Lemma~\ref{lemma:hhat-faithful-stronger}, if we assume that \(\Spin(V)\cdot S'=S\) then faithfulness is guaranteed -- see Remark~\ref{rem:smaller-poincaré-superalgebra} and Appendix~\ref{appendix} for commentary on this assumption.

Thus we can associate an odd-generated realisable filtered subdeformation \(\fatilde_{(S',\beta+\gamma+\rho)}\) of \(\fshat\) for every Lie pair \((S',\beta+\gamma+\rho)\) satisfying \(\Spin(V)\cdot S'=S\). 
With this established, one could now go on to develop notions of equivalence classes of Lie pairs and show that they classify the odd-generated realisable filtered subdeformations (and maximal extensions thereof) up to an appropriate notion of isomorphism in an entirely analogous manner to the case of subdeformations of \(\fs\). Since the treatment would be essentially formally identical, we will not repeat it here and simply refer the reader to \cite[\S4.5]{Beckett2024_ksa}.

\section{Highly supersymmetric Lorentzian spin-\(R\) manifolds}
\label{sec:high-susy-spin-mfld-rsymm}

We now consider the problem of reconstructing a highly supersymmetric background from its Killing superalgebra. In doing so, we generalise \cite[Sec.5]{Beckett2024_ksa} to the spin-\(R\) case. We will draw upon definitions and results from all previous parts of this work and also make extensive use of the notation and formalism for working with homogeneous spaces presented in \cite[\S{}A.4]{Beckett2024_ksa}.

\subsection{Homogeneous spin-\(R\) structures}
\label{sec:homogeneous-lorentz-spin-mfld-rsymm}

Let us first consider how to treat generalised spin structures on homogeneous spaces. For a more complete treatment, see \cite{Beckett2025_gen_spin}; similar constructions also appear in \cite{Artacho2023}. We generalise the classic notion of homogeneous \emph{spin} structures \cite{Bar1992,Cahen1993}; see also \cite[\S5.1]{Beckett2024_ksa}.

\subsubsection{Definitions}
\label{sec:homog-spin-R-defs}

Let \((G,H,\eta)\) be a metric Klein pair (a connected Lie group \(G\) and closed subgroup \(H\) with \(H\)-invariant (pseudo-)inner product \(\eta\) on \(\fg/\fh\)) with signature \((s,t)\), let \((M=G/H,g)\) be the associated homogeneous pseudo-Riemannian \(G\)-space (which we assume to be oriented) and let \(V=T_oM\cong \fg/\fh\) where \(o=H\in M\). By taking the push-forward \(g_*\) of the diffeomorphism assigned to each element \(g\in G\), the action of \(G\) on \(M\) lifts to an action on \(TM\); the latter action restricted to the isotropy group \(H\) preserves \((V,\eta)\), giving us the \emph{linear isotropy representation} \(\varphi:H\to\SO(V)\).\footnote{
	Without the assumption of orientation, the isotropy representation takes values in \(\Orth(V)\); since \(G\) is connected, \(M\) is orientable if and only if \(\varphi(H)\subseteq\SO(V)\), and in indefinite signature it is both orientable and time-orientable if and only if \(\varphi(H)\subseteq\SO_0(V)\).
	} 
The action of \(G\) on \(TM\) induces in turn an action on the special orthonormal frame bundle \(F_{SO}\), where
\begin{equation}\label{eq:g-f}
	g\cdot (f_i)_{i=1}^{\dim V} = (g\cdot f_i)_{i=1}^{\dim V} = (g_*f_i)_{i=1}^{\dim V}
\end{equation}
for all \(g\in G\) and all frames \(f=(f_i)_{i=1}^{\dim V}\). Fixing a frame \(f\) over \(o=H\), we may represent any \(A\in\SO(V)\) as a matrix \(\underline{A}\) with components \(A\indices{^i_j}\) in this frame where \(Af_i=A\indices{^i_j}f_i\). This gives us a Lie group isomorphism \(\SO(V)\xrightarrow{\cong}\SO(s,t)\), where \(A\) is mapped to the transformation on \(\RR^{s,t}\) represented in the standard basis by \(\underline{A}\), which allows us to view \(F_{SO}\) as an \(\SO(V)\)-principal bundle; we then have a \(G\)-equivariant isomorphism of \(\SO(V)\)-principal bundles
\begin{equation}\label{eq:homog-frame-bundle-assoc}
	G\times_\varphi\SO(V) \xlongrightarrow{\cong} F_{SO}
\end{equation}
given by \([g,A]\mapsto g\cdot f\cdot A\) where the left action of \(G\) on the first bundle is given by \(g\cdot[g',A]=[gg',A]\) and is compatible with the right action of \(\SO(V)\).

Suppose that there exists a lift of the linear isotropy representation \(\varphi:H\to\SO(V)\) along the canonical map \(\pihat:\Spin^R(V)\to\SO(V)\); that is, a Lie group morphism \(\varphihat:H\to\Spin^R(V)\) making
\begin{equation}\label{eq:lift-to-spin-R}
\begin{tikzcd}
	&	& \Spin^R(V) \ar[d,two heads,"\pihat"]	\\
	& H	\ar[r,"\varphi"] \ar[ur,dashed,"\varphihat"]& \SO(V)
\end{tikzcd}
\end{equation}
commute. Then \(\pihat:\Spin^R(V)\to\SO(V)\) induces a \(G\)-equivariant spin-\(R\) structure on \(M\),\footnote{
	A lift \(\hat{f}\in \Phat_o\) of \(f\in (F_{SO})_o\) induces an isomorphism \(\Spin^R(V)\cong\Spin^R(s,t)\) lifting the isomorphism \(\SO(V)\cong\SO(s,t)\) induced by \(f\), making \(\Phat\to M\) into a \(\Spin^R(s,t)\)-principal bundle.	
	}
\begin{equation}
	\varpihat: \Phat := G\times_{\varphihat}\Spin^R(V)\longrightarrow F_{SO}\cong G\times_\varphi\SO(V),
\end{equation}
where the left action of \(G\) on \(\Phat\) is given by \(g\cdot[g',A]=[gg',A]\) for \(g,g'\in G\), \(A\in\Spin^R(V)\) and is compatible with the right action of \(\Spin^R(V)\); clearly this lifts the action on \(F_{SO}\). We call this the \emph{homogeneous spin-\(R\) structure} associated to the lift \(\varphihat:H\to\Spin^R(V)\). We call the manifold \(M\) together with this spin structure the \emph{homogeneous spin-\(R\) manifold} associated to the metric Klein pair \((G,H,\eta)\) and the lift \(\varphihat\).

A classic result says that for \(G\) simply connected and \(\eta\) positive-definite, lifts of the isotropy representation along \(\pi:\Spin(V)\to\SO(V)\) classify \emph{spin} structures on a homogeneous \(G\)-space \cite{Bar1992,Cahen1993}. This result fails without both hypotheses because there are obstructions to lifting the action of \(G\) to the spin bundle (see the discussion following \cite[Lem.5.1]{Beckett2024_ksa}). Since a spin-\(R\) structure is not even a cover of the frame bundle, one cannot expect an analogous result even with the hypotheses. We therefore work in a more restrictive context, essentially assuming the existence of a lift of the action. See \cite{Beckett2025_gen_spin} for a characterisation of equivalence classes of such structures in terms of \(R\)-conjugacy classes of lifts \(\varphihat\).

\begin{definition}
	A \emph{homogeneous spin-\(R\) structure} on a homogeneous pseudo-Riemannian \(G\)-space \((M,g)\) is a spin-\(R\) structure \(\varpihat:\Phat\to F_{SO}\) on \(M\) equipped with a left action of \(G\) on \(\Phat\), compatible with the right action of \(\Spin^R(s,t)\), such that \(\varpihat\) is \(G\)-equivariant.
\end{definition}

\subsubsection{Associated bundles and connections}
\label{sec:homo-spin-R-assoc-conn}

In the homogeneous spin-\(R\) setting, associated bundles to \(\Phat\to M\), in particular the canonical \(\Rbar\)-bundle and the vector bundles discussed in \S\ref{sec:spin-R-reps-bundles}, can be viewed as associated bundles to \(G\to M\). Indeed, by composing the lift \(\varphihat\) with the canonical map \(\pi_R:\Spin^R(V)\to\Rbar\), we have the \(G\)-equivariant \(\Rbar\)-principal bundle map
\begin{equation}
	\Qbar = \Phat\times_{\pi_R}\Rbar 
		= (G\times_{\varphihat}\Spin^R(V)) \times_\varrho W
		\xlongrightarrow{\cong} G\times_{\pi_R\circ\varphihat}\Rbar;
\end{equation}
given by \([[g,A],\bar{r}]=[g,\pi_R(A)\bar{r}]\), while for any representation \(\varrho:\Spin^R(V)\to\GL(W)\), we have a natural \(G\)-equivariant vector bundle isomorphism
\begin{equation}
	\Phat\times_\varrho W = (G\times_{\varphihat}\Spin^R(V)) \times_\varrho W \xlongrightarrow{\cong} G\times_{\varrho\circ\varphihat} W
\end{equation}
given by \([g,A,w]\mapsto [g,\varrho(A)w]\). In particular,
\begin{equation}
\begin{gathered}
	TM \cong \Phat\times_{\pihat} V \cong G\times_\varphi V,	\\
	\ad\Qbar = \Phat\times_{\Ad^{\Rbar}\circ\pi_R}\fr 
	 \cong G\times_{\Ad^{\Rbar}\circ\pi_R\circ\varphihat}\fr,
\end{gathered}
\end{equation}
and similarly for other even associated bundles (\(T^*M\), \(\Wedge^pTM\), \(\ad F_{SO}\), etc.), and
\begin{equation}
	\Shatbundle := \Phat \times_{\sigmahat} S \cong G\times_{\sigmahat\circ\varphihat} S.
\end{equation}

Recall from \S\ref{sec:connections-rsymm} that on spin-\(R\) manifolds, we do not have a canonical lift of the Levi-Civita connection (as we do on a spin structure) but must instead employ a connection \(\eA\) which lifts the sum of the Levi-Civita principal connection \(\omega\) on the frame bundle \(F_{SO}\to M\) and a principal connection \(\alpha\) (for which there is no canonical choice in general) on the canonical \(\Rbar=R/\ZZ_2\)-bundle \(\Qbar\to M\) along the two-sheeted cover \(\varpihat\times\varpi_R:\Phat\to F_{SO}\times_M\Qbar\) as in equation \eqref{eq:principal-connections-add}. The following lemma, a straightforward application of Wang's Theorem \cite{Wang1958}\cite[Ch.II,Thm.11.5\&Prop.11.4]{Nomizu1963}\cite[Thm.A.6]{Beckett2024_ksa}, describes \emph{invariant} connections on homogeneous spin-\(R\) structures.

\begin{lemma}[Wang's Theorem for homogeneous spin-R structures]\label{lemma:Wang-spin-R}
	Let \((G,H,\eta)\) be a metric Klein pair with homogeneous spin-\(R\) structure \(\varpihat:\Phat=G\times_{\varphihat}\Spin^R(V)\to F_{SO}\) and canonical \(\Rbar\)-bundle \(\Qbar=G\times_{\pi_R\circ\varphihat}\Rbar\to M=G/H\) induced by a lift \(\varphihat:H\to\Spin^R(V)\) of the isotropy representation \(\varphi:H\to\SO(V)\), where \(V=\fg/\fh\). Then there is a one-to-one correspondence between \(G\)-invariant principal connections \(\eA\in\Omega^1(\Phat;\fso(V)\oplus\fr)\) on \(\Phat\to M\) and linear maps \(\Phihat:\fg\to\fso(V)\oplus\fr\) such that
	\begin{enumerate}
		\item \(\Phihat\circ\Ad^G_h = \Ad^{\Spin^R(V)}_{\varphihat(h)}\circ\Phihat\) for all \(h\in H\),
		\item \(\Phihat(X)=\varphihat_*(X)\) for all \(X\in\fh\);
	\end{enumerate}
	the maps \(\Phihat\) are known as \emph{Nomizu maps}. The correspondence is given by
	\begin{equation}
		\Phihat(X) = \eA_p\qty(\widehat\xi_X)
	\end{equation}
	where \(X\in \fg\), \(\widehat\xi_X\) is the fundamental vector field on \(\Phat\) associated to \(X\) and \(p=[1_G,\1]\) is the natural basepoint in \(\Phat_o\). Moreover, the compositions \(\pr_{\fso(V)}\circ\Phihat\) and \(\pr_{\fr}\circ\Phihat\) are Nomizu maps for the connections \(\omega\in\Omega^1(F_{SO};\fso(V))\) and \(\alpha\in\Omega^1(\Qbar;\fr)\) such that \(\eA = \varpihat^*\omega + \varpi_R^* \alpha\) as in \eqref{eq:principal-connections-add}. The connection \(\eA\) is torsion-free, or equivalently \(\omega\) is the Levi-Civita connection, if and only if
	\begin{equation}\label{eq:homog-spin-R-tosion-free}
			\pr_{\fso(V)}\qty(\Phihat(X))\cdot\overline{Y} - \pr_{\fso(V)}\qty(\Phihat(Y))\cdot\overline{X} - \overline{\comm{X}{Y}} = 0
	\end{equation}
	for all \(X,Y\in\fg\), where the \(X\mapsto\overline{X}\) denotes the natural map \(\fg\to\fg/\fh \cong V\).
	
	Finally, the curvature \(\Omega_\eA\in\Omega_G^2(\Phat;\fso(V)\oplus\fr)\) is given by
	\begin{equation}\label{eq:homog-spin-R-curvature}
		2(\Omega_\eA)_p\qty(\widehat\xi_X,\widehat\xi_Y) = \comm{\Phihat(X)}{\Phihat(Y)}_{\fso(V)\oplus\fr} - \Psi(\comm{X}{Y}_\fg),
	\end{equation}
	and by \eqref{eq:principal-curvatures-add}, the \(\fso(V)\)- and \(\fr\)-components of this expression give the curvature values \((\Omega_\omega)_p(\widehat\xi_X,\widehat\xi_Y)\) and \((\Omega_\alpha)_p(\widehat\xi_X,\widehat\xi_Y)\) respectively. These values determine the curvature 2-forms uniquely by homogeneity.
\end{lemma}

We note that in indefinite signature, if \(\varphihat(H)\subseteq \Spin^R_0(V)\), in particular if \(H\) is connected, then \(M\) is both oriented and time-oriented. The structure groups then reduce as in Lemma~\ref{lemma:time-orientable-spin-R}, and we have a \(\pihat_0\)-invariant bundle map \(\varpihat_0:\Phat_0 = G\times_{\varphihat}\Spin_0(V)\to F_{SO_0}\cong G\times_{\varphi}\SO_0(V)\); the action of \(G\) preserves both bundles and \(\varpihat_0\) is \(G\)-equivariant. Thus, given an invariant connection \(\eA\) on \(\Phat\) and a Dirac current \(\kappa:\Odot^2S\to V\) on \(S\), Lemma~\ref{lemma:bundle-dirac-current-exist-rsymm} gives us a \(G\)-invariant bundle Dirac current \(\kappabold:\Odot^2\Shatbundle\to TM\). In the present setting, this could also be defined via Frobenius reciprocity (see e.g. \cite[Prop.A.5]{Beckett2024_ksa}) by noting that \(\kappa\in\qty(\Odot^2S^*\otimes V)^H\).

\subsection{Reconstruction of highly supersymmetric backgrounds}
\label{sec:high-susy-spin-mfld-homogeneous-rsymm}

We say that a Lorentzian spin-\(R\) manifold \((M,g)\) with twisted spinor bundle \(\Shatbundle\) with fibre \(S\) and admissible pair \((D,\rhobold)\) is \emph{highly supersymmetric} if \(\dim\fShat_D>\frac{1}{2}\dim S\). Equivalently, \(\fKhat_{(D,\rhobold)}\) is highly supersymmetric as a subdeformation of \(\fshat\). As in the case without \(R\)-symmetry \cite{Figueroa-OFarrill2012,Hustler2016}\cite[Ex.5.4]{Beckett2024_ksa}, the Homogeneity Theorem (Thm.~\ref{thm:homogeneity}) implies that a highly supersymmetric Lorentzian spin-\(R\) manifold is locally a homogeneous space for \(\fKhat_{(D,\rhobold)}\). If \(M\) is geodesically complete, this is a global statement.

Our goal now is to begin with a highly supersymmetric filtered subdeformation \(\fg\) of the Poincaré superalgebra \(\fshat\) and construct from it a homogeneous Lorentzian space carrying a homogeneous spin-\(R\) structure with all the data needed to define a Killing superalgebra, and to show that \(\fg\) embeds into such an algebra, if it does indeed exist.

\subsubsection{Super Harish-Chandra pairs, connections and curvature}

We recall the following definition of a model for Lie supergroups -- see \cite[\S5.2.1]{Beckett2024_ksa} and the references therein for more relevant examples and discussion.	
	
\begin{definition}[Super Harish-Chandra pair]\label{def:super-Harish-Chandra}
	A \emph{super Harish-Chandra pair} is a pair \((G_{\overline0},\fg)\) where \(\fg\) is a finite-dimensional Lie superalgebra and \(G_{\overline0}\) is a connected Lie group integrating \(\fg_{\overline0}\) equipped with a representation \(\Ad: G_{\overline0}\to\GL(\fg)\) (which we call the \emph{adjoint representation} of the pair) integrating the restriction \(\ad:\fg_{\overline0}\to\fgl(\fg)\) of the adjoint representation of \(\fg\) to \(\fg_{\overline0}\).
\end{definition}

\begin{example}\label{example:poincaré-group-R-ext-HC-pair}	
	Consider the group \(\Spin^R(V)\ltimes V\), where we take the operation on \(V\) to be addition and the action of \(\Spin^R(V)\) on \(V\) to be the natural one via \(\pihat\), and note that its Lie algebra is \(\fshat_{\overline 0}=(\fso(V)\ltimes V)\oplus\fr\).
	This is an extension of the Poincaré group \(\SO(V)\ltimes V\) by \(R\); there is a short exact sequence
	\begin{equation}
	\begin{tikzcd}
		1 \ar[r] & R \ar[r] 
			&\Spin^R(V)\ltimes V \ar[r,"\pihat\times\1"]
			&\SO(V)\ltimes V \ar[r]
			& 1
	\end{tikzcd}
	\end{equation}
	which integrates the bottom row of the diagram in Lemma~\ref{lemma:hhat-ext}.
	We can define a super Harish-Chandra pair \((\Spin^R(V)_0\ltimes V,\fshat)\), where \(\Spin^R(V)_0=\Spin_0(V)\times_{\ZZ_2}R_0\) is the connected component of the identity of \(\Spin^R(V)\), by extending the adjoint representation of \(\Spin^R(V)_0\ltimes V\) on \(\fshat_{\overline 0}\) to \(\fshat\) by declaring that elements of \(\Spin^R(V)_0\) act on \(\fshat_{\overline 1} = S\) via the defining representation \(\sigmahat\) and those of \(V\) act trivially.
\end{example}

Now let \((V,\eta)\) be a Lorentzian inner product space and \(\fg\) a filtered deformation of a highly supersymmetric graded subalgebra \(\fa=V\oplus S'\oplus\fhhat\) of the extended Poincaré superalgebra \(\fshat\) as in Section~\ref{sec:highly-susy-subalg-rsymm}. 
Note that the zeroth filtered subspace \(\fg^0\subseteq\fg\), which is a subalgebra, is canonically isomorphic as a Lie algebra to \(\fa_0=\fhhat\), with which we now identify it. Its adjoint actions on \(\fg\) and \(\fa\) preserve the subspaces \(\fg^i\subseteq \fg\) and \(\fa_i\subseteq\fa\) as well as the natural mapping \(\fg^i \to \Gr_i\fg = \fg^i/\fg^{i+1} = \fa_i\). 

If \((G_{\overline0},\fg)\) is a super Harish-Chandra pair such that the connected Lie subgroup \(\Hhat\) of \(G_{\overline 0}\) generated by \(\fhhat\) is closed then \((G_{\overline0},\Hhat,\eta)\) is a metric Klein pair, where we identify \(V=\fg^{-2}/\fg^{-1}\cong \fg_{\overline 0}/\fhhat\) as \(\fhhat\)-modules. 
The following lemma gives sufficient conditions for there to exist a lift of the isotropy representation to \(\Spin^R(V)\), and thus a corresponding homogeneous spin-\(R\) structure on \((G_{\overline0},\Hhat,\eta)\) as in \S\ref{sec:homog-spin-R-defs}.
The proof is in the spirit of an analogous result for spin structures associated to filtered subdeformations of the Poincaré superalgebra without \(\fr\)-symmetry \cite[Lem.5.5]{Beckett2024_ksa}, but here we must use the additional hypothesis that \(\Spin(V)\cdot S'=S\).
Recall that this is a stronger version of our customary assumption that \(\fa\) is transitive (i.e. \(\fhhat\) acts faithfully on \(S'\)) and amounts to \(\fshat\) being the smallest Poincaré superalgebra containing \(\fa\)  (see Lemma~\ref{lemma:hhat-faithful-stronger} and Remark~\ref{rem:smaller-poincaré-superalgebra}).

\begin{lemma}\label{lemma:spinR-lift}
	Let \((G_{\overline0},\fg)\) be a super Harish-Chandra pair where \(\fg\) is a filtered deformation of a transitive highly supersymmetric graded subalgebra \(\fa\subseteq\fshat\) such that \(\Spin(V)\cdot S' = S\), and suppose that the connected subgroup \(\Hhat\) of \(G_{\overline0}\) generated by \(\fhhat\) is closed.
	Then the inclusion \(\fhhat\hookrightarrow\fso(V)\oplus\fr\) integrates to a unique morphism \(\varphihat:\Hhat\to\Spin^R(V)\) lifting the isotropy representation \(\varphi:\Hhat\to\SO(V)\) along \(\pihat:\Spin^R(V)\to \SO(V)\).
\end{lemma}

\begin{proof}
	Let \(\varphi:\Hhat\to\SO(V)\) denote the isotropy representation of the metric Klein pair \((G_{\overline0},\Hhat,\eta)\) and let \(\pi_{\Hhat}:\Htilde\to \Hhat\) denote the universal cover of \(\Hhat\). The inclusion map \(\fhhat\hookrightarrow\fso(V)\oplus\fr\) integrates to 
	a Lie group morphism \(\varphitilde:\Htilde\to\Spin^R(V)\) making the solid arrows of the diagram
		\begin{equation}
		\begin{tikzcd}
			& \Htilde \ar[d,"\pi_{\Hhat}"']\ar[r,"\varphitilde"]
				& \Spin^R(V) \ar[d,"\pihat"]	\\
			& \Hhat \ar[r,"\varphi"]\ar[ur,dashed,"\varphihat"]	& \SO(V)
		\end{tikzcd}
		\end{equation}
	commute. Any lift \(\varphihat\) of the isotropy representation \(\varphi\) makes the whole diagram commute and is the unique map which does so, being the factorisation of \(\varphitilde\) through \(\pi_{\Hhat}:\Htilde\to \Hhat\).
	Such a map exists if and only if \(\ker\pi_{\Hhat}\subseteq\ker\varphitilde\subseteq\Htilde\), and we will demonstrate that this is the case.

	Observe that \(S\) carries a representation of \(\Htilde\) via \(\varphitilde\) which integrates the representation of \(\fhhat\) on \(S\) and therefore must preserve \(S'\), since \(\Htilde\) is connected.
	On the other hand, \(\Htilde\) also acts on \(S'\) via \(\Ad|_{\Hhat}\circ\pi_{\Hhat}\), where \(\Ad\) is the adjoint representation of \((G_{\overline 0},\fg)\). These two actions integrate the same representation of \(\fhhat\) on \(S'\), so they must agree.
	Thus, if \(k\in\ker\pi_{\Hhat}\subseteq\Htilde\),
	\begin{equation}
		\varphitilde(k)\cdot s = \Ad_{\pi_{\Hhat}(k)}s = s,
	\end{equation}
	but then since \(\Spin(V)\cdot S'=S\), we find that 
	\(\varphitilde(k)=\1\) by Lemma~\ref{lemma:hhat-faithful-stronger},	so \(\ker\pi_{\Hhat}\subseteq \ker\varphitilde\) as claimed.
\end{proof}

We also get a canonical invariant connection on the induced homogeneous spin-\(R\) structure if \(\fg\) is geometrically realisable (Definition~\ref{def:geom-real-rsymm}).

\begin{lemma}\label{lemma:homog-spinR-connection}
	Let \((G_{\overline0},\fg)\) be a super Harish-Chandra pair where \(\fg\) is a geometrically realisable filtered deformation of a transitive highly supersymmetric graded subalgebra \(\fa\subseteq\fshat\) such that \(\Spin(V)\cdot S'=S\).
	Suppose that the connected subgroup \(\Hhat\) of \(G_{\overline0}\) generated by \(\fhhat\) is closed, so that \((G_{\overline0},\Hhat,\eta)\) is a metric Klein pair.
	Then there exists a \(\Hhat\)-equivariant map \(\Phihat:\fg_{\overline 0}\to\fso(V)\oplus\fr\) extending the inclusion \(\fhhat \xhookrightarrow{\quad} \fso(V)\oplus\fr\) which is the Nomizu map for a \(G_{\overline 0}\)-invariant principal connection \(\eA\) on the homogeneous spin-\(R\) structure on \(M=G_{\overline 0}/\Hhat\) which lifts the Levi-Civita connection. 
	The values of the Riemann curvature tensor and \(\fr\)-symmetry field strength at \(o\in\ M\) are then given by
	\begin{equation}\label{eq:curvatures-homogeneous}
		R\qty(\xi_X,\xi_Y)_o = -\thetatilde_1(\overline{X},\overline{Y}),
		\qquad
		F\qty(\xi_X,\xi_Y)_o = -\thetatilde_2(\overline{X},\overline{Y}),
	\end{equation}
	where \(X,Y\in\fg_{\overline 0}\), \(\xi_X,\xi_Y\) are their fundamental vector fields on \(M\), \(\overline{X},\overline{Y}\in V\) are their projections to \(V\), and \(\thetatilde\) is the degree-4 deformation map determined by the Spencer cohomology class associated to \(\fg\).
\end{lemma}

\begin{proof}
	By Lemma~\ref{lemma:spinR-lift}, we have a lift \(\varphihat\) of the isotropy representation, giving us a homogeneous spin-\(R\) structure \(\Phat:= G\times_{\varphihat}\Spin^R(V)\to F_{SO}\).
	Let us write the brackets of \(\fg\) as in  Theorem~\ref{thm:admiss-cocycle-integr-rsymm}; for brevity, we will not write the \(\fso(V)\)- and \(\fr\)-components of elements of \(\fhhat\) or \(\fshat_0\) separately.
	Writing \(\fg_{\overline 0}\simeq\fa_{\overline 0}=\fhhat\oplus V\) (as a vector space), we define \(\Phihat\) as follows:
	\begin{equation}\label{eq:Nomizu-map-Phihat-def}
		\Phihat(X+v) = X + \lambda(v)
	\end{equation}
	for \(v\in V\), \(X\in\fhhat\), where \(\lambda:V\to\fso(V)\oplus\fr\) is the map appearing in the parametrisation of \(\fg\); we claim that \(\Phihat\) is a Nomizu map for a torsion-free connection as described in Lemma~\ref{lemma:Wang-spin-R}.
	We immediately see by putting \(v=0\) that \(\Phihat|_{\fhhat} = \varphihat_*:\,\fhhat \xhookrightarrow{\quad} \fso(V)\oplus\fr\), which is indeed the inclusion map.
	For \(\Hhat\)-invariance, since \(\Hhat\) is connected, it is sufficient for us to show that \(\Phihat\) is \(\fhhat\)-invariant. 
	Indeed, for \(Y\in\fhhat\),
	\begin{equation}
	\begin{split}
		\Phihat\big(\ad^{\fg_{\overline 0}}_{Y}(X+v)\big)
			&= \Phihat\qty(\comm{Y}{X} + \comm{Y}{\lambda(v)} - \lambda(Y\cdot v) + Y\cdot v )	\\
			&= \comm{Y}{X + \lambda(v)}	\\
			&= \ad^{\fso(V)\oplus\fr}_Y\qty(\Phihat(X+v))
	\end{split}
	\end{equation}
	as required.
	It remains to check the torsion-free condition \eqref{eq:homog-spin-R-tosion-free}. Denoting projection \(\fg_{\overline 0}\to V\cong\fg_{\overline 0}/\fhhat\) by \(X\mapsto \overline{X}\), we of course have \(\overline{X+v}=v\), so we compute
	\begin{equation}
		\Phihat(X+v)\cdot w - \Phihat(Y+w)\cdot v
			= X\cdot w + \lambda(v)\cdot w - Y\cdot v - \lambda(w)\cdot v
			= \overline{\comm{X+v}{Y+w}},
	\end{equation}
	whence \(\Phihat\) is a Nomizu map for a connection \(\eA\) lifting the Levi-Civita connection.	The curvature \(\Omega_\eA\) can now be computed using equation~\eqref{eq:homog-spin-R-curvature}. For \(X,Y\in\fhhat\) and \(v,w\in V\), we have
	\begin{equation}
	\begin{split}
		2\Omega_\eA \qty(\widehat\xi_{X+v},\widehat\xi_{Y+w})_p
			&= \comm{\Phihat(X+v)}{\Phihat(Y+w)}_{\fso(V)\oplus\fr} - \Phihat\qty(\comm{X+v}{Y+w}_{\fg_{\overline 0}})	\\
			&= \comm{X+\lambda(v)}{Y+\lambda(w)}	\\
			&\qquad - \theta(v,w)
				- \comm{X}{\lambda(w)} + \lambda(X\cdot w)
				+  \comm{Y}{\lambda(v)} - \lambda(Y\cdot v)
				- \comm{X}{Y}	\\
			&\qquad - \lambda(\lambda_1(v)w-\lambda_1(w)v+X\cdot w-Y\cdot v)	\\
			&= - \theta(v,w) - \lambda(\lambda_1(v)w-\lambda_1(w)v) + \comm{\lambda(v)}{\lambda(w)}	\\
			&= - \thetatilde(v,w),
	\end{split}
	\end{equation}
	where the subscripts on the brackets indicate which Lie algebra they are to be taken in, and the last equality is just equation~\eqref{eq:theta-thetatilde-reln}, which is part of the definition of the presentation of the brackets. Equation~\eqref{eq:curvatures-homogeneous} follows upon taking \(\fso(V)\)- and \(\fr\)-components.
\end{proof}

We already saw the identification of the components of \(-\thetatilde\) with the curvature tensors when comparing the presentation of the brackets in Theorem~\ref{thm:admiss-cocycle-integr-rsymm} (equations~\eqref{eq:admiss-cocycle-brackets-rsymm} and \eqref{eq:theta-thetatilde-reln}) with those of a Killing superalgebra localised at a point (equations~\eqref{eq:KSA-deformed-brackets-rsymm}-\eqref{eq:KSA-deformed-brackets-theta2}).

\subsubsection{Reconstruction theorem}

We conclude by providing the following generalisation of the Reconstruction Theorem for geometrically realisable subdeformations of \(\fs\) \cite[Thm.5.6]{Beckett2024_ksa} to those of \(\fshat\).

\begin{theorem}[Reconstruction of highly supersymmetric backgrounds]\label{thm:reconstruction-homog-spin-mfld-rsymm}
	Let \(\fg\) be a geometrically realisable filtered deformation of a highly supersymmetric graded subalgebra \(\fa=V\oplus S'\oplus\fhhat\) of the \(\fr\)-extended Poincaré superalgebra \(\fshat\) such that \(\Spin(V)\cdot S'=S\).
	Let \(G_{\overline 0}\) be the connected and simply connected Lie group corresponding to the Lie algebra \(\fg_{\overline 0}\) and suppose that the connected subgroup \(\Hhat\subseteq G_{\overline 0}\) corresponding to the subalgebra \(\fhhat\subseteq \fg_{\overline 0}\) is closed.
	Then there exists a homogeneous spin-\(R\) structure \(\varpihat:\Phat\to F_{SO}\) on the oriented and time-oriented Lorentzian homogeneous space \(M=G_{\overline 0}/\Hhat\), a \(G_{\overline 0}\)-invariant connection \(D\) on the twisted spinor bundle \(\Shatbundle\to M\), a \(G_{\overline 0}\)-invariant Dirac current \(\kappabold\in\Hom(\Odot^2,TM)\), a \(G_{\overline 0}\)-invariant bundle map \(\rhobold\in\Hom(\Odot^2\Shatbundle,\ad\Qbar)\) (where \(\Qbar\to M\) is the canonical \(\Rbar\)-bundle), and an injective linear map \(\psi:\fg\hookrightarrow \fsymm(\varpihat,\eA)\oplus\fShat_D\) which restricts to a Lie algebra embedding \(\fg_{\overline 0}\hookrightarrow\fVhat_{(D,\rhobold)}\) satisfying
	\begin{align}\label{eq:reconstruct-brackets-rsymm}
		\psi(\comm{X+a}{s}) 
			&= \eLhat_{\xi_X}\psi(s) + a_X\cdot\psi(s),
		& \psi\qty(\comm{s}{s'}) 
%			&= \xi_{\comm{s}{s'}} + a_{\comm{s}{s'}}
			&= \kappabold\qty(\psi(s),\psi(s')) + \rhobold\qty(\psi(s),\psi(s')),
	\end{align}
	for all \(X\in\fg_{\overline 0}\) and \(s,s'\in \fg_{\overline 1}\), where \(\xi_X\in\fiso(M,g)\) and \(a_X\in\fR\) are the vectorial and \(\fr\)-symmetry parts of \(\psi(X)=\xi_X+a_X\).
	
	Furthermore, if \(\kappabold(\epsilon,\zeta)+\rhobold(\epsilon,\zeta)\) preserves \(\betabold\) and \(\rhobold\) and satisfies \(\imath_{\kappabold(\epsilon,\zeta)}F = \nablahat(\rhobold(\epsilon,\zeta))\) for all \(\epsilon,\zeta\in\fShat_D\) so that \((D,\rhobold)\) is an admissible pair and \(\fKhat_{(D,\rhobold)}\) is a Killing superalgebra (Def.~\ref{def:killing-spinor-rsymm}), in particular if \(\psi(\fg_{\overline 1})=\fShat_D\), then \(\psi\) is a strict embedding of \(\fg\) into \(\fKhat_{(D,\rhobold)}\) as a filtered Lie superalgebra.
\end{theorem}

\begin{proof}
	Since the map \(\fg\to\fX(M)\), \(X\mapsto\xi_X\) assigning elements of \(\fg\) to their fundamental vector fields on \(M\) is naturally a Lie algebra \emph{anti}-homomorphism (essentially because the fundamental vector fields for a \emph{left} \(G\)-action are obtained by pushing forward the \emph{right}-invariant vector fields on \(G\)), it is convenient to assume that \(\fg\) is in fact the \emph{opposite} of the filtered subdeformation in question and construct \(\psi\) as an \emph{anti}-embedding.\footnote{
			The opposite superalgebra has the same underlying vector space but bracket \(\comm{X}{Y}^{\mathrm{opp}}=\comm{Y}{X}\). The map \(X\mapsto -X\) gives a Lie algebra isomorphism \(\fg^{\mathrm{opp}}_{\overline 0}\cong \fg_{\overline 0}\), but this does not extend to an isomorphism of superalgebras in general.
			If \(\fg_{(\mu,\theta,\dots)}\) is the filtered deformation of \(\fa\subseteq\fshat\) with defining sequence \((\mu,\theta,\dots)\), there is an isomorphism \(\fg_{(\mu,\theta,\dots)}^{\mathrm{opp}}\cong\fg_{(-\mu,\theta,\dots)}\) given by \(v+s+X\mapsto v+s-X\), where \(v\in\fa_{-2}\), \(s\in \fa_{-1}\) and \(X\in\fa_0\).
		%	On the other hand, \(v+s+X\mapsto -v\pm s-X\) defines an isomorphism from \(\fg^{\mathrm{opp}}\) to a filtered subdeformation of the flat model with Dirac current \(-\kappa\) and \(\mu'=\beta-\gammahat+\partial\lambda\) (now assuming geometric realisability).
		}
	The result as stated then immediately follows.
	
	Geometric realisability means that the Spencer cohomology class \([\mu_-]\in \ssH^{2,2}(\fa_-;\fa)\) associated to the opposite of \(\fg\) is admissible (Defs.~\ref{def:admissible-class-rsymm}, \ref{def:geom-real-rsymm}); if we fix a representative \(\mu_-\in \ssZ^{2,2}(\fa_-;\fa)\), there exists an \(\fa_0=\fhhat\)-invariant cocycle of \(\fshat\), \(\beta+\gamma+\rho\in\cHhat^{2,2}\) (note that we no longer denote the components of this cocycle with tildes), and a linear map \(\lambda\in\ssC^{2,1}(\fa_-;\fshat)=\Hom(V,\fso(V)\oplus\fr)\) such that
	\begin{equation}
		i_*\mu_-=i^*(\beta+\gamma+\rho)+\partial\lambda\in \ssZ^{2,2}(\fa_-;\fshat).
	\end{equation}
	
	Since \(G_{\overline 0}\) is simply connected, the adjoint representation \(\ad^{\fg}|_{\fg_{\overline 0}}:\fg_{\overline 0}\to\fgl(\fg)\) integrates uniquely to a group representation \(\Ad:G_{\overline 0}\to\GL(\fg)\), making \((G_{\overline 0},\fg)\) into a super Harish-Chandra pair.
	By Lemma~\ref{lemma:spinR-lift}, there exists a lift \(\varphihat:\Hhat\to\Spin^R(V)\) of the linear isotropy representation inducing a homogeneous spin-\(R\) structure \(\varpihat:\Phat:=G_{\overline0}\times_{\varphihat}\Spin^R(V)\to F_{SO}\) on \(M=G/H\), and Lemma~\ref{lemma:homog-spinR-connection} gives us a	Nomizu map \(\Phihat:\fg_{\overline 0}\to\fso(V)\oplus\fr\) for a \(G_{\overline 0}\)-invariant torsion-free connection \(\eA\) on \(\Phat\to M\), except that since \(\fg\) is the opposite of a filtered deformation, \(\varphihat_*:\fhhat\to\fso(V)\oplus\fr\) must be \emph{minus} the inclusion, and \(\Phihat(X+v)=-X-\lambda(v)\).
	We denote by \(\nablahat\) the covariant derivative induced on \(\Shatbundle=G_{\overline 0}\times_{\sigmahat\circ\varphihat} S \to M\) (and other associated bundles) by \(\eA\). Since \(\Hhat\) is connected, the Dirac current \(\kappa\) chosen in the definition of \(\fshat\) induces a \(G_{\overline 0}\)-invariant, \(\nablahat\)-parallel bundle Dirac current \(\kappabold\) on \(\Shatbundle\) (see the discussion at the end of \S\ref{sec:homo-spin-R-assoc-conn}).
	
	The components \(\beta\) and \(\rho\) of the cocycle \(\mu_-\) are \(\fa_0=\fhhat\)-invariant, hence they are \(\Hhat\)-invariant since \(\Hhat\) is connected and has Lie algebra \(\fhhat^{\mathrm{opp}}\), so we can construct \(G_{\overline 0}\)-invariant bundle maps \(\betabold\in\Omega^1(M;\End\Shatbundle)\) and \(\rhobold\in\Hom(\Odot^2\Shatbundle,\ad\Qbar)\) via Frobenius reciprocity (see e.g. \cite[Prop.A.5]{Beckett2024_ksa}). We then have a \(G_{\overline 0}\)-invariant connection \(D=\nablahat-\betabold\) on \(\Shatbundle\). 
	
	Recall now the bundle \(\eEhat=TM\oplus_M\oplus\Shatbundle\oplus_M\ad\Phat\) from \S\ref{sec:localise-ksa-at-point-rsymm} and note that we have a natural isomorphism \(\eEhat\cong G_{\overline 0}\times_\varrho\fshat\), where \(\varrho:\Hhat\to\GL(\fshat)\) is composition of the natural action of \(\Spin^R(V)_0\) on \(\fshat\) with \(\varphihat\).	
	By standard theory of associated bundles, there is another natural isomorphism
	\begin{equation}
		\fX(M)\oplus\fShat\oplus\fso(M,g)\oplus\fR
		= \Gamma(\eEhat)
		\cong C^\infty_{\Hhat}\qty(G_{\overline 0};\fshat)
	\end{equation}
	where the latter is the space of \(\Hhat\)-equivariant maps \(\chi:G_{\overline 0}\to \fshat\) (i.e. those maps satisfying \(\chi(gh)=h^{-1}\cdot\chi(g)=\varrho(h)^{-1}\chi(g)\) for \(g\in G_{\overline 0}\), \(h\in \Hhat\)).
	We define a map \(\widetilde\Psi:\fg\to C^\infty_{\Hhat}\qty(G_{\overline 0};\fshat)\) by
	\begin{equation}
	\begin{aligned}
		\widetilde\Psi(X)(g)
			&= \overline{\Ad_{g^{-1}}X} - \Phihat(\Ad_{g^{-1}}X)
			= \Ad_{g^{-1}}X - \lambda(\overline{\Ad_{g^{-1}}X}),	\\
		\widetilde\Psi(s)(g)
			&= \Ad_{g^{-1}}s,
	\end{aligned}
	\end{equation}
	for \(X\in\fg_{\overline 0}\), \(s\in\fg_{\overline 1}=S'\) and \(g\in G_{\overline 0}\), where \(\Ad\) is the adjoint representation of the super Harish-Chandra pair \((G_{\overline 0},\fg)\), and we write \(\overline{X}\) for the projection of \(X\in\fg_{\overline 0}\) to \(V\simeq \fg_{\overline 0}/\fhhat\).
	This is clearly linear and injective.
	
	Now let \(\Psi(X+s)\) be the section of \(\eEhat\) corresponding to \(\widetilde\Psi(X+s)\). 	
	We can identify the \(\fX(M)\)-component of \(\Psi(X)\) as the fundamental vector field \(\xi_X\) (indeed, the latter corresponds to the \(\Hhat\)-invariant map \(G_{\overline 0}\to V\), \(g\mapsto\Ad_{g^{-1}}X\)), and write \(\Psi(X+s) = (\xi_X,\epsilon_s,A_X,a_X)\) for some \(\epsilon_s\in\fShat\), \(A_X\in\fso(M,g)\) and \(a_X\in\fR\).
	
	We first show that the sections \(\Psi(X+s)\in\Gamma(\eEhat)\) corresponding to \(\widetilde\Psi(X+s)\) are parallel with respect to the connection \(\eDhat\) on \(\eEhat\) defined by equation~\eqref{eq:eDhat-def}.
	By homogeneity, it is sufficient to show this at the basepoint \(o=\Hhat\in G/M\).
	Identifying the fibre \(\eEhat_o\) with \(\fshat\), a computation using Wang's theorem (Lemma~\ref{lemma:Wang-spin-R}) and the definitions of the maps \(\Phihat\) and \(\widetilde\Psi\) and the brackets on \(\fg\) (which are opposite those in \eqref{eq:admiss-cocycle-brackets-rsymm}) gives us
	\begin{equation}
	\begin{aligned}
		\nablahat_{\xi_Y}\Psi(X)\eval_o 
			&= \varrho\qty(\Phihat(Y))\qty(X + \lambda(\overline{X})) + \comm{X}{Y}_{\fg} + \lambda\qty(\overline{\comm{X}{Y}_{\fg}})
			&&= \Phihat(X)\cdot\overline{Y} - \thetatilde(\overline{X},\overline{Y}),	\\
		\nablahat_{\xi_Y}\Psi(s)\eval_o 
			&= \varrho\qty(\Phihat(Y))\qty(s) - \comm{Y}{s}_{\fg}
			&&= \beta(\overline{Y},s).
	\end{aligned}
	\end{equation}
	and using the curvature equations~\eqref{eq:curvatures-homogeneous} from Lemma~\ref{lemma:homog-spinR-connection}, the calculations above show that 
	\begin{equation}
		\nablahat_{\xi_Y}\qty(\xi_X,\epsilon_s,A_X,a_X)
			= \qty(-A_X\xi_Y, \betabold(\xi_Y)\epsilon_s, R(\xi_X,\xi_Y), F(\xi_X,\xi_Y)),
	\end{equation}
	meaning \(\eDhat\Psi(X+s)=0\) since the values of the fundamental vector fields span all tangent spaces. In particular, \(\xi_X\) is a Killing vector field, \(A_X=A_{\xi_X}=-\nabla\xi_X\), \(D\epsilon_s=0\) and \(\nablahat a_X=\iota_{\xi_X}F\), and \(\Psi(X+s)(o)=\widetilde\Psi(X+s)(1_{G_{\overline 0}})=(\overline{X},s,-\Phihat(X))\) is the Killing transport data of \(\Psi(X+s)\) at \(o\).	
	Moreover, by \(G_{\overline 0}\)-invariance of \(\betabold\) and \(\rhobold\), we have \(\eLhat_{\xi_X}\betabold+a_X\cdot\betabold=0\) and \(\eLhat_{\xi_X}\rhobold+a_X\cdot\rhobold=0\), so we can now define 
	\begin{equation}
		\psi(X+s)=(\xi_X+a_X)+\epsilon_s \in \fKhat_{(D,\rhobold)}=\fVhat_{(D,\rhobold)}\oplus\fShat_D,
	\end{equation}
	giving the map \(\psi:\fg_{\overline 0}\to\fKhat_{(D,\rhobold)}\).	
	It can then be checked that \(\psi\) is an anti-homomorphism by comparing the Killing transport data \(\Psi(\comm{X}{Y})(o)=\widetilde\Psi(\comm{X}{Y})(1_{G_{\overline 0}})\) with the transport data of \(\comm{\xi_X+a_X}{\xi_Y+a_Y}\) (see e.g. equation~\eqref{eq:phihat-x-XY}), and similarly for the equations~\eqref{eq:reconstruct-brackets-rsymm}.
	 
	 Let us now consider the question of whether \((D,\rhobold)\) is an admissible pair. Conditions \ref{item:admiss-pair-killing} and \ref{item:admiss-pair-jacobi} of  Definition~\ref{def:killing-spinor-rsymm} are satisfied for all \(\epsilon,\zeta,\vartheta\in\fShat\) precisely because \(\beta+\gamma+\rho\) is Spencer cocycle for \(\fshat\), while the remaining condition \ref{item:admiss-pair-bg} is satisfied for \(\epsilon,\zeta\in\psi(\fg_{\overline 1})\subseteq\fShat_D\); if it in fact holds for all \(\epsilon,\zeta\in\fShat_D\) as hypothesised, the final claim is immediate.
\end{proof}

\begin{remark}
	It is worth commenting on the uniqueness of the connections \(\eA\) and \(D\) defined in Lemma~\ref{lemma:homog-spinR-connection} and in the theorem above. Once the realisable subdeformation \(\fg\) (or indeed its isomorphism class) is fixed, the admissible Spencer cohomology class \([\mu_-]\) is uniquely determined (Prop.~\ref{prop:highly-susy-filtered-def-rsymm}), but the representative cocycle \(\mu_-\) is only determined up to addition of \(\partial\chi\), where \(\chi\in\Hom(V,\fhhat)\), and the cocycle \(\beta+\gamma+\rho\in\cHhat^{2,2}\) is only determined up to addition of an element in \(\cKhat(\fa_-) = \{\beta+\gamma+\rho \in \cHhat^{2,2} \, |\,\beta|_{V\otimes S'}=0,\,\rho|_{\Odot^2S'}=0	\}\) by Lemma~\ref{lemma:maps-in-cohomology-rsymm}.
	The former indeterminacy would naïvely seem to affect the definition of the Nomizu map \(\Phihat\) and the corresponding connection \(\eA\), since a different choice of representative induces a change in the definition of the map \(\lambda\) appearing in equation~\eqref{eq:Nomizu-map-Phihat-def}. 
	However, implicit in the presentation of the bracket of the Lie algebra \(\fg_{\overline 0}\) as a deformation of that of \(\fa_{\overline 0}=\fhhat\ltimes V\subseteq \fshat_{\overline 0}\) is a choice of (linear) splitting of the short exact sequence of \(\fhhat\)-modules \(0\to\fhhat\to\fg_{\overline 0}\to V\to 0\), and a change of cocycle representative \(\mu_-\mapsto \mu_-+\partial\chi\) corresponds to a change of splitting such that the element of \(\fg_{\overline0}\) previously written as \(v+A\) should now be written as \(v+A-\chi(v)\), which compensates for the replacement \(\lambda\mapsto\lambda+\chi\) in equation~\eqref{eq:Nomizu-map-Phihat-def}.
	Moreover, it was already observed in the discussion at the beginning of \S\ref{sec:integrability-rsymm} that a change in the choice of representative cocycle does not change \(\thetatilde\), so the curvature is explicitly invariant.
	The indeterminacy by elements of \(\cKhat(\fa_-)\)  is more meaningful; it means that if \(\cKhat(\fa_-)\neq 0\) then there is not a unique pair \((D,\rhobold)\) for which \(\fg\) embeds into \(\fKhat_{(D,\rhobold)}\) -- two such pairs need only agree on \(\psi(\fg_{\overline 1})\). A similar observation was made for realisable filtered subdeformations of \(\fs\) in \cite[Rem.5.8]{Beckett2024_ksa}.
\end{remark}

\begin{remark}
	Finally, we observe that Theorem~\ref{thm:reconstruction-homog-spin-mfld-rsymm} could be stated more simply if we generalised our definition of the Killing superalgebra by defining the odd part to be
	\begin{equation}
		\qty{\epsilon\in\fShat \, \middle| \,
			D\epsilon = 0, \,
			\eLhat_{\kappabold(\epsilon,\zeta)}\betabold + \comm{\rhobold(\epsilon,\zeta)}{\betabold} = 0, \,
			\eLhat_{\kappabold(\epsilon,\zeta)}\rhobold + \rhobold(\epsilon,\zeta)\cdot\rhobold = 0,\,
			\imath_{\kappabold(\epsilon,\zeta)}F = \nablahat(\rhobold(\epsilon,\zeta))	}
	\end{equation}
	and omitting part \ref{item:admiss-pair-bg}, which amounts to saying that the above is equal to \(\fShat_D\), from the definition of admissible pairs. The theorem then says that \(\psi\) is always an embedding of \(\fg\) into a Killing superalgebra.
\end{remark}

\section*{Acknowledgements}

This work is dedicated to Paul de Medeiros, whose insight inspired much of it and the larger project it is part of. His unique perspective is sorely missed.

The author would like to thank José Figueroa-O'Farrill, under whose supervision the majority of this work was done, for his guidance and patience. Thanks also to Andrea Santi, James Lucietti and C.S. Shahbazi for many enlightening conversations and helpful comments. 
%Finally, thanks to the anonymous reviewers for \textit{---}, whose comments and suggestions greatly improved the quality and clarity of this work.

This work was in part carried out with scholarship funding from the Science and Technologies Facilities Council (STFC) and the School of Mathematics at the University of Edinburgh.

\subsection*{Relation to other works}

An earlier version of this work previously appeared as part of the author's PhD thesis \cite[Ch.4]{Beckett2024_thesis}, but it has not been published elsewhere. Substantial revisions and additions have been made; we note the following in particular.
\begin{enumerate}
	\item A discussion of technicalities related to the existence of bundle Dirac currents has been added as \S\ref{sec:dirac-current-bundle-rsymm}; Lemmas \ref{lemma:time-orientable-spin-R} and \ref{lemma:bundle-dirac-current-exist-rsymm} are entirely new. 
	\item Definitions have been revised and improved to allow for Killing (super)algebras whose even subalgebras are not sums of subalgebras of \(\fiso(M,g)\) and \(\fR\) on the geometric side, and for filtered deformations of graded subalgebras \(\fa\) of \(\fshat\) for which \(\fa_0=\fhhat\) is not a sum of subalgebras of \(\fso(V)\) and \(\fr\) on the algebraic side -- see Remarks~\ref{rem:better-KSA-def}, \ref{rem:better-KSA-def-2} and \ref{rem:better-KSA-def-3}.
	New intermediary results Lemma~\ref{lemma:hhat-ext} and Corollary~\ref{coro:annihilators} have been added, other results and proofs have been edited accordingly, and \S\ref{sec:towards-classification-rsymm} has been slightly expanded. Example~\ref{example:poincaré-group-R-ext-HC-pair} has been added, the statement of Theorem~\ref{thm:reconstruction-homog-spin-mfld-rsymm} has been made more precise, and details have been added to its proof.
	\item Lemma 4.26 of \cite[Ch.4]{Beckett2024_thesis}, which also appeared in \S\ref{sec:homogeneity-faithful-rsymm} in an earlier version of this standalone work, is false. It was claimed that if the \(R\)-symmetry group is compact, any subgroup \(R'\subseteq R\) preserving a subspace of the spinor module \(S'\subset S\) is of the form \(R'=R''\times K\), where \(K\) is the kernel of the action of \(R'\) on \(S'\) and \(R''\) acts faithfully on \(S'\). In a basis adapted to a split \(S=S'\oplus(S')^\perp\) induced by an invariant inner product, elements of \(R'\) can be written as block-diagonal matrices \(\diag(A,B)\), and it was claimed that \(R''\) and \(K\) could be identified as groups consisting of block-diagonal matrices of the form \(\diag(A,\1)\) and \(\diag(\1,B)\) respectively. But this fails, since there may be elements \(\diag(A,B)\in R'\) for which \(\diag(A,\1),\diag(\1,B)\notin R'\). The false claim has been removed and replaced by the new result Lemma~\ref{lemma:hhat-faithful-stronger}, which requires the additional hypothesis that \(\Spin(V)\cdot S'=S\).
	The results and proofs in \S\ref{sec:high-susy-spin-mfld-homogeneous-rsymm} have been updated accordingly, and this new hypothesis is also used in the expanded \S\ref{sec:towards-classification-rsymm}.
	A discussion of related technicalities has been added as Appendix~\ref{appendix}.
	\item Background on generalised spin structures appears in abridged form here (\S\ref{sec:spin-R}, \S\ref{sec:symm-alg-spin-r}, \S\ref{sec:homogeneous-lorentz-spin-mfld-rsymm}) but is adapted and expanded upon in \cite{Beckett2025_gen_spin}.
\end{enumerate}

\begin{appendices}

\section{Passing to a smaller Poincaré superalgebra}
\label{appendix}

Let \(\fa=V\oplus S'\oplus\fhhat\) be a highly supersymmetric graded subalgebra of an extended Poincaré superalgebra \(\fshat=V\oplus S\oplus(\fso(V)\oplus\fr)\) as in Sections~\ref{sec:filtered-def-flat-model-rsymm} and \ref{sec:high-susy-spin-mfld-rsymm}.
We wish to address the question of whether \(\fa\) can be viewed as being contained in a proper subalgebra \(\fsch=V\oplus\Sch\oplus(\fso(V)\oplus\frch)\) which is itself an extended Poincaré superalgebra.
Morally speaking, if \(\fshat\) is \(N\)-extended then \(\fsch\) will be \(M\)-extended for some \(M<N\).

We define \(\Sch:=\Spin(V)\cdot S'\subseteq S\), noting that this is a spinor representation which is preserved by the actions of \(\fhhat\) and \(\fr':=\pr_{\fr}(\fhhat)\) (though the latter action may not preserve \(S'\)). 
Indeed, if \(A+a\in\fhhat\) with \(A\in\fso(V)\) and \(a\in\fr\), then clearly the action of \(A\) preserves \(\Sch\), while for all \(g\in\Spin(V)\) and \(s\in S'\),
	\begin{equation}
		a (g\cdot s) = g\cdot (a s) = g\cdot((A+a)\cdot s) - (\Ad^{\Spin(V)}_g A)\cdot(g\cdot s) \in \Sch,
	\end{equation}
since \(A+a\in\fhhat\) preserves \(S'\) and \(\Ad^{\Spin(V)}_gA\in\fso(V)\) preserves \(\Sch\). We denote the restriction of the Dirac current by \(\widecheck\kappa:=\kappa|_{\Odot^2\Sch}\) and the corresponding \(R\)-symmetry group and algebra (which lie in \(\GL(\Sch)\) and \(\fgl(\Sch)\) respectively) by \(\Rch\) and \(\frch\). If \(\kappa\) is symmetric and causal then so is \(\widecheck\kappa\), and \(\Rch\) is compact if \(R\) is. Moreover, \(\fsch\) as written above is indeed a Poincaré superalgebra (i.e. it is an \(\frch\)-extended flat model superalgebra in the sense of Definition~\ref{def:flat-model-alg-rsymm} in Lorentzian signature).

Two issues must now be resolved. The first is to find an embedding \(\frch\hookrightarrow\fr\) which will allow us to consider \(\fsch\) as a subalgebra of \(\fshat\), and the second is to show that \(\fa\) is contained in this subalgebra. In fact, the latter is not strictly true in general, but we will construct an appropriately-behaved re-embedding of \(\fhhat\) which will allow us to assume it essentially without loss of generality.

For the first issue, consider the action of \(\fstab_{\fr}(\Sch)\) on \(\Sch\); this must preserve \(\widecheck\kappa\), so it induces a short exact sequence
\begin{equation}
\begin{tikzcd}
	0 \ar[r] &\fann_{\fr}(\Sch) \ar[r] &\fstab_{\fr}(\Sch) \ar[r] &\frch \ar[r] &0
\end{tikzcd}
\end{equation}
which we wish to split.
If \(R\) is compact then so is \(\Stab_R(\Sch)\), so since \(\fk:=\fann_{\fr}(\Sch)\) is an ideal, there exists an \(\Ad\)-invariant inner product on the stabiliser algebra which decomposes it into a sum of ideals \(\fstab_{\fr}(\Sch)=\fk\oplus\fk^\perp\); this gives us an isomorphism of Lie algebras \(\fk^\perp\to\frch\) whose inverse splits the short exact sequence, giving us an embedding \(\frch\hookrightarrow\fstab_{\fr}(\Sch)\subseteq\fr\).
Alternatively, one might be interested in ensuring that we also have a compatible embedding of groups \(\Rch\hookrightarrow R\). This amounts to splitting the short exact sequence of groups
\begin{equation}
\begin{tikzcd}
	1 \ar[r] & \Stab_R^{\mathrm{pw}}(\Sch) \ar[r] &\Stab_R(\Sch) \ar[r] &\Rch \ar[r] &1
\end{tikzcd}
\end{equation}
where \(K = \Stab_R^{\mathrm{pw}}(\Sch) = \{r\in R\,|\,\forall s\in\Sch,\,rs=s\}\) is the pointwise stabiliser of \(\Sch\), the Lie algebra of which is \(\fann_\fr(\Sch)\); this sequence is similarly constructed from the action of the stabiliser group on \(\Sch\).
It is not guaranteed that the splitting constructed above integrates to the group level, but one can usually show in explicit examples that a splitting of the sequence of groups does in fact exist. Where this is the case, we generally find that the splitting induces a semi-direct product structure \(\Stab_R(\Sch)=\Rch\ltimes K\) which differentiates to \(\fstab_{\fr}(\Sch)=\frch\ltimes\fk\) rather than a direct sum of ideals.
Thus, while there is not a canonical splitting in general, in practice one can usually be found nonetheless. It is easy to see that the resulting embedding \(\frch\hookrightarrow\fr\) gives us an embedding \(\fsch\hookrightarrow\fshat\) as desired.

\begin{example}
	In the simplest case, if the \(R\)-symmetry is \(R=\Orth(N)\) (i.e. the spinors are Majorana), we can write \(S=S_1\otimes\Delta_N\), where \(S_1\) is the irreducible spinor representation and \(\Delta_N\simeq\RR^N\) is the defining representation of \(R\) carrying an invariant inner product.
	We then have \(\Sch=S_1\otimes\Delta_M\), where \(\Delta_M\) is an \(M<N\)-dimensional subspace of \(\Delta_N\), and \(\Rch\cong\Orth(\Delta_M)\cong\Orth(M)\), and we obtain a splitting
	\begin{equation}
		\Stab_R(\Sch) = \Stab_R(\RR^M) \cong \Orth(\Delta_M)\times\Orth(\Delta_M^\perp) \cong \Orth(M)\times\Orth(N-M)
	\end{equation}
	by choosing an adapted orthonormal basis for \(\Delta_N=\Delta_M\oplus\Delta_M^\perp\), in which matrices in the stabiliser take a block-diagonal form and we have \(K=\Orth(\Delta_M^\perp)\).
\end{example}

For the second issue, consider the fact that even if \(\fhhat\) acts faithfully on \(S'\) and thus on \(\Sch\), \(\fr'=\pr_{\fr}(\fhhat)\) may not.
Indeed, there may be elements \(A+a\in\fhhat\) with \(0\neq A\in\fso(V)\) (so \(A+a\) does not annihilate \(S'\)) but \(0\neq a\in\fann_\fr(S')\)  (thus necessarily \(a\notin\fhhat\)).
Apart from being somewhat undesirable by itself, this behaviour ensures that \(\fr'\) cannot embed in \(\frch\), since the latter acts faithfully on \(\Sch\). We therefore seek a subalgebra of \(\fso(V)\oplus\fr\) whose projection to \(\fr\) is contained in the image of \(\frch\) but which is sufficiently closely related to the original \(\fhhat\) that it is a reasonable substitute.

\begin{lemma}
	There exists a Lie subalgebra \(\fhch\subseteq\fso(V)\oplus\frch\) and a surjective Lie algebra morphism \(\phi:\fhhat\to\fhch\) such that 
	\begin{enumerate}
	\item the action of \(\fhhat\) on \(\Sch\) factors through that of \(\fhch\), i.e. \(\phi(X)\cdot s = X\cdot s\) for all \(X\in\fhhat\) and \(s\in \Sch\),
	\item \(\frch' := \pr_{\frch}(\fhch)\) coincides with the image of \(\fr'=\pr_{\fr}(\fhhat)\) in \(\frch\),
	\item \(\pr_{\fso(V)}(\fhch)=\pr_{\fso(V)}(\fhhat)\),
	\item \(\fhhat\cap\fso(V)\subseteq\fhch\cap\fso(V)\) and \(\phi\) restricts to an inclusion \(\fhhat\cap\fso(V)\hookrightarrow\fhch\cap\fso(V)\).\label{item:hcheck-4}
	\end{enumerate}
	Furthermore, \(\phi\) is an isomorphism if and only if \(\fhhat\) acts faithfully on \(S'\).
\end{lemma}

\begin{proof}
	We recall that the action of \(\fr'\) on \(S\) preserves \(\Sch\), i.e. \(\fr'\subseteq\fstab_{\fr}(\Sch)\), so we have a map \(\fr'\to\frch\) which is nothing but the restriction to \(S'\), where we recall that \(\fr\subseteq\fgl(S)\) and \(\frch\subseteq\fgl(\Sch)\). Thus, let us define
	\begin{equation}
	\begin{aligned}
		\phi&: \fhhat \longrightarrow 
			\fhch = \qty{ A + a|_{\Sch} \, \middle| \, A\in\fso(V), a\in\fr \text{ s.t. } A+a\in\fhhat},
			&& A+a\mapsto A+a|_{\Sch}.
	\end{aligned}
	\end{equation}
	Clearly \(\phi\) is a surjective Lie algebra morphism with the properties claimed; part~\ref{item:hcheck-4} follows from the observation that if \(A+a\in\fhhat\) and \(a\in\fann_{r}(\Sch)\) then \(\phi(A+a)=A\).
	We have \(\phi(A+a)=0\) if and only if \(A=0\) and \(a\) annihilates \(\Sch\), or equivalently \(S'\), whence \(\ker\phi=\fann_{\fhhat\cap\fr}(S')\) which by Lemma~\ref{coro:annihilators} is zero if and only if \(\fhhat\) acts faithfully on \(S'\).
\end{proof}

Now treating \(\frch\) and \(\fsch\) as subalgebras of \(\fr\) and \(\fshat\) respectively (while recalling that this identification is dependent on the splitting of the stabiliser), we note that while \(\fhch\) is not equal to \(\fhhat\) as a subalgebra of \(\fshat_0\), its actions on \(\Sch\) and \(S'\) are essentially indistinguishable, so that if \(\fa\) is transitive, it is isomorphic to the subalgebra \(\widecheck\fa = V\oplus S'\oplus\fhch\) of \(\fsch\) (or more generally, \(\fa\twoheadrightarrow\widecheck\fa\), with the kernel being \(\fann_{\fhhat\cap\fr}(S')\)). Thus we may choose to work with \(\widecheck\fa\subseteq\fsch\) instead of with \(\fa\subseteq\fshat\). In particular, this allows us to assume without loss of generality that \(\Spin(V)\cdot S'=S\) as required in \S\ref{sec:towards-classification-rsymm} and \S\ref{sec:high-susy-spin-mfld-homogeneous-rsymm}.

\end{appendices}

\printbibliography[heading=bibintoc,title={Bibliography}]

\end{document}